\documentclass{amsbook}
\usepackage{amsmath,amssymb}

\newtheorem{theorem}{Theorem}[chapter]

\theoremstyle{definition}

\theoremstyle{remark}

\numberwithin{section}{chapter}
\numberwithin{equation}{chapter}

\def\={\ =\ }
\def\+{\ +\ }
\def\-{\ -\ }

\def\b|{\big |}

\def\g1{\Gamma_1}

\def\llangle{\langle\!\langle}\def\rrangle{\rangle\!\rangle}\def\plainboxed#1{#1}
\def\twoheadrightarrow{\longrightarrow}\def\xrightarrow#1{\mathop{\longrightarrow}\limits^{#1}}

\def\begintheorem#1{\par\medskip\noindent\textbf{#1. }\itshape}
\def\endtheorem{\par\medskip\normalfont}
\long\def\proclaim#1#2\endproclaim{\par\medskip\noindent\textbf{#1. }{\itshape #2\par}\medskip}
\def\roster{\begin{list}{}{%
  \setlength{\leftmargin}{0pt}%
  \setlength{\rightmargin}{0pt}%
  \setlength{\labelwidth}{0pt}%
  \setlength{\labelsep}{0.5em}%
  \setlength{\itemindent}{0pt}%
  \setlength{\listparindent}{\parindent}%
  \setlength{\itemsep}{0.25em}}}
\def\endroster{\end{list}}
\renewenvironment{quote}{\par\medskip\noindent\ignorespaces}{\par\medskip}
\def\itemitem#1{\item[#1]}
\def\nfp{\par\medskip\noindent\textit{Proof.}\ }

\makeatletter
\newcommand{\fullwidthtext}{%
  \par\leftskip=0pt\relax\rightskip=0pt\relax\parshape=0\relax
  \@totalleftmargin=0pt\relax\linewidth=\textwidth\columnwidth=\textwidth\hsize=\textwidth}
\makeatother

\begin{document}
\frontmatter
\title{Word Equations: Pairs and their Makanin--Razborov Diagrams}
\author{Z. Sela}
\address{Hebrew University, Jerusalem 91904, Israel}
\maketitle
\begin{abstract}
Sets of solutions to systems of equations (varieties) over a free group can be understood using their associated Makanin--Razborov diagram that encodes them. In this monograph we construct an analogous diagram that is associated with a set of solutions to systems of equations over a free semigroup.

The diagram that we construct encodes a variety over a free semigroup, i.e., the set of solutions to a system of equations over a free semigroup. It has similar geometric properties as the diagram over a free group, and can be further applied to study some fragments of the first order theory of a free semigroup. In a forthcoming work, we use the techniques that we developed for varieties over free semigroups to analyze varieties over a free (non-commutative) associative algebra.
\end{abstract}
\mainmatter
\tableofcontents

\chapter{Introduction}
\fullwidthtext

\medskip

A free semigroup can be viewed as the most basic formal language. This
connection, and the analogy with Tarski's problem on the first order
theories of  free groups, led W. Quine to study the first order theory of a free
semigroup. In 1946 Quine proved that arithmetic can be defined in a free semigroup.
Hence, by Godel's incompleteness theorems, the theory is not axiomatizable, nor
decidable [Qu]. Later, smaller and smaller fragments of the theory have been shown
to be undecidable, including sentences with only two quantifiers,
by Durnev ([Du1],[Du2]), Marchenkov [Mar], and others.

On the positive side, the Diophantine theory of a free semigroup has been
 shown to be decidable by G. S. Makanin [Ma1]. Makanin presented an algorithm
that decides if a given system of equations in a free semigroup has a solution.
Several years later Makanin was able to modify his algorithm to decide if
a given system of equations in a free group has a solution [Ma2].

In 1987, A. A. Razborov managed to push forward Makanin's techniques and using quite involved
and technical combinatorial
arguments gave
a combinatorial description of the set of solutions to a system of equations
(variety) over a free group [Ra].
This description was further developed by
O. Kharlampovich and A. Myasnikov [Kh-My], and
a more geometric approach
was given by the author in [Se1].

The search for a description of the set of solutions to a system of equations
over free semigroups has a long history (see [Di]). For equations with one variable,
and then for equations with 3 variables, the structure of the set of solutions
in free semigroup was achieved before the analogous structures in free groups (see [Khm]).
In these cases the structure that was found was purely combinatorial.

Since the mid 1980's and in particular after Razborov's thesis, there have been
quite a few attempts to study sets of solutions in a free semigroup, usually
for particular families of systems, that are often
either with small number of
variables, or that are of rather particular type, mostly connected to quadratic
equations over a free group (e.g. [Ma3],[Ma4],[Ma5],[Ly],[Di] and their references).

In 2013 an algorithm  to enumerate and encode the
set of solutions to a general system of equations in a free semigroup was found  (by A. Jez).
This algorithm that enumerates the solutions efficiently is based on variants of
the Ziv-Lempel algorithm from information theory ([Je],[DJP]).

In this manuscript we present a geometric approach to study varieties in a free
semigroup. We use the combinatorial techniques that were introduced by Makanin
in proving the decidability of the Diophantine theory [Ma1], and we combine
them with geometric techniques that were used in the construction of the
JSJ decomposition of finitely presented groups ([Se7],[Ri-Se]), and with techniques that
appear in the solution to Tarski's problem.

Unfortunately, even though our approach is based on the construction of the JSJ decomposition
for groups, we were not able to find an analogue of the JSJ for studying varieties
in free semigroups (we were able to get such an analogue in some special cases).
However, we are able to find an analogue of Razborov's work over free groups, and
associate what we call a Makanin-Razborov diagram with each variety in a free
semigroup.

The MR diagram that we construct is not canonical, but it encodes all the points in
a given variety in a free semigroup. Furthermore, given a path in the diagram,
that we call a $resolution$, there exists a sequence of points in the associated variety
that factor through the resolution, and such that the sequence converges in the Gromov-Hausdorff
topology to an object from which the resolution can be reconstructed. Such sequences
are  viewed as generic points in the variety (a replacement for test sequences
that exist in a free group that are defined and used in [Se2]), and can be applied
to prove model theoretic results, that
include generalizations of Merzlyakov's theorem (over varieties) for free semigroups (see [Me] and [Se2]).

More recently, the construction of the MR diagram in free semigroups, and the generic sequences that factor
through its resolutions, were a main tool in the construction of a higher rank JSJ decomposition for studying
automorphisms of hierarchically hyperbolic groups (HHG) [Se8]. We believe that the combination of the techniques
from this recent work will enable one to make our constructions over free semigroups canonical, and in particular
will finally construct a (canonical) JSJ decomposition for semigroups (that will include finitely many canonical decompositions and
not just a single one as in groups).

The MR diagram in free  semigroups, its properties and the way it is constructed, suggest that
very basic questions about words (in formal languages)
are connected to concepts, objects and tools
from low dimensional topology. These include the JSJ decomposition, the
geometry and topology of surfaces, the dynamics of automorphisms of
surfaces and of free groups, and finally the dynamics of actions of groups
on real trees. We expect that this combination of techniques and structure
can be modified to describe sets of solutions to systems of equations in
free objects in other (associative and non-commutative) algebraic categories,
and we plan to continue in these directions.

Indeed, we started an ongoing  work on varieties
in  non-commutative, associative free algebras, that is partly based on our work on varieties in free semigroups (see [Se9],[At-Se]).
In analyzing varieties in these free algebras, our starting point is the analysis of sets of homogeneous solutions to homogeneous systems
of equations, and the techniques that we develop in this manuscript can be applied to perform this analysis and obtain a Makanin-Razborov
diagram that encodes these homogeneous solutions.

\medskip
Given a system of equations $\Sigma$ in  a free semigroup:
\begin{align*}
u_1(x_1,\ldots,x_n) &= v_1(x_1,\ldots,x_n) \\
& \vdots \hfill \\
u_s(x_1,\ldots,x_n) &= v_s(x_1,\ldots,x_n) \\
\end{align*}
we naturally associate with it a f.p.\ semigroup:
$$S(\Sigma) \ = \ \langle \, x_1,\ldots,x_n \, | \, u_1=v_1,\ldots,u_s=v_s \, \rangle $$

The set of solutions of $\Sigma$ over a free semigroup, $FS_k=<a_1,\ldots,a_k>$, is in bijection
with the set of semigroup homomorphisms: $\{h:S(\Sigma) \to FS_k\}$. Hence, studying the variety
of solutions to $\Sigma$ is equivalent to studying the structure of the set of homomorphisms from
$S(\Sigma)$ to the free semigroup $FS_k$.

With a f.g.\ semigroup $S$ we associate a group $Gr(S)$, that is obtained by forcing all the elements
in $S$ to have inverses, or alternatively, by looking at a presentation of $S$ as a semigroup as if
it is a presentation of a group. Naturally, $S$ is mapped into $Gr(S)$, but in general it is not embedded
in it. Let $\hat S$ be the (semigroup) image of $S$ in $Gr(S)$.

The free semigroup $FS_k$ embeds in the free group $F_k$ in a standard way.
A simple observation shows that every semigroup homomorphism from $S$ to $FS_k$ extends uniquely to a group
homomorphism from $Gr(S)$ to $F_k$ (see section 3). Hence, one can replace the set of semigroup
homomorphisms from $S$ to the free semigroup $FS_k$, with the set of pair homomorphisms:
$\{\eta:(\hat S,Gr(S)) \to (FS_k,F_k)\}$. i.e., those group homomorphisms from $Gr(S)$ to $F_k$ that map
$\hat S$ into the standard free semigroup $FS_k$.

Our whole approach is based on studying the structure of these pair homomorphisms. An immediate application
of the techniques that were used in free groups, shows that the set of pair homomorphisms
from $(\hat S,Gr(S))$ to $(FS_k,F_k)$ is canonically a finite union of the sets of pair homomorphisms from
pairs of the form $(S_i,L_i)$ to $(FS_k,F_k)$, where for each $i$,  $L_i$ is a  limit group, $S_i$ is a  subsemigroup that
generates $L_i$ as a group, and $(S_i,L_i)$ is a (limit)
quotient of the pair $(\hat S,Gr(S))$.

A free semigroup, and the  semigroups that we need to consider, have usually very few automorphisms.
The ability to replace a semigroup with a pair in which the ambient
group is a limit group, enables one to
work with a large group of automorphisms (of the limit group), that usually do not preserve the embedded semigroup.
In section 4 we describe an analogue of the shortening procedure for homomorphisms of pairs. Technically,
the shortening procedure for homomorphisms of pairs is much more involved than its analogue for
homomorphisms of groups.

In section 5 we describe a construction of a JSJ decomposition for a pair $(S,L)$, where $L$ is a limit group,
and $S$ generates it as a group. Unfortunately, the construction  applies only in special cases. Note that free products
exist in the categories of groups, of semigroups, and of pairs. However, it may be that a pair $(S,L)$, where
$L$ is a limit group,  is freely indecomposable as a pair, but the limit group $L$ is freely decomposable.
This simple fact implies that any attempt to borrow concepts from the JSJ
theory and from the construction of the MR diagram in free groups
must be further refined.
In section 8 we describe  the construction of the MR diagram for pairs
in the freely indecomposable case (theorem 8.8).
In section 9 we finally describe the construction of the MR diagram for
pairs in the general case (theorem 9.10).

As in free (and hyperbolic) groups, there are pairs, $(S_i,L_i)$,
associated with the nodes of the MR diagram,
together with their associated decompositions and modular groups. However, unlike the case of groups, the abelian decompositions
that are associated with the pairs that appear in the nodes of the diagram, need to recall not only the
algebraic structure of the group in question, but rather dynamical properties and the associated modular groups.
Hence, we need to extend the classes of vertex groups that we borrow from the JSJ theory of groups. In particular,
the abelian decompositions that appear along the MR diagram for pairs contain a new type of vertex groups that we call $Levitt$.
These Levitt vertex groups are free factors of the ambient limit groups that are connected to other vertex groups
by edges with trivial stabilizers, and each such vertex group contributes the group of its automorphisms (that are automorphisms
of free groups that are not
necessarily geometric automorphisms) to the modular group of the pair with which it is associated.

In the MR diagram in free groups, when we follow a path (resolution) in the diagram, the groups that are
associated with the nodes  along the path, form a  finite
sequence of successive proper quotients (maximal shortening quotients as they appear in [Se1]). In the MR diagram in
free semigroups that we construct in this manuscript this is not true. When we go over a resolution  in the diagram, the pairs (and groups)  that appear along
the nodes of the diagram are quotients, but not necessarily proper quotients.

The resolutions in the MR diagram in free groups end with  free groups
of various ranks. In the MR diagram in free semigroups, resolutions
end with marked graphs with directed edges, that have free groups as their fundamental groups. Group homomorphisms from a  given
free group into the coefficient free group are easily understood using substitutions, and they are in bijection with
a product set of the coefficient  group. In the case of semigroups, the set of homomorphisms that are associated with a terminal
node are also obtained  by substitutions, and these sets can be viewed as natural projections of affine spaces, i.e., natural
projections of product sets of the coefficient semigroup.

We added an introductory section (section 2) for the benefit of non-experts. In this section there is a survey of the basic notions
and objects that are used in studying varieties in free groups, semigroups and associative algebras, and the analogies between these objects and
basic notions in commutative algebraic geometry. In this section we also summarize the state of the art in the study of varieties in these basic
free non-commutative algebraic categories.

\smallskip
Finally I would like to thank Eliyahu Rips who encouraged me to proceed along  this long project, that somehow involved repeated
inaccuracies and misconceptions. I also thank the referees and the editor whose remarks made the presentation of this manuscript considerably
better.

\bigskip

\vglue 1.5pc
\chapter{Equations and the Structure of Their Solution Sets}
\fullwidthtext

\medskip

\section{Introduction}
\fullwidthtext
One of the oldest problems in mathematics is to solve a system of equations.  The meaning of ``solve'', however, depends strongly on the ambient algebraic category.  Given a system $\Sigma(X)=1$ over a group $G$, one may ask whether a solution exists, ask for an algorithm producing one, or ask the substantially deeper question:
\quote
{\it What is the structure of the entire set of solutions?\/}
\endquote
The last question is the subject of this chapter.

The model is classical algebraic geometry.  For polynomials $f_1,\dots,f_m\in K[x_1,\dots,x_n]$, one studies not merely one common zero but the algebraic set
$$
V(f_1,\dots,f_m)=\{a\in K^n:f_i(a)=0\hbox{ for all }i\},
$$
its irreducible components, dimension, coordinate ring, and maps to other varieties.  A similar transition - from individual solutions to the geometry of the complete solution set - occurred in the theory of equations over free groups.

For free groups the principal structural objects are limit groups, JSJ decompositions, and Makanin-Razborov (MR) diagrams.  Very roughly, limit groups arise as algebraic limits of families of solutions; JSJ decompositions describe their essential splittings and internal flexibility; and MR diagrams organize all possible families of solutions into a finite hierarchy.

We compare four settings: commutative rings, free semigroups and monoids, free groups, and free associative algebras.  The common theme is
$$
\plainboxed{\hbox{equations}\longrightarrow\hbox{coordinate object}\longrightarrow
\hbox{homomorphisms}\longrightarrow\hbox{structured families of solutions}.}
$$
The mechanisms behind the final arrow are quite different in the four categories.

\section{From equations to homomorphisms}
\fullwidthtext
\subsection{The commutative prototype}
\fullwidthtext
Let $K$ be a field and $I=(f_1,\dots,f_m)\subseteq K[x_1,\dots,x_n]$.  Every solution $a=(a_1,\dots,a_n)$ determines an evaluation homomorphism
$$
K[x_1,\dots,x_n]/I\longrightarrow K,\qquad x_i\longmapsto a_i.
$$
Over an algebraically closed field one replaces $I$ by the ideal $I(V)$ of all polynomials vanishing on the solution set.  Hilbert's Nullstellensatz identifies $I(V)$ with the radical of $I$.  Thus points can be viewed as homomorphisms from a coordinate algebra.

\smallskip
\subsection{Coordinate groups}
\fullwidthtext
Let $G$ be a group and $X=\{x_1,\dots,x_n\}$.  Equations with coefficients in $G$ are words in $G*F(X)$.  If $\Sigma=\{w_1,\dots,w_m\}$, put
$$
G_\Sigma=\frac{G*F(X)}{\llangle w_1,\dots,w_m\rrangle}.
$$
A solution in $G$ is exactly a homomorphism $G_\Sigma\to G$ restricting to the identity on the coefficient group $G$.  More intrinsically, define
$$
{ R}ad_G(\Sigma)=\{w\in G*F(X):w(\bar g)=1\hbox{ for every }\bar g\in Sol_G(\Sigma)\}
$$
and the coordinate group
$$
\Gamma_G(\Sigma)=\frac{G*F(X)}{{R}ad_G(\Sigma)}.
$$
Then
$$
Sol_G(\Sigma)\simeq Hom_G(\Gamma_G(\Sigma),G).
$$
This is the basic translation that makes a geometry of equations over groups possible.

\section{Several meanings of ``solve''}
\fullwidthtext
For a finite system $\Sigma$ one should distinguish:
\enumerate
\item the {\it decision problem\/}: is $Sol(\Sigma)$ nonempty?
\item the {\it search problem\/}: produce one solution;
\item the {\it enumeration problem\/}: list solutions effectively;
\item the {\it structural problem\/}: construct a finite or controlled object describing all solutions simultaneously.
\endenumerate
The fourth problem is mathematically the richest.  An enumeration may contain no visible explanation of why the solutions fall into particular families, while a structural description attempts to expose precisely that organization.

\section{The four categories at a glance}
\fullwidthtext
\subsection{Commutative rings}
\fullwidthtext
Over an algebraically closed field, radical ideals, prime ideals, irreducible decomposition, dimension theory, and coordinate rings describe polynomial solution sets.  Gr\"obner bases make many operations effective: ideal membership, elimination, dimension computations, and the solution of zero-dimensional systems.  Noetherianity is the central finiteness mechanism.

Over arithmetic rings the situation can be radically different.  The Davis-Putnam-Robinson-Matiyasevich theorem implies that no algorithm decides whether an arbitrary polynomial equation with integer coefficients has an integral solution.  Thus ``equations over rings'' does not have one uniform algorithmic answer.

\subsection{Free semigroups and monoids}
\fullwidthtext
Let $A$ be a finite alphabet.  A word equation is
$$
U(X,A)=V(X,A).
$$
A solution substitutes words for variables.  Makanin proved decidability of solvability of word equations.  Later work, including Plandowski's methods and recompression, gave more efficient and conceptually cleaner algorithms.  Complete solution sets admit finite descriptions; in particular, modern results represent them by EDT0L languages.

\subsection{Free groups}
\fullwidthtext
Let $F=F(a_1,\dots,a_r)$ be a nonabelian finitely generated free group.  Makanin proved decidability of equations over $F$.  Razborov went substantially further and constructed finite diagrams describing all solutions.  The geometric development of this theory uses limit groups, actions on {\bf R}-trees, JSJ decompositions, modular automorphisms, shortening arguments, and MR diagrams.

\subsection{Free associative algebras}
\fullwidthtext
Let $A=K\langle a_1,\dots,a_r\rangle$.  Here the unknowns themselves are noncommutative polynomials.  Gr\"obner-Shirshov bases provide important computational tools, but the unrestricted Diophantine problem for free associative algebras is undecidable.  Recent work develops structural MR/JSJ-type theories for significant classes of equations.  This is a developing research program and should not be presented as a theory of the same maturity as the free-group case.

\section{Historical route}
\fullwidthtext
The structural study of equations over free groups grew from classical decision problems in combinatorial group theory.  Makanin's work in the late 1970s and early 1980s established decidability first for free semigroups and then for free groups.  Razborov's decisive advance was to replace the yes/no answer by a finite description of the entire solution set.

In the 1990s and 2000s, work surrounding Tarski's problems transformed this into a form of noncommutative algebraic geometry.  Sela developed a geometric theory based on limit groups, {\bf R}-tree actions, JSJ decompositions, shortening, completions, and MR diagrams.  Kharlampovich and Myasnikov developed an algebraic-geometric and elimination-theoretic approach using fully residually free groups and fundamental sequences.  These theories led to deep results on the first-order theory of nonabelian free groups.

\section{Limit groups}
\fullwidthtext
\begintheorem{Definition}
A finitely generated group $L$ is {\it fully residually free\/} if for every finite set $S\subseteq L\setminus\{1\}$ there exists a homomorphism $h:L\to F$ that is injective on $S$.  Finitely generated fully residually free groups are called {\it limit groups\/}.
\endtheorem

They also arise from sequences of homomorphisms.  Let $h_n:G\to F$.  After passing to a suitable subsequence, define the stable kernel
$$
Ker_{st}(h_n)=\{g\in G:h_n(g)=1\hbox{ for all sufficiently large }n\}.
$$
The quotient $G/ Ker_{st}(h_n)$ is the algebraic limit of the sequence; in the situations central to the theory it is a limit group.

Geometrically, each $h_n$ gives an action of $G$ on the Cayley tree $T_F$.  If the lengths of the maps diverge, rescaling $T_F$ and taking a limit produces an action on an {\bf R}-tree.  The limiting action records asymptotic cancellation and growth.  Structure theory for such actions forces splittings of the acting group.  This is the bridge
$$
\hbox{families of solutions}\longrightarrow {R}\hbox{-tree actions}\longrightarrow\hbox{splittings}.
$$

\section{Bass-Serre theory and JSJ decompositions}
\fullwidthtext
Bass-Serre theory translates splittings into tree actions.  An amalgam
$$
G=A*_C B
$$
acts on a tree whose vertex stabilizers are conjugates of $A$ and $B$ and whose edge stabilizers are conjugates of $C$; conversely, a suitable action on a tree yields a graph-of-groups decomposition.

A group can admit many inequivalent splittings.  A JSJ decomposition organizes, in a canonical or deformation-invariant manner, the essential splittings over a prescribed class of edge groups.  For limit groups the relevant edge groups are usually cyclic or abelian, depending on the version of the theory.

Typical vertices are:
\itemize
\item {\it abelian vertices\/};
\item {\it surface-type\/} or quadratically hanging (QH) vertices;
\item {\it rigid vertices\/}, which admit no further relevant compatible splitting.
\enditemize
Surface vertices are fundamental groups of compact surfaces (with the incident edge groups represented by boundary subgroups, in the appropriate relative sense).  Their mapping class groups produce large groups of automorphisms.

\section{Modular automorphisms and shortening}
\fullwidthtext
A JSJ decomposition supplies natural automorphisms.  Splittings give Dehn twists; QH vertices give mapping-class-group automorphisms; abelian vertices contribute further automorphisms.  These generate a modular group $Mod(L)$, with precise definitions varying slightly among versions of JSJ theory.

Fix generating sets for $L$ and $F$ and define
$$
\|h\|=\max_{s\in S_L}|h(s)|_F.
$$
Within the equivalence class obtained by precomposition with modular automorphisms (and, when appropriate, conjugation in the target), choose a shortest representative.  If $h_n:L\to F$ is an infinite sequence of essentially distinct short maps, a limiting {\bf R}-tree action is obtained.  If the limit exhibits flexibility already encoded by the JSJ decomposition, a modular automorphism would shorten the maps, contradicting shortness.  In the relevant cases, the remaining genuine degeneration therefore factors through a proper quotient
$$
L\twoheadrightarrow L'.
$$
This is the basic shortening mechanism behind the construction of MR diagrams.

\section{Makanin-Razborov diagrams}
\fullwidthtext
An MR diagram is a finite directed diagram encoding all homomorphisms from a coordinate group to the target free group.  Schematically,
\[
\begin{array}{ccccc}
&&L&&\\
&\swarrow&\downarrow&\searrow&\\
L_1&&L_2&&L_3\\
\downarrow&&\swarrow\searrow&&\downarrow\\
\cdots&&\cdots&&\cdots
\end{array}
\]
with edges represented by quotient maps.  A branch
$$
L=L_0\twoheadrightarrow L_1\twoheadrightarrow\cdots\twoheadrightarrow L_k
$$
represents not one homomorphism but a family.  At intermediate levels one may precompose with modular automorphisms:
$$
L_0\xrightarrow{\alpha_0}L_0\xrightarrow{\pi_1}L_1
\xrightarrow{\alpha_1}L_1\longrightarrow\cdots\longrightarrow L_k\xrightarrow{\phi}F,
\qquad \alpha_i\in Mod(L_i).
$$
The finiteness theorem is highly nontrivial: suitable complexity decreases along proper shortening quotients, and descending chains of relevant limit quotients cannot continue indefinitely.

\section{JSJ versus MR}
\fullwidthtext
The two constructions have different roles.  A JSJ decomposition is attached to one group and describes its splittings and natural symmetries.  An MR diagram is attached to a family of homomorphisms and organizes the different quotient types through which those homomorphisms factor.  The conceptual chain is
$$
\plainboxed{
\hbox{JSJ}\Rightarrow\hbox{modular automorphisms}\Rightarrow\hbox{shortening}
\Rightarrow\hbox{proper limit quotients}\Rightarrow\hbox{MR diagram}.}
$$
A useful slogan is that JSJ theory describes {\it horizontal variation\/} within a family, whereas the MR diagram records {\it vertical degeneration\/} among families.

\section{A worked surface-type example}
\fullwidthtext

We now replace the elementary abelian example by a genuinely surface-type equation.  The purpose is to display, concretely, the role of a QH vertex and of its modular automorphisms.  The example is deliberately chosen so that the geometry can be seen directly.

Let
$$
F=F(a,b)
$$
and put
$$
c=[a,b]=aba^{-1}b^{-1}.
$$
Consider the quadratic equation over $F$
\equation
[x,y]=c.
\endequation
The left side is the boundary word of an orientable surface of genus one with one boundary component.  Thus (1) is naturally a once-punctured-torus equation.

\subsection{The coordinate group as a surface attachment}
\fullwidthtext
Before taking the radical, the group presented by the equation is
$$
G=\langle F,x,y\mid [x,y]=c\rangle.
$$
Equivalently,
\equation
G\cong F(a,b)*_{\langle c\rangle=\langle [x,y]\rangle}F(x,y).
\endequation
The second vertex group $F(x,y)$ is the fundamental group of a once-punctured torus, and its boundary subgroup is $\langle[x,y]\rangle$.  The first vertex is the coefficient group.  The graph of groups is therefore
\medskip
\centerline{\hbox{$F(a,b)\quad -_{\langle c\rangle=\langle[x,y]\rangle}\quad \pi_1(\Sigma_{1,1})=F(x,y)$}}
\medskip
This is exactly the kind of surface (QH) vertex that occurs in JSJ theory.  The boundary subgroup is distinguished because it is glued to the coefficient group.

A solution of (1) is an $F$-homomorphism
$$
h:G\to F,
$$
meaning that $h$ is the identity on the coefficient vertex and that $h(x),h(y)$ satisfy $[h(x),h(y)]=c$.

\subsection{A base solution}
\fullwidthtext
There is an obvious solution
$$
h_0(x)=a,\qquad h_0(y)=b,
$$
because $[a,b]=c$.  The important point is that this one solution sits inside an infinite structured family generated by surface automorphisms.

\subsection{Dehn twists on the punctured torus}
\fullwidthtext
The mapping class group of the once-punctured torus acts on $F(x,y)$ by outer automorphisms preserving the boundary conjugacy class.  Two elementary Nielsen transformations already exhibit the phenomenon.  Define
$$
\alpha_n:x\mapsto xy^n,\quad y\mapsto y,
$$
and
$$
\beta_n:x\mapsto x,\quad y\mapsto yx^n,
\qquad n\in{Z}.
$$
With the convention $[u,v]=uvu^{-1}v^{-1}$, one checks directly that
$$
[xy^n,y]=[x,y],
\qquad
[x,yx^n]=[x,y].
$$
Hence both automorphisms fix the boundary word $[x,y]$ exactly and therefore extend to automorphisms of $G$ fixing the coefficient group $F$ pointwise.  Topologically these are Dehn twists about the two standard essential curves on the punctured torus.

Precomposing $h_0$ with these automorphisms produces explicit infinite families of solutions.  For example,
$$
h_0\circ\alpha_n:
\qquad x\mapsto ab^n,\qquad y\mapsto b,
$$
and indeed
$$
[ab^n,b]=[a,b]=c.
$$
Similarly,
$$
h_0\circ\beta_n:
\qquad x\mapsto a,\qquad y\mapsto ba^n
$$
is a solution for every $n\in{Z}$.

More generally, words in the transformations $\alpha_n$ and $\beta_n$ generate a large modular orbit of $h_0$.  Thus a single point of the solution set expands into an infinite family for a geometric reason: the surface vertex has a nontrivial mapping class group.

\subsection{The surface branch of the MR picture}
\fullwidthtext
For this family, the relevant branch of the MR description should not be pictured as immediately collapsing the surface group to an abelian quotient.  Instead the surface vertex survives and its modular group carries the parameters.  Schematically:
\medskip
\begin{center}
\begin{tabular}{c}
$G=F*_{\langle c\rangle}\pi_1(\Sigma_{1,1})$\\
$\downarrow$\\
$\alpha\in\mathop{Mod}_F(G)$ (surface mapping classes / Dehn twists)\\
$\downarrow$\\
$G$\\
$\downarrow^{h_0}$\\
$F$
\end{tabular}
\end{center}
\medskip
The displayed branch represents the maps
$$
\{h_0\circ\alpha:\alpha\in Mod_F(G)\}.
$$
This is already an infinite family, although the underlying diagram is finite.  It illustrates one of the central points of MR theory: parameters need not appear as elements of the diagram; they can live in the modular automorphism groups attached to its vertices.

\subsection{Where shortening enters}
\fullwidthtext
The preceding modular orbit is only one family of solutions of the quadratic equation; a full canonical MR diagram must account for every $F$--homomorphism $G\to F$.  We do {\it not\/} claim that the one branch drawn above is the complete canonical MR diagram of (1).  That stronger statement would require the general shortening and limit-group machinery.

Suppose, therefore, that $h_i:G\to F$ is an infinite sequence of solutions of (1).  The surface JSJ tells us first to use mapping classes of $\Sigma_{1,1}$ to choose short representatives of the modular equivalence classes of the $h_i$.  If infinitely many essentially different short representatives remain and their lengths diverge, rescaling the Cayley tree of $F$ gives a limiting {\bf R}-tree action.  The shortening argument separates surface flexibility, which has already been absorbed by $Mod_F(G)$, from genuine degeneration.  A genuine degeneration produces a proper limit quotient
$$
G\twoheadrightarrow L.
$$
The canonical MR construction adds the finitely many relevant shortening quotients as lower vertices and repeats the procedure there.

Thus the worked example exhibits both parts of the general mechanism:
\[
\begin{gathered}
\text{surface equation }[x,y]=c\\
\downarrow\\
\text{QH vertex }\pi_1(\Sigma_{1,1})\\
\downarrow\\
\text{mapping-class / Dehn-twist automorphisms}\\
\downarrow\\
\text{infinite modular families of solutions}\\
\downarrow\\
\text{shorten modulo these symmetries}\\
\downarrow\\
\text{remaining degenerations factor through proper limit quotients.}
\end{gathered}
\]
This is more representative of the role of JSJ theory in an MR diagram than an abelian centralizer example: the topology of a surface directly generates the internal parameters of a family of solutions.

\begintheorem{Remark}
There is a subtle but important distinction between the displayed group $G=\langle F,x,y\mid[x,y]=c\rangle$ and the coordinate group defined using the full radical of the solution set.  In a general discussion of equations one should use the radical coordinate group.  For the worked example, the presentation above is the natural relative surface group through which the displayed solutions factor and is the correct object for illustrating the QH/JSJ mechanism.  Establishing all properties of the radical coordinate group, or writing the full canonical MR diagram, belongs to the deeper theory rather than to this introductory calculation.
\endtheorem

\section{Quadratic equations and surfaces more generally}
\fullwidthtext
Equations of the form
$$
[x_1,y_1]\cdots[x_g,y_g]=w
$$
are orientable quadratic equations and are naturally associated with surfaces of genus $g$ with boundary.  Nonorientable quadratic equations involve products of squares.  Surface mapping class groups act on the corresponding surface groups, and their automorphisms generate large families of solutions.  This is why QH vertices are fundamental in JSJ decompositions of limit groups: the topology is not decorative but supplies actual algebraic parameters.

\section{Comparison with irreducible decomposition}
\fullwidthtext
There is a useful, but imperfect, analogy with algebraic geometry.  If
$$
\sqrt I=\mathfrak p_1\cap\cdots\cap\mathfrak p_r,
$$
then
$$
V(I)=V(\mathfrak p_1)\cup\cdots\cup V(\mathfrak p_r).
$$
Over free groups, irreducible algebraic sets are closely related to limit-group coordinate objects, and maximal limit quotients play a role analogous to irreducible components.  But an MR diagram goes further: after passing to one limit quotient, its homomorphisms may themselves degenerate in several ways.  The diagram is therefore hierarchical and is better thought of as combining decomposition, resolution, and parametrization.

\section{Algorithms over free groups}
\fullwidthtext
Makanin's theorem decides whether a finite system over a free group has a solution.  Once existence is decidable, a solution can in principle be found by effective search, although modern procedures are more direct.  Razborov's theory addresses the stronger structural problem: it gives a finite diagram encoding all solutions.  Sela's canonical MR construction and related developments organize this description geometrically through limit groups and shortening quotients.

A second, complementary modern description is language-theoretic.  If elements of $F$ are represented by freely reduced words and a tuple is encoded as $w_1\#\cdots\#w_n$, then the language of all reduced-word solutions is effectively EDT0L.  This is a finite generative description of the same underlying set of tuples, but it emphasizes different structure.  EDT0L theory sees words and finite-state control; MR theory sees homomorphisms, splittings, modular symmetries, and degenerations.

\section{Word equations in free semigroups}
\fullwidthtext
For $U(X)=V(X)$ in a free monoid, variables may overlap and exhibit periodicity even though there are no inverses.  Makanin proved decidability, and later methods substantially improved the algorithms.  Recompression transforms equations by repeatedly replacing pairs and blocks while preserving a controlled correspondence between solution sets.  Modern constructions yield finite graphs and EDT0L descriptions of all solutions.  In this book we develop an MR-type approach whose aim is to organize solution families by limiting and quotient structures, suggesting that hierarchical parametrization is not peculiar to group inverses.

\section{Polynomial equations and Gr\"obner bases}
\fullwidthtext
For $I=(f_1,\dots,f_m)\subseteq K[x_1,\dots,x_n]$, a Gr\"obner basis replaces the given generators by a finite set whose leading terms control the entire leading-term ideal.  Buchberger's algorithm computes such bases.  Elimination orders expose projections, while primary and prime decompositions expose components.  In this setting Noetherian ideal theory supplies the finiteness that, in free-group geometry, is supplied by the rigidity of limit groups, JSJ structure, and descending complexity of quotient resolutions.

\section{Associative algebras and the noncommutative frontier}
\fullwidthtext
Let $A=K\langle a_1,\dots,a_r\rangle$ and consider $P_i(X)=0$.  One may again form a coordinate algebra and view solutions as homomorphisms back to $A$.  Gr\"obner--Shirshov theory provides rewriting and ideal-theoretic tools, but there is no general Makanin-style algorithm because unrestricted Diophantine solvability is undecidable.

This does not preclude structural theories for restricted classes.  Recent work develops noncommutative algebraic geometry in which limit objects, JSJ-type decompositions, and MR-type diagrams are constructed for selected systems.  At present this should be viewed as an active frontier rather than a completed analogue of free-group Diophantine geometry.

\section{A unified comparison}
\fullwidthtext
\medskip
\centerline{\vbox{\halign{\strut #\hfil & \quad \vtop{\hsize=4.15in\noindent #\strut}\cr
{\bf Ambient object}&{\bf Coordinate object / structure / algorithmic tools}\cr
\noalign{\smallskip\hrule\smallskip}
Algebraically closed field&$K[X]/I(V)$; components and dimension; Gr\"obner bases and elimination.\cr
\noalign{\smallskip}
Free monoid&Word-equation data; finite families and EDT0L descriptions; Makanin's algorithm and recompression.\cr
\noalign{\smallskip}
Free group&$\Gamma_F(\Sigma)$; limit groups, JSJ decompositions, and MR diagrams; Makanin--Razborov methods and shortening.\cr
\noalign{\smallskip}
Free associative algebra&Coordinate algebra; developing MR/JSJ theory; Gr\"obner--Shirshov methods.\cr}}}
\medskip

The recurring pattern is
$$
\hbox{equations}\to\hbox{coordinate objects}\to\hbox{homomorphisms}\to\hbox{finite structural description},
$$
but the finiteness mechanism changes from category to category.

\section{From equations to logic}
\fullwidthtext
An existential formula has the form
$$
\exists X\;\Sigma(X,Y)=1.
$$
When $Y$ is a parameter, one is studying a projection of a solution set.  Alternating quantifiers, for example
$$
\forall y\,\exists x\;\Sigma(x,y)=1,
$$
require understanding how families of solutions vary with parameters.  Graded MR diagrams, completions, closures, formal solutions, and rigid and solid limit groups were developed to control such phenomena.  This leads from equations to definable sets and ultimately to the first-order theory of free groups.

\section{State of the art}
\fullwidthtext
For nonabelian finitely generated free groups, the theory is mature: solvability is decidable; solution sets admit finite MR descriptions; limit groups and JSJ decompositions organize their geometry; and reduced-word solutions admit effective EDT0L descriptions.  The structural theory extends, with important modifications, to torsion-free hyperbolic groups and several related classes.

For free semigroups and monoids, existence is decidable and complete solution sets admit powerful finite language-theoretic descriptions.  For polynomial equations over algebraically closed fields, classical algebraic geometry and computational commutative algebra give a much older and broader structural theory, whereas arithmetic rings exhibit severe undecidability phenomena.  For free associative algebras, general Diophantine solvability is undecidable, while a newer program investigates structural MR/JSJ phenomena for important restricted classes.

\section{Summary}
\fullwidthtext
The deepest form of solving an equation is usually not to find one point but to understand the geometry of all points.  The coordinate-object viewpoint turns solutions into homomorphisms.  Over free groups, sequences of homomorphisms lead to limit groups and {\bf R}-tree actions.  JSJ decompositions reveal splittings and modular automorphisms; shortening removes this internal flexibility; genuine degeneration produces proper limit quotients; and iteration yields an MR diagram.

The surface equation
$$
[x,y]=[a,b]
$$
shows why surface vertices are central.  Its once-punctured-torus vertex has Dehn-twist automorphisms that generate infinite families of solutions while preserving the boundary equation.  In a full MR construction, one first quotients conceptually by this modular freedom and only then regards the remaining limiting behavior as genuine degeneration.  Thus
$$
\plainboxed{\hbox{JSJ describes internal flexibility; MR organizes all families and degenerations.}}
$$

\section{Further reading}
\fullwidthtext
For the original decision procedures, see Makanin's papers on equations in free semigroups and free groups ([Ma1],[Ma2]).  Razborov developed the first systematic finite descriptions of all solutions
[Ra].  Sela's [Se1] presents a canonical MR diagram encoding the solution set and begins the geometric series leading to the solution of Tarski's problems.  Kharlampovich and Myasnikov give an  algebraic-geometric development ([Kh1],[Kh2]).  Serre's [Ser] is the standard starting point for Bass-Serre theory.
For formal-language descriptions, see work of Ciobanu, Diekert, Elder, Je\.z, Plandowski and collaborators ([CDE],[DJP]).
For associative algebras, Gr\"obner-Shirshov theory is classical [Bo-Ch], while the MR/JSJ viewpoint is considerably more recent.

\vglue 1.5pc
\chapter{Maximal Pairs}
\fullwidthtext
\medskip
In a similar way to the study of equations over groups [Ra1], with a finite system of equations $\Phi$
 over a free semigroup $FS_k =$  $<a_1,\ldots,a_k>$ it is
natural to associate a f.p.\ semigroup $S(\Phi)$. If the system $\Phi$ is defined
by the coefficients $a_1,\ldots,a_k$, the unknowns
$x_1,\ldots,x_n$ and the   equations:
\begin{align*}
u_1(a_1,\ldots,a_k,  x_1,\ldots,x_n) &=
v_1(a_1,\ldots,a_k,  x_1,\ldots,x_n)  \\
& \vdots \hfill \\
u_s(a_1,\ldots,a_k,  x_1,\ldots,x_n) &=
v_s(a_1,\ldots,a_k,  x_1,\ldots,x_n)  \\
\end{align*}
we set the associated f.p.\ semigroup $S(\Phi)$ to be:
$$
S(\Phi) \ = \ < \, a_1,\ldots,a_k,x_1,\ldots,x_n \, | \, u_1=v_1,\ldots,u_s=v_s \, >
$$

\noindent
Clearly, every solution of the system $\Phi$ corresponds to a homomorphism
(of semigroups)
$h:S(\Phi) \to FS_k$ for which $h(a_i)=a_i$, and every such homomorphism
corresponds to a solution of the system $\Phi$. Therefore, the study of
sets of solutions to systems of equations in a free semigroup is equivalent to
the study of all homomorphisms from a fixed f.p.\ semigroup $S$ into a free semigroup
$FS_k$,
for which a given prescribed set of elements in $S$ is mapped to a fixed
basis of the free semigroup $FS_k$.

Hence, as in studying sets of solutions to systems of equations over a free or a hyperbolic group
[Se1], to study sets of solutions to systems of equations over a free semigroup, we fix a f.p.\
(or even a f.g.) semigroup $S$, and study the structure of its set of homomorphisms into
a free semigroup, $FS_k$, that we denote, $Hom(S,FS_k)$.

Given a f.g.\ semigroup, $S$, we can naturally associate a group with it. Given a presentation
of $S$ as a semigroup, we set the f.g.\ group $Gr(S)$ to be the group with the presentation of $S$
interpreted as a presentation of a group. Clearly, the semigroup $S$ is naturally mapped into the
group, $Gr(S)$, and the image of $S$ in $Gr(S)$ generates $Gr(S)$. We set  $\eta_S: S \to Gr(S)$ to
be this natural homomorphism  of semigroups.

\noindent
The free semigroup, $FS_k$, naturally embeds into a free group, $F_k$. By the construction of the group, $Gr(S)$,
 every homomorphism of semigroups,
$h:S \to FS_k$, extends to a homomorphism of groups, $h_G: Gr(S) \to F_k$, so that: $h=h_G \circ \eta_S$.

Our goal is to study the structure of the set of homomorphisms (of semigroups), $Hom(S,FS_k)$. By construction, every
homomorphism (of semigroups), $h:S \to FS_k$, extends to a homomorphism (of groups), $h_G:Gr(S) \to F_k$. Therefore,
the study of the structure of $Hom(S,FS_k)$, is equivalent to the study of the structure of the collection of
homomorphisms of groups, $Hom(Gr(S),F_k)$, that restrict to homomorphisms of (the semigroup) $S$ into the
free semigroup (the $positive$ $cone$), $FS_k$.

By the techniques of section 5 in [Se1], with any given collection of homomorphisms of a f.g.\ group into a free
group, we can associate its Zariski closure, and with the Zariski closure one can associate canonically a
finite collection of limit groups, that are all (maximal) limit quotients of the given f.g.\ group,
so that every homomorphism from the given collection factors through at least one of the quotient maps from the
given f.g.\ group into the (finitely many) limit quotients. By (canonically) associating a finite collection
of maximal limit quotients with the set of homomorphisms, $Hom(Gr(S),F_k)$,  that restrict to
(semigroup) homomorphisms from $S$ to $FS_k$, we get the following basic theorem, which is the basis for our
approach to study the structure of $Hom(S,FS_k)$.

\vglue 1.5pc
\proclaim{Theorem 3.1} Let $S$ be a f.g.\ semigroup, and let $Gr(S)$ be the f.g.\ group that is associated
with $S$, by interpreting a semigroup presentation of $S$, as a presentation of a group. Let $\eta_S: S \to Gr(S)$
be the natural semigroup homomorphism, and note that $\eta_S(S)$ generates $Gr(S)$ as a group.

There exists a canonical finite  collection of (limit) pairs, $(S_1,L_1),\ldots,(S_m,L_m)$, where the $L_i$'s are
limit quotients of $Gr(S)$, and the semigroups, $S_i$, are quotients of the semigroup $S$ that generate
the limit groups $L_i$ as groups, with the following
properties:
\roster
\item"{(1)}" for each index $i$, $1 \leq i \leq m$,
there exists a (canonical) quotient map of pairs, $\tau_i:(\eta_S(S),Gr(S)) \to (S_i,L_i)$.

\item"{(2)}" by construction, every homomorphism of semigroups, $h:S \to FS_k$, extends to a map of pairs,
$h_P:(\eta_S(S),Gr(S)) \to (FS_k,F_k)$. For each such homomorphism of pairs, there exists an index $i$, $1 \leq i \leq m$,
and a homomorphism of pairs, $u_h:(S_i,L_i) \to (FS_k,F_k)$, for which: $h_P=u_h \circ \tau_i$.
\endroster
\endproclaim

\nfp Given a group $G$ its residually free quotient $RF(G)$ is a canonical quotient of G, where the kernel of the quotient map, $G \to RF(G)$, is the collection of
elements in $G$ that are mapped to the identity by every homomorphism: $f:G \to F_k$. Given a pair $(\eta(S),Gr(S))$ we denote
by $RF(\eta(S),Gr(S))$  the pair $(RF(\eta(S)),RF(Gr(S)))$, where $RF(\eta(S))$ is the image of $\eta(S)$ in the residually
free quotient, $RF(Gr(S))$.

Every pair homomorphism, $h:(\eta(S),Gr(S)) \to (FS_k,F_k)$,  factors through the residually free quotient: $RF(\eta(S),Gr(S))$. Hence,
analyzing all the pair homomorphisms, $h:(\eta(S),Gr(S)) \to (FS_k,F_k)$, is equivalent to
analyzing all the pair homomorphisms, $h:RF(\eta_S,Gr(S)) \to (FS_k,F_k)$.

This last collection of homomorphisms can be analyzed
precisely as it is done in the proof of theorem 7.2 in [Se1], and the conclusion of theorem 3.1 is identical to the statement of that theorem.

\hfill$\square$

In theorem 7.2 in [Se1], it is shown that with each f.g.\ group it is possible to associate
a canonical finite collection of limit groups, so that each homomorphism from the f.g.\ group into a
free group factors through at least one of the finitely many limit quotients. Theorem 3.1 is the analogue
of that theorem for semigroups. It reduces (canonically) the study of the structure of the set of semigroup
homomorphisms from a given semigroup to a free semigroup, $Hom(S,FS_k)$, to the study of the structure of
homomorphisms of finitely many pairs $\{(S_i,L_i)\}$ into $(FS_k,F_k)$, where the $L_i$'s are limit groups, and the
$S_i$'s are subsemigroups of the $L_i$'s that generate the limit groups $L_i$'s as groups.


\vglue 1.5pc
\chapter{Positive cones, their embeddings in standard cones, and shortenings}
\fullwidthtext
\medskip

To analyze homomorphisms of a f.g.\ semigroup into the free semigroup
we analyzed sequences of homomorphisms of
pairs $(S,G)$ into the standard pair $(FS_k,F_k)$. In studying such homomorphisms of pairs, the subsemigroup
$S$ is viewed as the positive cone in the ambient group $G$. To get a structure theory for the
entire collection of pair homomorphisms we are led to replacing the given cone with a standard cone. The structure
of a standard cone depends on the structure of the ambient group, and more specifically on the structure
of the tree that is obtained as a limit of a sequence of homomorphisms, and the dynamics of the action
of the ambient group on that tree.

In general, a standard cone can not be embedded into the ambient group $G$, but rather into an extension of $G$.
Still, a standard cone is essential in applying the shortening argument for semigroups, for any attempt to construct
JSJ decompositions, for separating free factors, and in general for constructing a Makanin-Razborov diagram
that is associated with a pair $(S,L)$.

Let $(S,L)$ be a pair of a limit group, $L$, and its subsemigroup, $S$, that generates it as a group. Let
$\{h_n:(S,L) \to (FS_k,F_k)\}$ be a sequence of pair homomorphisms that converges into a faithful action
of $L$ on a real tree $Y$. In this section we show how to extend the subsemigroup $S$ to a bigger semigroup,
by adding to it standard generating sets of the various components of the real tree $Y$. We further make
sure that the (finitely many) generators of the original subsemigroup $S$ can be expressed as positive words
in the standard generating sets that are associated with the components of the limit tree $Y$.

The standard generating sets that we add depend on the dynamics of the action of the limit group
$L$ on the limit tree $Y$. We start with the rather basic case, of an axial action of a free
abelian group on a line.

\vglue 1pc
\proclaim{Theorem 4.1} Let $(S,A)$ be a pair of a free abelian group group $A$ and a subsemigroup $S$ that generates $A$.
Let $\{h_n\}$ be a sequence of pair homomorphisms of $(S,A)$ into $(FS_k,F_k)$ that converges into an  axial action of
$A$ on a line $Y$. Let $A_0$ be the direct summand of $A$ that acts trivially on $Y$, and suppose that:
$rk(A) - rk(A_0)=\ell$.

Then there exists a $positive$ collection of generators, $a_1,\ldots,a_{\ell}$, in $A \setminus A_0$ so that:
\roster
\item"{(1)}" $A=A_0+<a_1>+\ldots+<a_{\ell}>$.

\item"{(2)}" there exists some index $n_0$, such that for every $n>n_0$: $h_n(a_1),\ldots,h_n(a_{\ell}) \in FS_k$.

\item"{(3)}" for each of the generators $s_1,\ldots,s_r$ of the semigroup $S$, $s_j$ can be expressed as a positive word in
$a_1,\ldots,a_{\ell}$ modulo an element in $A_0$. i.e., for every $j$, $1 \leq j \leq r$:
$$s_j=a_1^{m^j_1} \ldots a_{\ell}^{m^j_{\ell}} a_0(j)$$
where $m^j_i \geq 0$ and $a_0(j) \in A_0$.

\item"{(4)}" Wlog we can assume  that for every $j$, $1 \leq j \leq r$, and every $n>n_0$, $h_n(a_0(j)) \in FS_k$.
\endroster
\endproclaim

\nfp Let
$a_1,\ldots,a_{\ell}$ be elements in $A$, for which $A=A_0+<a_1>+\ldots+<a_{\ell}>$, and so that the elements $a_1,\ldots,a_{\ell}$
translate along the axis of $A$ (the real tree $Y$)
in a positive direction. Since the action of $A$ on its axis is an axial action, $a_j$ translates
a positive distance $\alpha_j$, and the real numbers $\alpha_j$, $1 \leq j \leq \ell$, are independent over the rationals.

Each of the fixed set of generators of the semigroup $S$, $s_1,\ldots,s_r$, can be written as a word in the elements $a_1,\ldots,a_{\ell}$
times an element in $A_0$.
We modify the elements $a_1,\ldots,a_{\ell}$ iteratively, so that the elements $s_1,\ldots,s_r$ can be represented as positive
words in the modified set of generators. If an element $s_j$ is already a positive word in the generators,
then the modifications of the generators that we perform, will change the positive word to another positive word in
the new generators.

Suppose that after rearranging the order of the given generators, $s$, one of the elements, $s_1,\ldots,s_f$,
is represented by the  word:
$$s \, = \, \Sigma_{j=1}^m \, k_ja_j \ - \ \Sigma_{j=m+1}^{\ell}t_ja_j$$
where the equality is in $A$ modulo its subgroup $A_0$,  the coefficients $k_j$  are non-negative, and $t_j$ are positive.
Wlog we may assume that $t_{m+1}$ is one of the maximal elements from the set $t_{m+1},\ldots,t_{\ell}$.
To prove the theorem we show that after finitely many modifications of the set of generators, $s$ can be represented
by another word in a new  set of positive generators, so that the absolute value of the negative coefficients is still
bounded by $t_{m+1}$, and the number of
appearances of the coefficient $-t_{m+1}$ is strictly smaller than the number of appearances of it in the original word.

Suppose that there exists an index $j$, $1 \leq j \leq m$, for which $k_j$ is positive, and  $\alpha_j>\alpha_{m+1}$. In that case we replace the
(positive) generator $a_j$ with a positive generator $\hat a_j=a_j - a_{m+1}$. With respect to the new
set of generators the coefficients $k_j$ are unchanged, $-t_{m+1}$ is replaced by $-t_{m+1}+k_j$ and $t_j$,
$j=m+2,\ldots, \ell$ are unchanged. In particular, all the new negative coefficients are bounded below by
the previous $-t_{m+1}$, and the number of negative coefficients that are equal to $-t_{m+1}$ is reduced by 1.

Suppose that for all the indices $j$, $j=1,\ldots,m$, for which $k_j$ is positive, $\alpha_j<\alpha_{m+1}$. Since $s$ is in the positive cone,
some of the coefficients $k_j$ are positive. After rearranging, we may suppose  that $k_1>0$.
We start by replacing the
generator $a_{m+1}$ with $\hat a_{m+1}=a_{m+1}-a_1$. This replaces $k_1$ by $\hat k_1=k_1-t_{m+1}$ and
keeps all the other coefficients
unchanged. If $\hat k_j=0$ we have reduced the number of generators that participate in the word that represents the element
$s$. If $\hat k_j>0$ we do the following:
\roster
\item"{(1)}" If for the new set of generators and coefficients for $s$, there exists an index $j$, $j=1,\ldots,m$,
for which $k_j$ is positive and $\alpha_j> \hat \alpha_{m+1}$, we replace the generator
$a_j$ with $a_j-a_{m+1}$, and hence replace the coefficient $-t_{m+1}$ with $-t_{m+1}+k_j$. This reduces the number of negative
coefficients that are equal to the previous $-t_{m+1}$ by 1.

\item"{(2)}" Suppose that for the new set of generators and for all $j$, $j=1,\ldots,m$, for which $k_j$ is positive, $\alpha_j< \hat \alpha_{m+1}$.
In that case we replace $\hat a_{m+1}$ by
$\hat a_{m+1}-a_1$, and hence replace $\hat k_1$ by $\hat k_1 - t_{m+1}$.
\endroster
After finitely many steps either the coefficient of the first generator becomes 0 or negative, or
we reduced the number of negative coefficients that are the biggest in their absolute value (i.e., equal to
$-t_{m+1}$).

If $\hat k_1 <0$ then we reduced the number of generators with positive coefficients in the representation of the
element $s$. Hence, we repeat the same steps where the elements that have non-negative coefficients are only
$a_2,\ldots,a_m$. Therefore, after finitely many steps the number of generators with negative coefficient that are the
biggest in their absolute value is strictly smaller. Finally, after finitely many steps we obtain a new set of
generators for which the element $s$ is represented by a positive word.

Note that if all the coefficients in an element $s$ were non-negative, then all the manipulations that we do keep them non-negative.
Hence, continuing to the other elements in the initial set:
$s_1,\ldots,s_r$, we get a generating set for which all these elements are represented by positive words (modulo $A_0$).

\hfill$\square$

Theorem 4.1 constructs a positive basis in case the ambient limit group is abelian and its limit action is axial.
The next case that we consider is the case of a pair $(S,L)$, in which the ambient limit group $L$ is a surface group, and $L$ acts
freely  on the limit (real) tree. In this case of an IET action of the surface group $L$,
we were not able to prove the existence of a positive standard cone. Instead we prove a weaker statement that enables one to apply
shortening arguments in the sequel  (shortening arguments are essential in the constructions of the JSJ decompositions, and
the Makanin-Razborov diagrams in the sequel).

\vglue 1pc
\proclaim{Theorem 4.2} Let $Q$ be a (closed) surface group, and let  $(S,Q)$ be a pair where $S=<s_1,\ldots,s_r>$ is a f.g.\
subsemigroup that generates $Q$
as a group.
Let $\{h_n\}$ be a sequence of pair homomorphisms of $(S,Q)$ into $(FS_k,F_k)$ that converges into a free  (IET)
action of $Q$
on a real tree $Y$.


There exists a sequence of automorphisms $\{\varphi_{\ell} \in Aut(Q)\}_{\ell=1}^{\infty}$ with the following properties:
\roster
\item"{(1)}" for each index $\ell$, there exists $n_{\ell}>0$, so that for every $n>n_{\ell}$, and all the generators of
the subsemigroup $S$, $s_1,\ldots,s_r$, $h_n \circ \varphi_{\ell}(s_j) \in FS_k$.
(note that the automorphisms $\varphi_{\ell}$ do
not preserve the subsemigroup $S$, but the image of $S$ under the twisted homomorphisms $h_n \circ \varphi_{\ell}$
remains
in the standard positive cone in $F_k$ for large enough $n$).

\item"{(2)}" for every $n>n_{\ell}$ and every $j$, $1 \leq j \leq r$:
$$\frac {d_T( h_n \circ \varphi_{\ell}(s_j), id)}
{\max_{1 \leq j \leq r}d_T(h_n(s_j),id.)} < \frac {1}{\ell}$$
where $T$ is a  Cayley graph of the coefficient free group $F_k$ with respect to a fixed set of generators.
\endroster
\endproclaim

\nfp We view $FS_k$ as the positive cone inside $F_k$. This gives a positive orientation to any element in  $FS_k$. With every element $w \in F_k$ we
naturally associate the shortest path $p_w$ from the identity to $w$ in the Cayley graph of $F_k$ w.r.t.\ to the given set
of free generators, $a_1,\ldots,a_k$. Clearly each such path can be cut canonically into finitely many positive and negative  subpaths.

Since the elements $h_n(s_i) \in FS_k$ for every $n$, they are in the positive cone. Hence, each of the (limit) segments $[y_0,s_i(y_0)] \subset Y$,
where $y_0$ is the basepoint in $Y$,  is the limit of positive elements or paths in the Cayley graph of $F_k$, so we define the segments
$[y_0,s_i(y_0)]$ to be
positively oriented.

Since $s_1,\ldots,s_r$ generate $Q$ as a group, the orbit of the segments $[y_0,s_i(y_0] \subset Y$, under the action of $Q$, cover the entire tree $Y$.
Hence, every interval $[y_1,y_2] \subset Y$ can be cut canonically into finitely many subintervals, such that every subinterval is either positively
or negatively oriented. i.e., each such subsegment is the limit of either positive or negative elements.

The generators $s_1,\ldots,s_r$ of the subsemigroup $S$ generate the surface group $Q$. The action of $Q$ on the limit tree $Y$
is an IET action, which is in particular a geometric action. Hence, there exists a finite subtree $R$ of $Y$, such that if one applies
the Rips machine to the action of the pseudogroup generated by the restrictions of the actions of the generators
$s_1,\ldots,s_r$ to the finite tree $R$, the output is
a standard IET pseudogroup. i.e., the output is supported on a union of finitely many intervals that are all part of one positively
oriented  interval, and these finitely many intervals
are divided into finitely many subintervals,  a permutation of these subintervals, and a standard set of generators of
$Q$ that map
each subinterval to its image (dual) under the permutation. Note that the generators of $Q$ that are associated with the IET
presentation, are obtained by gluing the finitely many intervals on which the IET is supported into a (finite) tree, and
then the generators of $Q$ are the elements that map a subinterval into its dual. Hence, the generators of $Q$
that are associated with the IET transformation need not be positively oriented (when acting on the tree $Y$)  in case the IET is supported on more than one
connected interval. This is in contrast with the positive orientation of the  original generators $s_1,\ldots,s_r$, and of
each of the
subintervals that define the IET transformation.


For presentation purposes we first assume that the endpoints of the intervals in the interval exchange (that is obtained by the Rips
machine from the action of $Q$ on $Y$) belong to a single orbit under $Q$ (i.e., that the generators of the interval exchange act
transitively on the endpoints of the intervals in the IET), and that the base point in $Y$ (i.e., the point in $Y$ which is
a limit of the identity element in the Cayley graph of $F_k$) belongs to that orbit.



\smallskip
The action of the surface group $Q$ on the tree $Y$ is reduced to an interval exchange transformation that is based on
positively oriented intervals $I_1,\ldots,I_g$.
By our assumption, the points that are at the beginnings and at the ends of the subintervals that define the interval exchange
are all in the orbit of the basepoint in $Y$. Let $q_0,\ldots,q_f \in I$ be the points that are at the beginnings and at the ends
of the subintervals that define the IET transformation on the intervals $I_1,\ldots,I_g$.
We look at those indices $i$, $1 \leq i \leq f$, for which the segment $[q_{i-1},q_i]$ is supported on one of the segments
$I_1,\ldots,I_g$. For each such index $i$, there is a unique element in $Q$ that maps the subinterval $[q_{i-1},q_i]$
to a
(positively oriented) subinterval $[y_0,y_i]$, where $y_0$ is the basepoint in $Y$. For each $i$, $1 \leq i \leq f$, we set
$v_i \in Q$ to be the element that maps $y_0$ to $y_i$.

By construction, for large enough index $n$, $h_n(v_i) \in FS_k$,
for every $i$, $1 \leq i \leq f$. The elements in $Q$, that restrict to the elements that are associated with
the interval exchange transformation, i.e., the elements that map  subintervals to their companions
(after gluing the intervals $I_1,\ldots,I_g$ and obtaining a finite tree), generate the surface group $Q$.
Each of these elements is in the subgroup that is generated by  $v_1,\ldots,v_f$. Hence, $v_1,\dots,v_f$ generate $Q$.

For each $j$, $1 \leq j \leq r$, the segment $[y_0,s_j(y_0)]$ is positively oriented, and $s_j$ can be presented as a word in the
elements $v_1,\ldots,v_f$, $s_j=w_j(v_1,\ldots,v_f)$.
The words $w_j$ need not be positive words in the elements $v_1,\ldots,v_f$.
For each index $j$, $w_j(v_1,\ldots,v_f)$ represents a finite path in the tree
$Y$, that starts at $y_0$ and ends at $s_j(y_0)$. We denote this path $p_{w_j}$.
Since the words $w_j$ need not be positive words, $p_{w_j}$ need not be positively oriented.

Each path $p_{w_j}$ is supported on a finite subtree of the real tree $Y$ that we denote $T_{w_j}$. We view $T_{w_j}$ as
a combinatorial
tree and not as a metric tree. To $T_{w_j}$ we add a finite collection of vertices:
\roster
\item"{(1)}" a vertex for the base point of $Y$ (the initial point of $p_{w_j}$), and for the point
$s_j(y_0)$ (the terminal point of $p_{w_j}$).

\item"{(2)}" a vertex for each root and each branching point in $T_{w_j}$.

\item"{(3)}" each segment in $T_{w_j}$, being a subsegment of $Y$, can be divided into finitely many segments with either positive
or negative orientation. If there exists an edge in $T_{w_j}$  which is not oriented (i.e., the edge can be cut into subsegments
with both positive and negative orientations), then the word $w_j$ can be strictly shortened as a word in $v_1,\ldots,v_f$, and
still represents $s_j \in Q$. Therefore, we can assume that every edge in $T_{w_j}$ is oriented.

\item"{(4)}" $p_{w_j}$ starts at the base point in $Y$, and ends in
$s_j(y_0)$. While moving along $p_{w_j}$, at certain segments of the path the distance to the base point $y_0 \in Y$ increases,
and at other segments it decreases.  At each point which is the boundary of such (increasing or decreasing) segments
along
$p_{w_j}$ (i.e., points in which the distance to $y_0$ changes from increasing to decreasing or vice versa), we add a vertex to
$T_{w_j}$.
\endroster
In the sequel we denote the paths $p_{w_j}$ that are associated with the words $w_j$ on $T_{w_j}$.

We assumed that  finitely many steps of Makanin's procedure enable one to reduce the (generalized) equation that is obtained
from an original system of equations over $FS_k$ to an interval exchange transformation on the union of finitely many
intervals. i.e., to a pseudogroup of partial isometries of the union of intervals, that is defined by dividing each interval in two ways
to finitely many subsegments, dividing the obtained segments into disjoint pairs of the same lengths, and identifying the two segments
(that are called bases) in each pair.

Given one of the finitely many intervals, that we denote $I$,  we can define a Dehn twist. The interval $I$ has a (positive) orientation, and we
perform the Dehn twist at the positive end of $I$, that we denote $e$. There are two bases $J_1$ and $J_2$ that are supported by subintervals of $I$,
and are adjacent to $e$. If the length of $J_1$ is equal to that of $J_2$, we can identify the base that $J_1$ is paired with to the base that
$J_2$ is paired with, and cut $I$ at the other end of $J_1$ and $J_2$.

Hence, we may assume that the length of $J_1$ is strictly smaller than the length of $J_2$. In that case we cut $J_2$ and its dual base
at the other end of $J_1$,
and transfer $J_1$ to the base that is paired with $J_2$. Finally we cut $I$ at the point in the other end of the original base $J_1$.

This is called a Dehn twist. Note that the number of bases remains identical after a Dehn twist, and the length of one of the supporting intervals and
the length of a pair of bases is strictly reduced.

At this point we apply a sequence of Dehn twists on the generators of the interval exchange (which is the output of the Rips
machine from the action of the surface group $Q$ on the real tree $Y$). We perform the sequence of Dehn twists from the positive
side of the interval on which the IET is based. Hence, the sequence of IETs that we obtain are supported on a
decreasing sequence of
subintervals of the original interval $I$ that supports the original IET. All these supporting subintervals share the endpoint
$y_0$ which is the base point in $Y$. Since the action of $Q$ on $Y$ was assumed to be a minimal IET action, the lengths of
the supporting intervals has to approach 0.

With each of the obtained IETs we associate a positive generating set for a subsemigroup of $Q$, in a similar way to the
association of $v_1,\ldots,v_f$ with the original IET. We denote the corresponding set of generators $v^t_1,\ldots,v^t_f$ (note
that the number of generators of the corresponding semigroups can only decrease, so by omitting finitely many of them, we can
assume that their number is fixed).

\vglue 1pc
\proclaim{Lemma 4.3} For every index $t \geq 1$, each of the elements $v_1,\ldots,v_f$ is contained in the subsemigroup
that is generated by: $v^t_1,\ldots,v^t_f$.
\endproclaim

\nfp By the definition of a Dehn twist that is performed from the positive side of the supporting interval, for every
index $t \geq 1$, each of the elements $v^{t}_1,\ldots,v^{t}_f$ can be written as a (concrete) positive word in the
elements $v^{t+1}_1,\ldots,v^{t+1}_f$. Hence, the claim of the lemma follows by induction.

\hfill$\square$

By lemma 4.3 each of the elements $v_1,\ldots,v_f$ can be written as a positive word in the elements
$v^t_1,\ldots,v^t_f$. By substituting these positive words in the words $w_j$ we obtain presentations of the elements
$s_j$, $1 \leq j \leq r$, in terms of the  generating sets $v^t_1,\ldots,v^t_f$:
$$s_j \, = \, w_j(v_1,\ldots,v_f) \, = \, w_j^t(v^t_1,\ldots,v^t_f).$$
Since the elements $v_1,\ldots,v_f$ and $v^t_1,\ldots,v^t_f$ are positively oriented, and the elements $v_j$ are presented
as positive words in the generators $v^t_1,\ldots,v^t_f$, the paths that are associated with the words $w_j^t$ in the tree $Y$,
that we denote
$p_{w_j^t}$, are identical to the paths $p_{w_j}$, that are associated with the words $w_j$.

While shortening the lengths of the elements $v^t_1,\ldots,v^t_f$, the ratios between these lengths  may not be bounded. To prove theorem
4.2 we first need to replace the elements $v^t_1,\ldots,v^t_f$ by elements that are not longer than them, and for which the ratios
between their lengths are universally bounded. We start by proving the existence of such generators assuming a global
bound on the periodicity
of the images of the given set of generators of the semigroup $S$, $s_1,\ldots,s_r$, under the sequence of
homomorphisms $\{h_n\}$. In the sequel we prove a general version of the proposition omitting the bounded periodicity
assumption.

\vglue 1pc
\proclaim{Proposition 4.4} Suppose that there exists an integer $c_p$, such that the periodicity of the elements
$h_n(s_1),\ldots,h_n(s_r)$ is bounded by $c_p$ for all integers $n$.
i.e., for every $n$, and every $j$, $1 \leq j \leq r$, $h_n(s_j)$ can not be written as a word of the form:
$s_j=w_1 \alpha^{c_p+1} w_2$, where $w_1,w_2,\alpha \in FS_k$, and $\alpha$ is not the empty word.

For every $t \geq 1$, the elements
$v^t_1,\ldots,v^t_f \in Q$ can be replaced by elements $u^t_1,\ldots,u^t_{g_t} \in Q$, with the following properties:
\roster
\item"{(1)}" $f \leq g_t \leq e(f)$ (where $e(f)$ is a positive integer that depends only on $f$ and not on $t$).

\item"{(2)}" for every index $t$, the semigroup that is generated by $u^t_1,\ldots,u^t_{g_t}$ in $Q$ contains the semigroup that
is generated by $v^t_1,\ldots,v^t_f$ in $Q$.

\item"{(3)}" for every $t \geq 1$, and for large enough index $n$, $h_n(u^t_i) \in FS_k$ for every index $i$, $1 \leq i \leq g_t$.

\item"{(4)}" the tuples $u^t_1,\ldots,u^t_{g_t}$ belong to finitely many isomorphism classes (under the action of $Aut(Q)$).

\item"{(5)}" because of (2) each of the elements $v^t_1,\ldots,v^t_f$ can be represented as a positive word in
the elements $u^t_1,\ldots,u^t_{g_t}$. Substituting these words in the words $w^t_j$, for each $t \geq 1$ and each $j$, $1 \leq j \leq r$,
$s_j=z^t_j(u^t_1,\ldots,u^t_{g_t})$. With each word $z^t_j(u^t_1,\ldots,u^t_{g_t})$ we can associate a path, $p_{z^t_j}$
 in the tree $Y$, and this path is
identical to the path $p_{w^t_j}$, that is identical to $p_{w_j}$.

\item"{(6)}" there exist positive constants $d_1,d_2$ that depend only on $f$ and $c_p$, such that for every $t \geq 1$,
there exist a sequence of indices (that depend on $t$) $1 \leq i_1 < i_2< \ldots < i_{b(t)} \leq g_t$, such that
for every $1 \leq  m_1 < m_2 \leq  b(t)$:
$$d_1 \cdot length (u^t_{i_{m_1}})  \, \leq \, length (u^t_{i_{m_2}}) \, \leq \,
 d_2 \cdot length (u^t_{i_{m_1}}).$$
Furthermore,  for every index $i$, $1 \leq i \leq g_t$ for which $i \neq i_m$, $m=1,\ldots,b(t)$:
$$ 10g_t \cdot d_2 \cdot length (u^t_i)  \, \leq \, length (u^t_{i_1}) $$

\item"{(7)}"
For each index $t$,
in the words $z^t_j$, $1 \leq j \leq r$, in a  distance bounded by some global constant $q=q(f)$ (that does not depend on $t$),
either before or after the occurrence of an
element $u^t_i$, for $i \neq i_1,\ldots,i_{b(t)}$, appears one of the elements $u^t_{i_m}$, $1 \leq m \leq b(t)$.

\endroster
\endproclaim

\nfp To get a new set of elements that satisfy the conclusion of the proposition we need to modify the standard
sequence of Dehn twists that
we used to get the elements $v^t_1,\ldots,v^t_f$ (which are Dehn twists that are performed on the pair of  bases that
are adjacent
to the vertex at the positive end of the interval on which the interval exchange transformation is supported).

We start by dividing the generators $v^t_1,\ldots,v^t_f$ into finitely many sets according
to their length. We order the elements from the longest to the shortest. We place a separator between consecutive sets
whenever there is a pair of consecutive elements (ordered according to length) that satisfy:
$$length(v^t_i) \geq c_1(f,c_p) \cdot length(v^t_{i+1})$$
where $c_1(f,c_p)=4fc_p$ and the length of an element is its displacement in the real tree $Y$.

Clearly the elements $v^t_i$ are divided into at most $f$ sets. If there are no separators, we set
$g_t=f$, and $u^t_i=v^t_i$, $i=1,\ldots,f$. Suppose that there is a separator. In that case our goal is to construct
a procedure that will iteratively lead to elimination of the  separators. We call the elements that are in the shortest set $short$, and all the other elements $long$.

The interval exchange transformation is defined by dividing the intervals that support it into bases, where each interval is divided into bases in two
different ways. The bases are divided into disjoint pairs, and the transformation is defined by identifying pairs of bases. By definition,
the elements are defined to be the group elements that map the end or the beginning of a base to the next end or beginning of base on the same interval.
By construction every base is a product of (distinct) elements.

Note that along the process we modify pairs of bases, and compare their lengths. Only at the end of the procedure we define the set of new elements
from the set of obtained bases, and using
the same method that was used to construct the original elements from the original bases, we construct the new elements from the obtained bases, and the
constructed new elements
satisfy the conclusion of the proposition.

\smallskip
We call a pair of bases $long$ if it is a product of
elements that at least one of them is long (equivalently, if one, and hence both, of the appearance of the pair of bases supports a long element).
The other pairs of bases are called  $short$.

If all the bases are short, the proposition follows. Hence, there are long pairs of bases. Note that it may that there are no short pairs of bases.

The goal of the procedure to prove proposition 4.4 is to gradually reduce the number of long pairs of bases. We do that by iteratively modifying
the sets of long and short pairs of bases, where the modifications cut long pairs of bases and may introduce new short pairs of bases. Short pairs of bases are
unchanged along the procedure, and the number of new short pairs of bases that are created along the procedure is universally bounded by the number of
original (long and short) pairs.

In addition, some of the (created) new short pairs are called $secondary$ $short$. Their length can be much smaller than the lengths of the original short elements.
Also, when the procedure terminates, the new bases define a new set of elements that satisfy the conclusions of the proposition,
and this set  may contain secondary short elements, that are much smaller than the original short elements.
The procedure guarantees that in a bounded
(combinatorial) distance from the appearance of these secondary short elements
in the words $z^t_j$ there will always appear a short  element. These secondary short elements are controlled in the second part of part (6) in
the statement of the proposition.

At this point we modify the standard Dehn twists that are performed on an interval exchange basis, so that the
performed Dehn twists do not change short pairs of bases.  Long pairs of bases are gradually shortened until
they are replaced by short and secondary short pairs of bases.

Let $I$ be the interval that supports the IET transformation, and
let $pv$ be the vertex at the positive end of $I$. On the given IET transformation we perform the following operations:
\roster
\item"{(1)}" Suppose that  the two bases that are adjacent to $pv$ are long.
Let $b_1$ be the longer base that is adjacent to
$pv$.
Suppose that the other end of $b_1$ is not covered by a long base. In this case we cut $b_1$ at the end of the last long base
that completely overlaps with $b_1$  and perform Makanin's $entire$ $transformation$ over the base $b_1$
(see [Ca-Ka] for Makanin's entire transformation).

Let $length_m$ be the minimal length of a short element and let $length_M$ be the maximal length of a short base or element.
If the length of what is left from $b_1$ after the
entire transformation is bounded above by $\frac {length_m} {c_1(f,c_p)}$ we set the pair of bases that is left from $b_1$ to be secondary short. If it is at least
$\frac {length_m} {c_1(f,c_p)}$ it is set to be short.

In this case (case (1)) we reduced the number of pairs of long bases by 1 and added a pair of short or secondary short bases. If the pair of bases that were left from
the cut of $b_1$ is short, then the set of short bases has changed, and we change accordingly the separator between long and short pairs of bases. i.e., we
start with the set of short bases and the original minimal length elements, and define the new separators between the new set of bases accordingly.
It may be that some pairs that were considered
previously to be long may become short. We also change the minimal and maximal lengths of short elements and bases, $length_m$ and $length_M$,
 according to the new set
of short bases and original short elements.

\item"{(2)}" Suppose that  the other end (the beginning) of $b_1$ is covered by another long base. In this case
we perform Makanin's entire transformation
over the carrying base $b_1$,  and continue as in part (1).

If the length of what is left from $b_1$ after the
entire transformation is bounded above by $\frac {length_m} {c_1(f,c_p)}$ we set the pair of bases that is left from $b_1$ to be secondary short. If it is at least
$\frac {length_m} {c_1(f,c_p)}$ and at most $c_1(f,c_p) \cdot lenth_M$ it is set to be short. If it is bigger than $c_1(f,c_p) \cdot length_M$ it is set to be long.

If what left from $b_1$ is short we update the set of short bases and $length_m$ and $length_M$ as we did in part (1).
If what left from $b_1$ is short or secondary short, the number of
long pairs of bases is reduced by 1. If what left is long, the number of long, short and secondary short  pairs did not change and we continue to the next step.

\item"{(3)}" Suppose that the bases that are adjacent to $pv$ are a short and a long base. Suppose that the length of the
long base that is
adjacent to $pv$, that we denote $b_1$,  is bigger than the sum of the lengths of the short bases that are adjacent to $pv$
and the first long base that is adjacent to them.
In that case we perform Makanin's entire transformation precisely as in cases (1) and (2).

\item"{(4)}" Suppose that a long and a short base are adjacent to $pv$, and that the length of the long base that is adjacent to $pv$,
that we denote $b_1$,
is smaller than the sum of the short bases that are adjacent to $pv$ and the first long base that is adjacent to these short
bases. We set $b_2$ to be the first long base that is adjacent to the short bases that are adjacent to $pv$.

If $b_1$ and $b_2$ are not a pair of bases we do the following.
We first transfer the short bases that are adjacent to $pv$ using the base $b_1$, i.e., we move the short bases that are adjacent to $pv$ to the
dual base of $b_1$. Then we cut $b_1$ after these short bases,
and  perform an entire transformation
over the base $b_2$ precisely as we did in cases (1) and (2). Note that if what left from $b_1$, after we transfer the short bases that are
adjacent to $pv$ and $b_2$, has length smaller than $c_1(f,c_p)$ times the length of the maximal  short element or base, we set this remainder of $b_1$ and its pair
to be $short$. Otherwise it remains $long$.

Suppose that $b_1$ and $b_2$ are a pair of bases. This implies that $b_1$ (and $b_2$) contains a large periodic word, which is
a power of
the chain of short bases that are adjacent to $pv$. We assumed that the periodicity of the words $h_n(s_j)$ is
bounded by $c_p$ for all $n$, and all $j$, $1 \leq j \leq r$. Hence, the length of $b_1$ (and $b_2$) is bounded by
$f \cdot (c_p+1)$ times the maximal length of a short element. This implies that there can not be a separator between the
length of an element that is supported by $b_1$ or $b_2$ and the collection of short elements. Therefore,
in this case $b_1$ (and $b_2$) can not be long bases, a contradiction. Hence, in case (4) $b_1$ and $b_2$ can not be a pair
of long bases (assuming bounded periodicity).


\item"{(5)}" Suppose that two short bases are adjacent to $pv$. Let $b_1$ and $b_2$ be the two long bases that are adjacent to
the two sequences of short bases that are adjacent to $pv$. $b_1$ and $b_2$ can not be a pair of bases because of our bounded
periodicity assumption, as we argued in part (4).
 Suppose that $b_1$ is longer than $b_2$, and that the sequence of
short bases that are adjacent to $pv$ and $b_1$ is longer than the sequence that is adjacent to $pv$ and $b_2$. In that case,
we first cut $b_2$ (and the base that is paired with it) into two bases at the endpoint of $b_1$.

We declare the part of $b_2$ that overlaps with short bases
before it overlaps with $b_1$ to be  $short$ or secondary short. What left from $b_2$ is long if its length is bigger than $c_1(f,c_p)$ times the length of the
maximal short element, and short otherwise. Then we perform an entire transformation along $b_1$ precisely as we
did in cases (1) and (2). In this case we either reduce the number of long pairs of bases (if what left from $b_2$ is short or if
the other end of $b_1$
overlaps with a
short base), or we added a short or a secondary short base (the part that left from $b_2$ after cutting it at the beginning of $b_1$),
 and left behind at least two short bases, that are not going to take part in the rest of the procedure,
so the number of short and secondary short
bases that participate
in the next iterations in this part (before we change the set of long bases) is reduced by
at least 1.

\item"{(6)}" With the notation and the assumptions of part (5),
suppose that the sequence of
short bases that are adjacent to $pv$ and $b_1$ (the longer base) is shorter than the sequence that is adjacent to $pv$ and $b_2$.
In that case,
we first cut $b_1$ (and the base that is paired with it) into two bases at the endpoint of $b_2$.
We declare the part of $b_1$ that overlaps with short bases
before it overlaps with $b_2$ to be either short or secondary short, depending whether its length is bigger than $\frac {length_m} {c_1(f,c_p)} length_m$ or not.
Then we perform an entire transformation along the longer between the
remaining of $b_1$ and $b_2$, precisely as we
did in cases (1) and (2). As in part (5), in this case we either reduce the number of long pairs of bases,
or we added a  short or a secondary short  base and left behind at least two short bases, so the number of short and secondary short bases that take
part in the rest
of the procedure is reduced by
at least 1.

\endroster

In step (1) the number of pairs of long bases is reduced by at least one, and the number of short and secondary short bases is increased by
at most the reduction in the number of long bases. In step (2) there may be no change in the number of bases, and if there is
a change the number of short or secondary short bases is increased by at most the reduction in the number of long bases. The outcome of
steps (3) and (4) in terms of the number of bases is identical to that of steps (1) and (2). The outcome of steps (5) and (6)
on the number of bases is similar to that of steps (1) and (2) with an additional increase of the number of short or secondary short
bases by 1, but an  additional reduction of  the number of active
short or secondary short bases (i.e., the number of bases that take part in the next steps of the iterative procedure) by at least 1.


We run this procedure iteratively until there are no long bases. i.e., until all the bases are either short
or secondary short, the number of bases can grow to at most $2f$. At each step of the procedure the set of pairs of bases define a collection
of elements, that are simply the elements that translate between the consecutive endpoints of bases. Since there are at most $2f$ pairs of bases, there are at
most $4f$ elements that are associated with  the set of bases at each step of the procedure.

We denote the collection of elements that are associated with the set of bases in the terminal step of the procedure, $u^t_1,\ldots,u^t_{g_t}$,
where $g_t \leq 4f$. The pairs of bases at each step of the procedure are obtained from the pairs of bases in the previous step using one of the moves
(1)-(6). Each of the elements that are defined by the previous pairs of bases, can be written as a (positive) word in the
elements that are defined by pairs of bases that are obtained using each of these moves (1)-(6). Hence, using a finite induction,
the semigroup that is generated by the final collection of elements
contains the semigroup that is generated by the initial collection of elements.

For large $n$, $h_n(u^t_i) \in FS_k$ since the new elements were obtained using a sequence of the moves (1)-(6). Since the new elements define an IET, the
closed surface group $Q$ is generated (as a group) by the terminal elements
$u^t_1,\ldots,u^t_{g_t}$. Since they define an IET and the action of $Q$ is free, $Q$ admits a presentation with the terminal generators, and boundedly many
relations that are obtained from the generators of the IET on the endpoints of the bases. Hence, there are only finitely many possible presentations of $Q$
using the terminal generators, so the terminal generators belong to finitely many orbits under the action of $Aut(Q)$.

The terminal set of pairs of bases can be either short or secondary short. We say that an element $u^t_i$ is short if its length is at least
$\frac {1} {2{c_1(f,c_p)}^2}$ times the minimal length of a short base. Otherwise an element $u^t_i$ is a secondary short element. Since there are short
pairs of bases, there are short terminal elements $u^t_i$.

We set the indices $i_1,\ldots,i_b$ in part (6) of the proposition
to be the indices of the short elements. Since there are boundedly many of them, their lengths clearly satisfy the
first inequality in part (6) of the proposition, and the lengths of the secondary short elements satisfy the second inequality in part (6).

To prove part (7) of the proposition note that if there are no secondary short elements the claim is empty. Also note that all the elements
that the procedure starts with, i.e., all the elements, $v^t_1,\ldots,v^t_{f(t)}$,
that are defined using the original bases are either short or long, so the conclusion of part
(7) holds for the original set of elements.

A short (and not secondary short) element, $\hat u$, that is created along the iterative procedure to construct the final elements $\{u^t_i\}$,
can be cut to at most the number of final elements
new elements, so at least one of these new elements has to be short (and not secondary short), and part (7) follows for all the elements that appear in the
presentation of $\hat u$ as a product of final elements from the set $\{u^t_i\}$.

Let $\hat u$ be a long element that is created along the iterative process to construct the elements $\{u^t_i\}$. If $\hat u$ is replaced in the next step
in the process by a product of short and secondary short elements, then at least one of these elements has to be short, and their number is bounded by
the final number of elements, so part (7) follows for all the elements that appear in the presentation of $\hat u$ as a product in the final set of elements
$\{u^t_i\}$.

Otherwise either $\hat u$ is unchanged in the next step of the procedure, or it is replaced by a product that contains a new long element. Suppose that
$\hat u$ is replaced by a product that starts and ends with a long element. In that case the number of short and secondary short elements between
two long elements in the product is bounded by the total number of elements $\{u^t_i\}$, which is universally bounded by $e(f)$ - a function of the original
number of elements, so part (7) follows by an iterative application of the claim for long elements.

Suppose that $\hat u$ is replaced by a product that contains a long element, and starts with a short or a secondary short element. Since we assumed that
the images of the elements $s_1,\ldots,s_r$ contain only bounded periodicity, a long base can not overlap with itself with a short (or secondary short)
displacement. Hence, when the product that replaces the long element $\hat u$ starts with a short or a secondary short element, the effective number of
long pairs of bases is reduced by at least 1. This means that when an entire transformation (Dehn twist) is operated on the base that supports the
long element with which the product starts, part (2) of the procedure is applied, and the number of long pairs is reduced by 1.

Let $ba$ be the maximal number of pairs of bases along the whole procedure, i.e., the maximal number of pairs of bases for all $t$,
and all the finite iterative
procedures that starts at each step $t$.
The number of times
such a reduction can happen is bounded by the initial number of long pairs, which is bounded by $ba$. Hence, we get a universal
bound on the number of secondary short elements that appear before a short or a long element in  the beginnings of
the presentations of the long element $\hat u$ as a product
of elements that are created in all the next steps of the procedure.

Suppose that $\hat u$ is replaced by a product that terminates with a long element followed by a product of secondary short elements. The number of short
elements is bounded by $4ba$.
To prove part (7) of the proposition, suppose that in one of the next steps,
the terminating long element in the product that
represents $\hat u$ is replaced with a product that terminates with a long element followed by a product of secondary short elements. Again, there is a
universal bound on the
number of these new secondary short elements, so our goal is to prove that such consecutive replacements can occur a bounded number of times.

Recall that in the iterative process that produces the elements $\{u^t_i\}$, steps (5) and (6) can occur at most once before the number of long
pairs is reduced by at least 1. In a part of the process in which consecutively a long base is adjacent to the terminal point (that we denoted $pv$),
usually after steps (5) or (6) were
applied, and until the number of pairs of long elements is reduced, the number of times that secondary short bases are adjacent to the terminal point
is bounded by ${ba}^2$, and the number of times that the terminal part of a long base is replaced by a product that terminates
with a secondary short base is bounded by the number of pairs of long bases, that is bounded by $ba$.

Therefore, in each step of the procedure along such a part of the procedure,
a long element $\hat u$ is replaced by a product, and if the product contains a long
element, then this long element is followed by a consecutive sequence of secondary short elements before a short or a long element,
and the length of such consecutive sequence of secondary short elements is bounded by $4{ba}^2$.

This  bound is for a part of the process in which the number of long pairs of bases is reduced by 1, and the number of such
parts is bounded by the initial number of long pairs of bases.
Hence, altogether, in each step of the procedure, the long element $\hat u$ is replaced by a product in which a consecutive sequence of secondary short elements
has length that is bounded by  $4{ba}^3$.
This proves part (7) of the proposition.

\hfill$\square$

The next proposition is the key in our proof of theorem 4.2. Since we proved proposition 4.4 under a bounded
periodicity assumption, and assuming that the beginnings and ends of the subintervals in the IET transformation belong to the same orbit under
the action of the surface group $Q$,  we state and prove it under the same  assumptions, and generalize it in the sequel.
It implies
that the words $z_j^t$ (in part (5) of the statement of proposition 4.4) can be replaced by words $\hat z^t_j$
with uniformly bounded
cancellation.

\vglue 1pc
\proclaim{Proposition 4.5}
Suppose that there exists an integer $c_p$, such that the periodicity of the elements
$h_n(s_1),\ldots,h_n(s_r)$ is bounded by $c_p$ for all integers $n$.

With the notation of proposition 4.4, there exists a constant $C>0$, so that for every
index $t \geq 1$ the words $z^t_j$, $1 \leq j \leq r$, can be replaced by words: $\hat z^t_j$ with the following properties:
\roster
\item"{(1)}"
 As elements in $Q$:
$\hat z^t_j(u_1^{t},\ldots,u^{t}_{g_t})=z^{t}_j(u_1^{t},\ldots,u^{t}_{g_t})$.

\item"{(2)}" $\hat z^t_j$ is obtained from $z^{t}_j$ by eliminating distinct pairs of subwords. Each pair of  eliminated
subwords corresponds to two subpaths of the path
$p_{z^{t}_j}$ that lie over the same segment in the tree $T_{w_j}$, where the two subpaths have opposite orientations.

\item"{(3)}" With the word $\hat z^t_j(u^t_1,\ldots,u^t_{g_t})$ we can naturally associate a path in the
tree $Y$, that we denote, $p_{\hat z^t_j}$. The path $p_{\hat z^t_j}$ can be naturally divided into
subsegments according to the appearances of the subwords $u^t_i$ in the word $\hat z^t_j$.

Let $DB_{\hat z^t_j}$ be the number of such subsegments that are associated with subwords $u^t_i$ in $p_{\hat z^t_j}$, that at least
part of them is covered more than once by the path $p_{\hat z^t_j}$. Then for every $t \geq 1$ and every $j$, $1 \leq j \leq r$,
$DB_{\hat z^t_j} \leq C$.
\endroster
\endproclaim

\nfp Suppose that such a constant $C$ does not exist. Then for every positive integer $m$, there exists an index $j_m$,
$1 \leq j_m \leq r$, and an index $t_m>1$, so that for every possible choice of words $\hat z^{t_m}_{j_m}$ that satisfy
parts (1) and (2), part
(3) is false for the constant $C=m$.

By passing to a subsequence, we may assume that $j_m$ is fixed, and we denote it $j$. We can further assume that for
the subsequence the integer $g_{t_m}$, that counts the number of generators of the semigroup that is constructed in step
$t_m$ according to the procedure that is described in proposition 4.4 is fixed and we denote it $g$. We may also
assume that the integer  $b(t_m)$ and the sequences of indices $i_1,\ldots,i_{b(t_m)}$, that are associated with
the sets of generators of the semigroups that are constructed according to proposition 4.4 along the subsequence
$\{t_m\}$ are fixed,  and we denote them $b$ and $i_1,\ldots,i_b$.

The tree $T_{w_j}$ has finitely many
edges. By traveling along the path $p_{w_j}$ we pair subsegments of the path that cover the same edge in $T_{w_j}$ with
opposite orientations. We start from the basepoint in $T_{w_j}$. Given an edge in $T_{w_j}$ that is covered more than once by
the path $p_{w_j}$, we pair the first and the second subsegments of $p_{w_j}$ that pass through this edge. Note that since
$T_{w_j}$ is a tree, these
two subsegments have to be of opposite orientations. If $p_{w_j}$ passes more than 3 times through an edge in $T_{w_j}$,
we further pair the 3rd and 4th subsegments of $p_{w_j}$ that pass through such an edge and so on.

The number of pairs of subsegments of $p_{w_j}$ that we obtained depends only on the original tree $T_{w_j}$, and the path
$p_{w_j}$. The paths $p_{z^t_j}$ are the same as the path $p_{w_j}$ (as paths in the real tree $Y$) for every index
$t$. Our goal is to show that there exists a subsequence of the given sequence,
so that for every $m$ in the subsequence, the subwords of $z^{t_m}_j$ that are associated with subsegments of
the path $p_{z^{t_m}_j}$ that were paired together, can be replaced by eliminating distinct pair of subwords (part (2) in the
statement of the proposition), to subwords of uniformly bounded word length,
 and so that by eliminating these
distinct pairs the element $z^{t_m}_j$ is
replaced by an element $\hat z^{t_m}_j$ that represents the same element in the surface group $Q$  and satisfies part (3) of the
proposition for some constant $C$.

Let $e$ be an edge in $T_{w_j}$, and let $p_1$ and $p_2$ be subpaths of $p_{w_j}$ that were paired together and are
supported on $e$. Let $p^{t_m}_1$ and $p^{t_m}_2$ be the corresponding subpaths of $p_{z^{t_m}_j}$. By construction, the word
lengths of the paths
$p^{t_m}_1$ and $p^{t_m}_2$ are not bounded.

For each index $m$, we set $sc_m$ to be the maximal length of the elements $u^{t_m}_{i_1},\ldots,u^{t_m}_{i_b}$ (note that
$sc_m$ approaches 0 when $m$ grows to infinity). We set $Y_m$ to be the real tree $Y$ equipped with the metric that is obtained
from the metric on $Y$ by dividing it by $sc_m$. $Q$ acts isometrically on $Y$, so it naturally acts isometrically
on $Y_m$. By proposition 4.4 the length of each of the elements $u^{t_m}_{i_1},\ldots,u^{t_m}_{i_b}$
is bounded above by 1, and below by $d_1$ (the constant $d_1$ is defined in part (6) of proposition 4.4).

From the actions of $Q$ on the trees $Y_m$ it is possible to extract a subsequence (still denoted $\{t_m\}$) that
converges into a (non-trivial) action of $Q$ on some real tree $Y_{\infty}$. $Q$ acts faithfully on $Y_{\infty}$ and the length of
each of the elements $u_{i_1},\ldots,u_{i_b}$ on $Y_{\infty}$
is bounded below by $d_1$ and above by $1$. Since the action of $Q$ on $Y_{\infty}$ is faithful, and $Q$ is a surface group,
the action of $Q$ on $Y_{\infty}$ contains only IET and discrete components. Since we assumed a global bound on the periodicity of the
images: $h_n(s_j)$ for every $n$, and every $j$, $1 \leq j \leq r$, and the elements $s_1,\ldots,s_r$ generate $Q$, the
action of $Q$ on the limit tree $Y_{\infty}$ contains no discrete components. In particular the stabilizer of every
non-degenerate segment in $Y_{\infty}$ has to be trivial.

Let $e$ be an edge in $T_{w_j}$ (in the tree $Y$) with subpaths $p_1$ and $p_2$ that were paired together, and let $p^{t_m}_1$ and $p^{t_m}_2$ be
the corresponding subpaths of $p_{\hat z^{t_m}_j}$. As the lengths of the elements $u^{t_m}_{i_1},\ldots,u^{t_m}_{i_b}$ (when
acting on the tree
$Y$) approaches 0 with $m$, the combinatorial lengths of the paths
$p^{t_m}_1$ and $p^{t_m}_2$ grows to $\infty$ with $m$.

Given the appearance of a generator $u^{t_m}_{i_a}$, $a=1,\ldots,b$ in $p^{t_m}_1$, and
the appearance of a generator $(u^{t_m}_{i_{\hat a}})^{-1}$ in $p^{t_m}_2$ that overlap in a non-degenerate segment, we define their
$dual$ $position$ to be the subsegment in which they overlap, and its image in the segments that are associated with both
$u^{t_m}_{i_a}$  and
$u^{t_m}_{i_{\hat a}}$. Although the combinatorial lengths of the paths
$p^{t_m}_1$ and $p^{t_m}_2$ grows to infinity with $m$, the number of possible dual positions between the appearances of the
various generators remain bounded.

\vglue 1pc
\proclaim{Lemma 4.6} With the assumptions of proposition 4.5 (in particular, the bounded periodicity of the images
$h_n(s_j)$), there exists a global bound $R$,
such that for every index $m$ (from the chosen convergent subsequence),
the number of dual positions of overlapping subwords $u^{t_m}_{i_a}$  in $p^{t_m}_1$ and $(u^{t_m}_{i_{\hat a}})^{-1}$ in $p^{t_m}_2$
($a,\hat a=1,\ldots,b$) is
bounded by $R$.
\endproclaim

\nfp Suppose that the number of dual positions is not universally bounded. Then there exists a subsequence (still denoted
$\{t_m\}$), pair of indices $i_a$ and $i_{\hat a}$, $a,\hat a=1,\ldots,b$, such that the number of dual positions of
$u^{t_m}_{i_a}$  in $p^{t_m}_1$ and $(u^{t_m}_{i_{\hat a}})^{-1}$ in $p^{t_m}_2$ is bigger than $m$.

In that case, a simple pigeon hole argument implies that in the action of $u^{t_m}_{i_a}$ on $Y_m$,
there exists a subinterval $J_m \subset [y_m,u^{t_m}_{i_a}(y_m)]$ (where $y_m$ is the base point of $Y_m$) so that:
\roster
\item"{(i)}" the length of $J_m$ satisfies: $1 \geq length(J_m) > \epsilon_0 > 0$ for every $m$.

\item"{(ii)}" there exists some non-trivial element $q_m \in Q$, for which $q_m(J_m)$ overlaps with $J_m$ in an interval of length
that is at least $(1- \delta_m) \cdot length(J_m)$, and the sequence $\delta_m$ approaches 0 with $m$.
\endroster

Parts
(i) and (ii) clearly imply that the periodicity of the elements $h_n(s_j)$ can not be globally bounded,
a contradiction to the assumptions of propositions 4.4 and 4.5.

\hfill$\square$

Lemma 4.6 proves that in the chosen subsequence of words, $z^{t_m}_j$,
the number of possible dual positions between the appearances of the possible pairs of
generators, $u^{t_m}_{i_{a}}$ and $(u^{t_m}_{i_{\hat a}})^{-1}$, in overlapping paths, $p^{t_m}_1$ and $p^{t_m}_2$, remains bounded.
Given a word $z^{t_m}_j$, its associated
path $p_{z^{t_m}_j}$, and two of its overlapping subpaths $p^{t_m}_1$ and $p^{t_m}_2$, we use the bound on
the number of dual positions to replace $z^{t_m}_j$ with a
shorter word $\tilde z^{t_m}_j$, so that:
$$z^{t_m}_j(u^{t_m}_1,\ldots,u^{t_m}_{g_{t_m}}) \, = \,
\tilde z^{t_m}_j(u^{t_m}_1,\ldots,u^{t_m}_{g_{t_m}}).$$
When two appearances of a pair $u^{t_m}_{i_a}$ in $p^{t_m}_1$
and $(u^{t_m}_{i_{\hat a}})^{-1}$ in $p^{t_m}_2$ belong to the same dual position, we can
trim the paths $p^{t_m}_1$ and $p^{t_m}_2$,
by erasing the (identical) subpaths between these two appearances from $p^{t_m}_1$ and $p^{t_m}_2$. We set $\tilde z^{t_m}_j$ to be
the word that is obtained from $z^{t_m}_j$ after erasing the corresponding identical subpaths. We further
 repeat this erasing for all the
appearances of repeating pairs in the same dual position in all the (paired) overlapping subpaths $p^{t_m}_1$ and $p^{t_m}_2$
along the path
$p_{z^{t_m}_j}$. We denote the words that are obtained after this erasing, $\hat z^{t_m}_j$.
Since by lemma 4.6 there are at most $R$ dual positions for every pair of generators, and there are at most $(2f)^2$
pairs of generators, overlapping subpaths in the path $p_{\hat z^{t_m}_j}$ have combinatorial length bounded by $R \cdot (2f)^2$.
As the number of overlapping subpaths in $p_{\hat z^{t_m}_j}$ is bounded by  the number of overlapping subpaths in $p_{z^{t_m}_j}$,
which is identical to the number of overlapping subpaths in $p_{w_j}$, the total combinatorial length of the overlapping
subpaths in $p_{\hat z^{t_m}_j}$ is universally bounded. Therefore, the words $\hat z^{t_m}_j$ satisfy the conclusions
of proposition 4.5.

\hfill$\square$

By part (4) of proposition 4.4 the tuples $u^{t_m}_1,\ldots,u^{t_m}_{g_{t_m}}$ belong to finitely many isomorphism classes. By proposition
4.5 the total combinatorial lengths of the overlapping subpaths in the paths $p_{\hat z^{t_m}_j}$  are universally bounded. Hence,
we can pass to a further subsequence (that we still denote $\{t_m\}$) for which the isomorphism class of the tuples
$u^{t_m}_1,\ldots,u^{t_m}_{g_{t_m}}$ is identical. In particular $g_{t_m}=g$ is fixed, and the overlapping subpaths in the paths
$p_{\hat z^{t_m}_j}$  represent the same words in the generators $u^{t_m}_1,\ldots,u^{t_m}_{g_{t_m}}$ (in correspondence).

For a fixed index $t_{m_0}$ and each index $t_m$ from this subsequence we set:
$$\hat s^{t_m}_j= \hat z^{t_{m_0}}_j(u^{t_m}_1,\ldots,u^{t_m}_g).$$
Since the tuples $u^{t_m}_1,\ldots,u^{t_m}_j$ belong to the same isomorphism class, the tuples $\hat s^{t_m}_1,\ldots,\hat s^{t_m}_r$
belong to the same isomorphism class as the tuple $s_1,\ldots,s_r$ as $s_j=\hat s^{t_{m_0}}_j$. Since the overlapping subpaths in
the paths
$p_{\hat z^{t_m}_j}$  represent the same words in the generators $u^{t_m}_1,\ldots,u^{t_m}_g$, the paths
$[y_0,\hat s^{t_m}_j]$ are all positively oriented. Since the lengths of the elements $u^{t_m}_i$ approaches 0 when $m$ grows to
$\infty$, so are the lengths of the elements $\hat s^{t_m}_j$. Therefore, if we denote the automorphism that maps the tuple
$s_1,\ldots,s_r$ to $\hat s^{t_m}_1,\ldots,\hat s^{t_m}_r$ by $\varphi_m$, then a subsequence of the  automorphisms
$\{\varphi_m\}$ that we denote $\{\varphi_{\ell}\}$ satisfy the conclusions of theorem 4.2.



\medskip
Suppose that the periodicity of the images $h_n(s_j)$, for every $n$ and $j$, $1 \leq j \leq r$, is globally bounded,
and not all the endpoints of the intervals
belong to the same orbit.

In that case the group that is generated by the elements $v_1,\ldots,v_f$ is a limit group $L$  that
  strictly contains the surface group $Q$, and the elements $u^t_1,\ldots,u^t_{g_t}$ generate the limit group $L$ that is generated by
$v_1,\ldots,v_f$. However, for each index $t$, we obtained an IET transformation, where the generators of this IET transformation generate
the surface group $Q$, and since the combinatorics of these IET transformations belong to finitely many classes, there are finitely many
presentations of the surface group $Q$ that are associated with the sequence of tuples of elements
$u^t_1,\ldots,u^t_{g_t}$.

Hence, the arguments in the proofs of lemmas 4.3-4.6 remain valid in case the beginnings and ends of the subintervals in the IET transformation belong
to finitely many distinct orbits, and so do
the conclusion of theorem 4.2, so far under the
bounded periodicity assumption.



\medskip
To prove theorem 4.2 omitting the bounded periodicity assumption, we use the same strategy of proof, but modify the statements
and the arguments to include long periodic subwords, or in the limit, to include non-degenerate segments with non-trivial
stabilizers (that are contained in the discrete or simplicial part of the limit action).

As in the bounded periodicity case, for presentation purposes we first assume that the beginnings and ends of the subintervals in the
original IET transformation belong to the same orbit under the action of the surface group $Q$. We start with a generalization of
proposition 4.4.

Recall that we started the proof of theorem 4.2 with an iterative process of Dehn twists.
$v^t_1,\ldots,v^t_f$ are the (positive) elements between consecutive initial and final points of the bases that generate the
IET after $t$ Dehn twist iterations. The aim of proposition 4.4 was to replace these generators by a new (possibly larger set of)
generators so that the ratios between their lengths is globally bounded.

In the sequel a $semiperiodic$ element $u$ is an element of the form $u=\alpha w^p$, where $p>10$, $w$ is positively oriented, and  $\alpha$ is a suffix of $w$.

\vglue 1pc
\proclaim{Proposition 4.7}
For every $t \geq 1$, the elements
$v^t_1,\ldots,v^t_f \in Q$ can be replaced by elements $u^t_1,\ldots,u^t_{g_t} \in Q$ that satisfy properties (1)-(5)
in proposition 4.4. Properties (6) and (7) in proposition 4.4 are replaced by the following properties:

\roster
\item"{(6)}" there exist a pair of indices $1 \leq \ell(t) \leq b(t) \leq g_t \leq e(f)$, for which $b(t)>1$. The elements $u^t_1,\ldots,u^t_{\ell(t)}$
are semi-periodic,  the elements $u^t_{\ell(t)+1},\ldots,u^t_{b(t)}$ are short, and the elements $u^t_{b(t)+1},\ldots,u^t_{g_t}$ are secondary short.
Note that some of these sets may be empty, but either the set of semiperiodic elements or the set of short elements must be non-empty.

\item"{(7)}"  for some $d_2>1$ and every $1 \leq i_1 \leq b(t)$ and $b(t)<i_2 \leq g_t$,
$$ 10g_t \cdot d_2 \cdot length (u^t_{i_2})  \, \leq \, length (u^t_{i_1}) $$

\item"{(8)}" there exists some positive $d_1$ such that for every  $\ell(t)+1 \leq i_1,i_2 \leq b(t)$:
$$d_1 \cdot length (u^t_{i_1})  \, \leq \, length (u^t_{i_2}) \, \leq \,
 d_2 \cdot length (u^t_{i_1})$$
For some $d_3>1$ and every $1 \leq i_1 \leq \ell(t)$ and $\ell(t)+1 \leq i_2 \leq b(t)$:
$$d_3 \cdot length (u^t_{i_2})  \, \leq \, length (u^t_{i_1}) $$

\item"{(9)}" For each $t$, and every index $i$, $1 \leq i \leq \ell(t)$, there exist distinct indices
$\ell(t)+1  \leq j_1,\ldots,j_{e_i(t)} \leq g_t$
(that depend on $t$),
such that: $w_i=u^t_{j_1} \ldots u^t_{j_{e_m(t)}}$, and
$u^t_{i}=\alpha_i w_i^{p_i}$ where $\alpha_i$ is a suffix of $w_i$.

\item"{(10)}" for each index $t$, in the words $z^t_j$, $1 \leq j \leq r$, in a bounded distance (where the bound depends only on $f$)
either before or after each occurrence of an
element $u^t_{i_2}$, for $b(t)+1 \leq i_2 \leq g_t$, appears one of the elements $u^t_{i_1}$, for some $1 \leq i_1 \leq b(t)$.
\endroster
\endproclaim

\nfp We use the same iterative procedure that was used to prove proposition 4.4. We divide the elements $v^t_1,\ldots,v^t_f$
into finitely many sets according to their length.
We order the elements from the longest to the shortest. We place a separator between consecutive sets
whenever there is a pair of consecutive elements (ordered according to length) that satisfy:
$length(v^t_i) \geq c_1(f) \cdot length(v^t_{i+1})$
where $c_1(f)=4f$.
Note that the constant $c_1(f)$ in the current proposition depends only
on the number of generators $f$, and not on the periodicity bound, that we didn't assume exists.

We use the same iterative procedure that we used in the proof of proposition 4.4, and we refer to its various cases
and notions as it appears in the proof of proposition 4.4.
In cases (1)-(3) of this procedure we proceed precisely as in proposition 4.4 (the bounded periodicity case).

In case (4) of the procedure, we assumed that a long and a short base are adjacent to $pv$, the terminal point of
the interval that supports the IET, and that the length of the long base that is adjacent to $pv$,
that we denote $b_1$,
is smaller than the sum of the short bases that are adjacent to $pv$ and the first long base that is adjacent to these short
bases. We set $b_2$ to be the first long base that is adjacent to the short bases that are adjacent to $pv$.

If $b_1$ and $b_2$ are not a pair of bases we do what we did in proposition 4.4 (essentially what we did in
cases (1) and (2)).
Suppose that $b_1$ and $b_2$ are a pair of bases. In this case $b_1$ and $b_2$ are a pair of long bases,
such that $b_1$ is obtained by a shift of $b_2$ by a short (or secondary short) subinterval. Hence, $b_1$ (and $b_2$)
can be presented as a short word times
a high power of another short word, precisely as described in case (9). In this case we add a short or a secondary short pair of bases at the ends of
$b_1$ and $b_2$.
The pair that remains from $b_1$ and $b_2$ satisfies the properties in case (9).  We continue with the interval that complements the support of $b_1$ together
with the support of the short (and secondary short) bases that are adjacent to $pv$.

The number of pairs of long bases is reduced by 1, and we added a pair of short or secondary short bases.
The number of active short and secondary short bases, i.e., the number of short and secondary short bases that
participate in the rest of the procedure is either unchanged or is reduced by 1.

In case (5) let $b_1$ and $b_2$ be the two long bases that are adjacent to
the two sequences of short and secondary short bases that are adjacent to $pv$. If $b_1$ and $b_2$ is not a pair of bases we do what we did in the
proof of proposition 4.4. If $b_1$ and $b_2$ are a pair of bases we do what we did in case (4) when $b_1$ and $b_2$
are a pair of bases. Once again the number of pairs of long bases is reduced by 1, and the number of  short and secondary short bases that participate in the rest
of the procedure is either unchanged or is reduced by 1.

In case (6) we act precisely as we did in the proof of proposition 4.4 if the two long bases are not a pair of bases, and as
we did in part (5) above if the long bases are a pair of bases.

As in proposition 4.4,
when we run this procedure until there are no long elements. i.e., until all the elements are either short,
secondary short, or a product of short or secondary short elements  and a high power of a product of short and secondary short elements
(part (9) of the proposition),
 the number of bases can grow to at most $2f$, and therefore the number of elements can grow
to at most $4f$.

The conclusions of the proposition follow from the termination of the iterative procedure, precisely as in the
proof of proposition 4.4.

\hfill$\square$

Proposition 4.7 replaces proposition 4.4 in the general case (i.e., when there is no periodicity assumption).
To obtain the same conclusions as in proposition 4.5, we further modify the tuples, $u^t_1,\ldots,u^t_{g_t}$.

\vglue 1pc
\proclaim{Proposition 4.8}
With the notation of proposition 4.7, it is possible to further modify the tuple of elements
$u^t_1,\ldots,u^t_{g_t}$, by performing Dehn twists on some of the semiperiodic elements (the elements
$u^t_1,\ldots,u^t_{\ell(t)}$ that satisfy part (9) in proposition 4.7),
so that there exists a constant $C>0$, for which for the modified tuples, that we still denote: $u^t_1,\ldots,u^t_{g_t}$,
for every index $t \geq 1$ the words $z^t_j$, $1 \leq j \leq r$,
can be replaced by words: $\hat z^t_j$ that satisfy properties (1)-(3)
in proposition 4.5.
\endproclaim

\nfp Suppose that such a constant $C$ does not exist (for any possible application of Dehn twists on the semiperiodic
elements in the tuples: $u^t_1,\ldots,u^t_{g_t}$). Then for every positive integer $m$,
there exists an index $t_m>1$, so that for every possible choice of (applications of) Dehn twists to the
semiperiodic elements in the tuple, $u^{t_m}_1,\ldots,u^{t_m}_{g_{t_m}}$, at least
one of the  words $\hat z^{t_m}_{j}$ that satisfy
(3) is false for the constant $C=m$.

By passing to a subsequence, we may assume that for
the subsequence the integer $g_{t_m}$, that counts the number of generators of the semigroup that is constructed in step
$t_m$ according to the procedure that is described in proposition 4.7 is fixed and we denote it $g$. We may also
assume that the integers,  $\ell(t_m)$ and
$b(t_m)$, that are associated with
the sets of generators of the semigroups that are constructed according to proposition 4.7 along the subsequence
$\{t_m\}$ are fixed,  and we denote them $\ell$ and $b$.

For each index $t_m$, we denote by $length_m$, the minimal length of a short element (i.e., the elements
$u^{t_m}_{\ell+1},\ldots,u^{t_m}_b$). For each semiperiodic element $u^{t_m}_1,\ldots,u^{t_m}_{\ell}$ (that are necessarily long)  we denote the
length of its
period by $lper^i_m$.

For each index $m$, and every $i$, $1 \leq i \leq \ell$, we look at the ratios: $\frac {lper^i_m} {length_m}$.
We can pass to a subsequence of the indices $m$, for which (up to a change of order of indices):
 $0 < \epsilon < \frac {lper^i_m} {length_m}$ for some positive $\epsilon >0$, and $i=1,\ldots, \ell'$. And for every
$i$, $\ell' < i \leq \ell$,
the ratios $\frac {lper^i_m} {length_m}$ approach 0.
We perform Dehn twists along the semiperiodic elements $u^{t_m}_1,\ldots,u^{t_m}_{\ell'}$
(from the subsequence of indices $m$), so that all these semiperiodic elements have lengths bounded (from below and above) by
constants times the minimal length of a short element. These elements will be treated as short elements and not as semiperiodic elements
in the sequel.

First, suppose that $\ell'=\ell$, i.e., that there exists an $\epsilon>0$, such that for every
$i$, $1 \leq i \leq \ell$,
 $0 < \epsilon < \frac {lper^i_m} {length_m}$. In that case, after applying Dehn twists to the semiperiodic
elements, all the elements $u^{t_m}_1,\ldots,u^{t_m}_g$ are either short or secondary short. By the argument that was
used to prove proposition 4.5, either:
\roster
\item"{(i)}"  the number of dual positions of the different elements in cancelling subpaths along the paths $p_{z^t_j}$ is globally
bounded (for the entire subsequence $\{t_m\}$).

\item"{(ii)}" there exists a subsequence (still denoted $\{t_m\}$), and a fixed
positive word in some of the short elements and possibly some of the secondary short
elements, which is a periodic word, and the ratio between the length of the period and the length
of the element that is represented by the positive word approaches 0.
\endroster
Part (ii) implies that some positive word in short and possibly secondary short elements is periodic, and the length of the period divided by the minimal length
of a short element approaches 0. The existence of such a periodic word implies that  the surface group $Q$ contains a free abelian group of rank at least 2,
or that
$Q$ is freely decomposable, and we get a contradiction. Hence, in case $\ell=\ell'$ there is a global bound
on the number of dual positions of the different elements for the entire subsequence $\{t_m\}$, and the conclusion of
the proposition follows by the same argument that was used to prove proposition 4.5.
The same argument remains valid if $\ell'<\ell$ but the lengths of the semiperiodic elements
$u^{t_m}_1,\ldots,u^{t_m}_{\ell}$ can be bounded from above by a constant times the length of a short element.

Suppose that there exists a subsequence of indices (still denoted $\{t_m\}$), for which along a paired subpaths
$p_1$ and $p_2$ (of opposite orientations), at least one of the appearances of a semiperiodic element $u^{t_m}_i$, $\ell' < i \leq \ell$
(along $p_1$ or $p_2$), overlaps
with an unbounded number of elements (along $p_2$ or $p_1$ in correspondence). In that case, like in part (ii)
in the case in which the
lengths of the semiperiodic elements are  bounded from above by a (global) constant times the length of a short element, we get:
\roster
\item"{(ii')}" there exists a subsequence (still denoted $\{t_m\}$), and a fixed
positive word in either short elements or semiperiodic elements and possibly some secondary short
elements, which is a periodic word, and the ratio between the length of the period and the length
of the element in $Q$ that is represented by the positive word approaches 0.
\endroster

Hence,  either the surface group $Q$ contains a free abelian group of rank at least 2, or
$Q$ is freely decomposable, and we get a contradiction. Therefore, there exists a global bound on the number of elements
that overlap with a semiperiodic element that appears along a paired subpaths $p_1$ and $p_2$.

In this last case, if there is no global bound on the number of dual positions between overlapping elements, once again  the surface group
$Q$ must contain either a free abelian group of rank at least 2, or it is freely decomposable, and we get a contradiction.

Finally, in all the cases we obtained a global bound on the number of dual positions of the different elements,
so the proposition follows by the same argument that was used to prove proposition 4.5.

\hfill$\square$

Given propositions 4.7 and 4.8, that generalize propositions 4.4 and 4.5, the rest of the proof
of theorem 4.2 follows precisely as in the bounded periodicity case. Hence, theorem 4.2 holds
without the assumption of bounded periodicity, but still with the assumption that the beginnings and ends
of the subintervals in the original IET transformation belong to a single orbit under the action of the surface group $Q$.

Suppose that the beginnings and the ends of the subintervals in the original IET transformation do not belong to the same
orbit under the action of $Q$ and there is no uniform bound on the periodicity of the images $h_n(s_i)$. The proof of proposition 4.7 remains valid in this case,
since we did not use the unique orbit assumption of the beginnings and ends of the subintervals in the original IET transformation.

The proof of proposition 4.8 uses the fact that the subgroup generated by $u^t_1,\ldots,u^t_{g_t}$ is the surface group $Q$, and a surface group
is freely indecomposable and contains no non-cyclic free abelian group. As in the bounded periodicity case, if the beginnings and ends of the
subintervals in the initial IET transformation belong to finitely many orbits, the groups that is generated by $v_1,\ldots,v_f$, and hence the group
that is generated by the sequences of elements, $u^t_1,\ldots,u^t_{g_t}$, contains the surface group $Q$.

However, there are only finitely many distinct  presentations of the surface subgroup $Q$ that are associated with the IET transformations that are associated with
the sequences of elements $u^t_1,\ldots,u^t_{g_t}$. Furthermore, the arguments that are used in the proof of proposition 4.8 construct
a free splitting of the subgroup $Q$ or show that $Q$ contains a free non-cyclic subgroup in case the conclusions of
the proposition do not hold.

If the conclusions of proposition 4.8 do not hold in case the beginnings and ends of the subintervals in the original IET transformation
belong to more than a single orbit, then the same arguments that were used in the proof of the proposition show that the surface subgroup $Q$,
and not only the group that is generated by $u^t_1,\ldots,u^t_{g_t}$, admits a free splitting or contains a non-cyclic free abelian group.

Therefore, the proof of proposition 4.8 remains valid in case the beginnings and ends of the original IET transformation belong to more
than a single orbit. This completes the proof of theorem 4.2.



\hfill$\square$

\medskip
So far we studied axial and IET components in a limit action of a pair on a real tree.
The next theorem constructs a positive free basis in case the ambient group is free (non-abelian), and its action on
the limiting tree is free and is either discrete or non-geometric. The remaining case, which is the Levitt action of a free
group, will be analyzed in the sequel.
The existence of such a positive free basis is crucial in our general approach
to the structure of
varieties over a free semigroup.

\vglue 1pc
\proclaim{Theorem 4.9} Let $F$ be a free group, and let $\{h_n:(S,F) \to (FS_k,F_k)\}$ be a sequence of homomorphisms
of pairs that converges into a
free action of the
limit pair $(S,F)$ on a real tree $Y$. Suppose that the action of $F$ on the limit tree $Y$ does not contain any
Levitt components.

Then there exists a directed finite graph $\Theta$ with the following properties:
\roster
\item"{(1)}" $\Theta$ contains a base point, and the free group $F$
is identified with the fundamental group of the graph $\Theta$. With each (positively oriented) edge in $\Theta$
we associate a label,
and each of the given set of generators of the subsemigroup $S$ corresponds to a positive path in $\Theta$, that starts and
ends at the base point, and  can be expressed
by a positive word in the labels that are associated with its edges.

\item"{(2)}" there exists an index $n_0$ such that for every index $n>n_0$ the homomorphisms $\{h_n\}$ are encoded by the graph $\Theta$.
 i.e., each of these  homomorphisms  is obtained by
substituting elements from  the free semigroup $FS_k$ to the elements that are associated with the various
positive edges in $\Theta$
(these substitutions define a homomorphism from the fundamental group $F$ of $\Theta$ to the coefficient free group $F_k$,
and this homomorphism is precisely the homomorphism from the subsequence).
\endroster
\endproclaim

\nfp If the action of $F$ on the limit tree $Y$ is discrete (simplicial) then the conclusion of the theorem is immediate. Hence, we may assume that
$Y$ contains non-discrete parts and no Levitt components. In particular, the action is not geometric
(Levitt (or thin))  and
non-geometric actions on a real
tree are defined in [Be-Fe1]).

Note that since the tree $Y$ is the limit of a sequence of homomorphisms $\{h_n:(S,L) \to (FS_k,F_k)\}$ it inherits an equivariant orientation
from that of the pair $(FS_k,F_k)$. A segment that connects the base point in $Y$ to its image under the action of one of the generators in $S$ is
positively oriented. The finite union of orbits of these segments cover the entire tree $Y$. Hence, every segment in $Y$ can be cut into finitely
many segments so that each of the obtained segments is either positively or negatively oriented.

Suppose that the action of $F$ on $Y$ is free and contains only discrete and non-geometric components in
the sense of Bestvina-Feighn  (i.e., there are  no Levitt components). In that case every resolution
(in the sense of [Be-Fe1]) of the action
of $F$ on the real tree $Y$ is discrete (simplicial).

Suppose that $F \times T \to T$ is such a resolution. In particular,  the action of $F$
on $T$ is free and discrete, and there exists an $F$-equivariant map from $T$ onto $Y$. By
the construction of a resolution (see [Be-Fe1]), for a given finite subtree of $U \subset Y$,
it is possible to construct a resolution so that a lift of $U$ in the resolution is mapped isometrically onto $U$ in $T$.

Hence, we can construct a  resolution
so that
the (finitely many) segments that connect between the base point in $T$ and the
images (in $T$) of the base point under the action of the given (finite) set
of generators of the semigroup $S$ are embedded isometrically into the limit tree $Y$.
Since the finite union of the orbits of these segments in both trees $T$ and $Y$ cover
these trees, and each of these segments is positively oriented,
the (equivariant) orientation of segments in $Y$ lifts to an equivariant orientation of segments in
$T$ (i.e., every segment in $T$ can be cut into finitely many segments, and each of the obtained segments is either positively or negatively
oriented).

The action of $F$ on $T$ is free and discrete. Hence, by Bass-Serre theory it is
possible to associate with this action a finite graph of groups with a basepoint. We
denote this graph $\Theta$. Since the action is free, $\Theta$ has trivial vertex and edge
stabilizers, and its fundamental group is identified with $F$. Furthermore, the segments in $T$
are oriented equivariantly, and every edge in $\Theta$ is contained in (an orbit of) a segment that
is associated with one of the given generators of the semigroup $S$, hence, the orientation of segments in $T$
gives an orientation of
the edges in $\Theta$. Finally, since the segments that connect the base point to the images of the base
point under the action of the given set of generators of $S$ in $T$ are all positively oriented, the
loops in $\Theta$ that correspond to elements in the semigroup $S$ are all positively oriented, so
$\Theta$ satisfies the properties that are listed in part (1).

With each positive segment in the graph $\Theta$ we associate a label. Every substitution of values from
the free semigroup $FS_k$ to these labels gives a homomorphism of pairs from $(S,L)$ to $(FS_k,F_k)$.
Since the action of $F$ on $T$ resolves the action of $F$ on $Y$, and the action of $F$ on
$T$ is
discrete and free,
given a sequence of homomorphisms that converges into the action of $F$ on the tree $Y$, there is an index
$n_0$,
for which for every $n>n_0$ positively oriented segments in $\Theta$ are mapped to elements in $FS_k$.
Therefore, these homomorphisms  are obtained from $\Theta$ by substituting
these values to the labels that are associated with positively oriented edges, and we get part (2).

\hfill$\square$

To conclude this section we need to analyze Levitt components, that play an essential role in our analysis
of homomorphisms of pairs. The analysis of Levitt components that we use, is based on our analysis of
IET components, as it appears in the proof of theorem 4.2.

In analyzing Levitt components we don't use automorphisms of the ambient free
group that is associated with a Levitt component to shorten homomorphisms, but we rather embed this free group in a bigger
residually free group, and use isomorphisms between  these bigger groups to shorten homomorphisms.
i.e., we use isomorphisms between  residually free groups
that are associated with generalized equations in Razborov's work, residually free groups that appear along  Razborov fundamental sequences (see [Ra] or [Ca-Ka]),
to shorten homomorphisms that converge into an indecomposable free action of a free group with Levitt's dynamics.

In the bounded periodicity case, i.e., in case there is a uniform bound on the periodicity of the images of the fixed generators
of the given semigroup under the set of homomorphisms that
we consider, we argue that after passing to a subsequence these automorphisms of an ambient residually free group $RF$ restrict to automorphisms of
the original free group $F$ that is a subgroup of $RF$.

As we did in the IET case (theorem 4.2), we first prove the theorem in case there is a bound on the periodicity of the images of the
fixed set of generators of the semigroup $S$ that we are
looking at, and then discuss the general case.

\vglue 1pc
\proclaim{Theorem 4.10} Let $F$ be a f.g.\ free group, and let  $(S,F)$ be a pair where $S$ is a f.g.\
subsemigroup that generates $F$
as a group.
Let $\{h_n\}$ be a sequence of pair homomorphisms of $(S,F)$ into $(FS_k,F_k)$ that converges into an indecomposable free
action of $F$
on a real tree $Y$  of Levitt (thin)  type (for the notion of an indecomposable action see definition 1.17 in [Gu]).


Suppose that the periodicity of the images of the generators of $F$ under the homomorphisms $\{h_n\}$ are bounded by some
constant $c_p$.
Then the conclusions of theorem 4.2 (for the IET case) are valid for the Levitt case. There exists a sequence of automorphisms
$\{\varphi_{\ell} \in Aut(F)\}_{\ell=1}^{\infty}$ with the following properties:
\roster
\item"{(1)}" for each index $\ell$, there exists $n_{\ell}>0$, so that for every $n>n_{\ell}$, and all the images of the generators of
the subsemigroup $S$, $s_1,\ldots,s_r$, $h_n \circ \varphi_{\ell}(s_j) \in FS_k$.
(note that the automorphisms $\varphi_{\ell}$ do
not preserve  $\nu(S)$ in  general, but the images of $S$ under the twisted homomorphisms $h_n \circ \varphi_{\ell}$ remain
in the standard positive cone $FS_k$ in $F_k$).

\item"{(2)}" for every $n>n_{\ell}$:
$$\frac {d_T(h_n \circ \varphi_{\ell}(s_j), id)}
{\max_{1 \leq j \leq r}
d_T(h_n(s_j),id.)} < \frac {1}{\ell}$$
where $T$ is a  Cayley graph of the coefficient free group $F_k$ with respect to a fixed set of generators.
\endroster
\endproclaim

\nfp The argument that we use is an adaptation to the Levitt case of the argument in the IET case (theorem 4.2).
Let $s_1,\ldots,s_r$ be the generators of the subsemigroup $S$ that generates the free group $F$ that acts freely
on the limit tree $Y$, and the limit action is of Levitt type.
Hence, there exists a finite subtree $R$ of $Y$, such that if one applies
the Rips machine to the action of the pseudogroup generated by the restrictions of the actions of the generators
$s_1,\ldots,s_r$ to the finite tree $R$, the output is
a Levitt type pseudogroup (see [Be-Fe1]).

A Levitt type pseudogroup that is based on an  oriented interval $I$, is generated by finitely many pairs
of bases, where the union of the supports of these bases is a finite union of subintervals of $I$.
Since it is
a Levitt type pseudogroup, there are subintervals in $I$ that are covered only once by the union of bases.
Since every segment in the tree $Y$ can be divided into finitely many oriented segments, and since a Levitt pseudogroup is mixing
in the sense of [Mo] (note that an indecomposable action is always mixing), we can assume that the Levitt pseudogroup is supported on an oriented subinterval $I$ of $Y$.

Levitt pseudogroups are analyzed in [Be-Fe1], [GLP] and in [Ra]. Given a Levitt pseudogroup that is  supported on some (subintervals of an)
oriented interval
$I$, there are some subintervals that are covered exactly once by the bases of the pseudogroup.
Starting with the original f.g.\ pseudogroup one applies to it a sequence of moves
(see [Be-Fe1]). There are finitely many subintervals that are covered only once. In each move, one cuts such a subinterval that is
contained in a base. Hence,
possibly cut the base that is supported on this subinterval (if the subinterval does not contain an endpoint of the base),
and cuts a corresponding subinterval from the paired base.

Following section 7 in [Be-Fe1], subintervals that are covered exactly once and removed along the process are divided into 3 classes:
\roster
\item"{(1)}" an isolated base,
that is a base that its interior cuts no other bases. In this case the number of pairs of bases after removing the isolated base
decreases by 1.

\item"{(2)}" a semi-isolated base, that is a subinterval that is covered only once, contained in a base and contains an endpoint of that base.
Removing such a subinterval does not increase the number of bases.

\item"{(3)}" a subinterval that is covered only once and is contained in the interior of a base. In this case the number of
pairs of bases after removing
the subinterval increases by 1, and the number of connected components of subintervals that are covered by bases
increases by at least 1.
\endroster

In addition to removing subintervals that are covered exactly once, one performs the following operation:
\roster
\item"{(4)}" if there exists a subinterval that supports exactly two bases, and the support of both is the entire subinterval,
one removes that subinterval and the bases that it supports. If the bases are paired, they are erased and the
 number of pairs of bases reduces
by 1. If these bases are not paired,  one further pairs the  the two bases
that were previously paired with the bases that were supported on that subinterval. In that case the number of paired bases is also reduces by 1.
\endroster

The Rips' machine applies moves of types (1) and (4) as long as possible (i.e., removes subintervals of type (1) and (4)). Then
applies moves of type (2) as long as possible. If there are no
more moves of types (1), (4) and (2), then one applies a move of type (3). Since we assumed that the action of the given group
$F$ on the tree $Y$ is of Levitt type, and the given pseudogroup generates the action, moves of type (1) can not occur (since if a
subinterval of type (1) exists, the action is not indecomposable). Furthermore,
in all moves of type (4), the two bases that are supported on that subinterval can not be a pair of bases (If there exists a subinterval that
supports exactly a pair of bases, the action of $F$ can not have a dense orbit, and hence can not be indecomposable).

By proposition 7.2 in
[Be-Fe1] moves of type (3) occur infinitely many times, i.e., a sequence of moves of type (4) and (2) has always a finite length.
The Rips' machine does not give a priority to the order of the subintervals of type (3) that are being treated, as long as any such
subinterval is treated after a finite time.

With the pseudogroup and the interval that supports it we can naturally associate a graph. The vertices in this graph are maximal
subintervals that are covered at least once by bases of the pseudogroup. Each such maximal subinterval starts and ends
 with an endpoint
of one of the bases that is not contained in the interior of another base. The edges in the graph are associated with
the pairs of bases of the pseudogroup. For each pair of bases there is an edge that connects between a maximal subinterval
that supports a base to a maximal subinterval (possibly the same one) that supports its paired base. Clearly the Euler characteristic
of the graph that is associated with the initial Levitt pseudogroup is negative and is bounded below by $1-b$,
where $b$ is the number of paired
bases.

The graph that is associated with a pseudogroup that is obtained from the original one after a sequence of moves of
types (2)-(4) can have only a bigger (smaller in absolute value) Euler characteristic. After performing a further finite sequence of all
the possible moves of types (2) and (4) (there are finitely many such by proposition 7.2 in [Be-Fe1]), the valence of each vertex
in the associated graph is at least 3, so at each step of the process in which all the moves of types (2) and (4) were
performed,
 the number of generators of the obtained pseudogroup is universally bounded (in terms of the original Euler characteristic).

The sum of the  lengths of the bases in the pseudogroups along the process that starts with a Levitt pseudogroup and applies moves of type
(2)-(4) approaches 0 (see proposition 8.12 in [Be-Fe1]). Since the number of bases is totally bounded, the (infinite) intersection of the unions of the
subintervals that support these bases
consists of finitely many points.

\vglue 1pc
\proclaim{Lemma 4.11} Let $U$ be a f.g.\ Levitt pseudogroup that is based on some oriented interval $I$ that generates
an indecomposable  Levitt type action of a free group $F$ on some real tree $Y$. Let $U_1,\ldots$ be
the pseudogroups that are obtained from $U$ by applying the moves (2)-(4), where in each step we apply a single move of type (3) (if possible)
and then all the possible moves of types (2) and (4). Then for some index $t_0$ and for all $t>t_0$:
\roster
\item"{(1)}" the Euler characteristics of the graphs that are associated with the pseudogroups $U_t$ are equal to the
Euler characteristic of the graph that is associated with $U_{t_0}$.

\item"{(2)}" there is no nontrivial word $w$ in the generators of the pseudogroup $U_t$, that is defined on a non-degenerate
subinterval of $I$, and acts trivially on that subinterval.

\item"{(3)}" only moves of types (3) and (4) are applied in the process at each step $t$.
\endroster
\endproclaim

\nfp For every index $t$, each of the the operations (2)-(4) do not reduce the absolute value of the Euler characteristic of the graph
that is associated with the pseudogroup $U_t$. i.e., these operations can not increase the absolute value of the Euler characteristic.
Hence, after finitely many steps the Euler characteristic stabilizes and part (1) follows.

Suppose that at some step $t$, there is a nontrivial word $w$ in the generators of $U_t$, that is defined on some non-degenerate
subsegment of $I$, and $w$ fixes that subinterval pointwise. In that case there exists a point $q \in I$, that is contained in the
interior of a subinterval $J$,
$J \subset I$, such that the interval J has a finite orbit under the action of the  finite set of prefixes of $w$.
We can further assume that all
the translates of $J$ in that finite orbit are disjoint.

Suppose that the disjoint translates of $J$ in that finite orbit are $J_0,\ldots,J_{\ell}$. All these subintervals are covered by at least
2 bases from $U_t$, and for each index $i$ there exists a pair of bases that map $J_i$ to $J_{i+1}$, $i=0,\ldots,\ell-1$, and $J_{\ell}$ to $J_0$.
The subinterval $J$ (or rather the subinterval to which it is going to move using moves of type (4)) will stay covered by at least
two bases along the process, unless at some step $t_1>t$, at least one of the subintervals in its orbit will be covered
twice by a base and its paired base, and these two paired bases have the same support (i.e., the map from that base to its paired base is the
identity map). Since the lengths of the bases approaches 0 when $t \to \infty$, such a step $t_1>t$ must exist.
In this case we can remove this pair of bases. Since the two paired bases that we remove can not be supported
on a subinterval that does not support another base, the removal of this pair increases the Euler characteristic of the graph that is
associated with the pseudogroup $U_{t_1}$. By part (1) such a reduction can occur only finitely many times, and it does not happen for
$t>t_0$, and part (2) follows.

By parts (1) and (2) there is a step $t_0$ such that for every $t>t_0$ the Euler characteristic of the graph that is associated with
the generating sets of the pseudogroups $U_t$ is constant, and there is no non-degenerate subsegment of the interval $I$ that
is fixed pointwise by a non-trivial word $w$ in the pseudogroup $U_t$. In particular, for all $t>t_0$ the pseudogroup $U_t$
can not contain a pair of bases that are supported on the same subinterval of $I$ (as otherwise this pair can be removed and increase
the Euler characteristic, or leave a subinterval of $I$ with trivial (simplicial) dynamics).

By the structure of the process, at step $t_0+1$ we start with  $U_{t_0}$ and apply to it all the possible
moves of types (2) and (4) (there are
finitely many such by proposition 7.2 in [Be-Fe1]). At this point the endpoints of all the bases in the obtained pseudogroups
are covered at least twice (otherwise a move of type (2) is still possible).

Now we apply move (3) along a subinterval $J_1 \subset I$ that is covered only once by a base $b_1$, and since it is a move
of type (3), $J_1$ is contained in the interior of the subinterval that supports the base $b_1$. The base $b_1$ is paired with a base
$b_2$, so by move (3) we cut $b_1$ along $J_1$ and $b_2$ along a (paired) subinterval $J_2 \subset I$, and increase the number of paired bases by 1.
We denoted the obtained pseudogroup $U_{t_0+1}$.

The subintervals $J_1$ and $J_2$ have to be disjoint. $J_1$ was covered exactly once before we erased the corresponding part of
$b_1$, so every point in $J_2$ was covered by at least two bases before we erased the corresponding part of $b_2$, since otherwise there
is a subinterval of $J_2$  on which the pseudogroup $U_{t_0}$ act discretely, and a Levitt pseudogroup is mixing.

The endpoints of all the bases that are not supported by $J_2$ are covered at least twice, since they were covered twice before the move
of type (3) was applied, and this move does not effect these endpoints. Suppose that there exists an endpoint of a base that is supported by $J_2$,
and this endpoint is supported by a single base in the pseudogroup $U_{t_0+1}$. Since every point in the subinterval $J_2$ is covered at least
once,  in case there is a point in $J_2$ (including its endpoints) that is covered exactly once, there must exist a point in $J_2$ (including its endpoints)
that is not contained in the interior of a base. Hence, the graph that is associated with $U_{t_0+1}$ contains at least 2 new vertices of valence at least
2 in addition to the vertices that existed in the graph that was associated with the pseudogroup before the move of type (3). Since the number
of paired bases was increased only by 1, the Euler characteristic of $U_{t_0+1}$ is bigger by at least 1 from that of $U_{t_0}$, and we got
a contradiction to part (1).

Therefore, all the endpoints of bases in $U_{t_0+1}$ are covered by at least 2 bases. Moves of type (4) do not change this property, so no move
of type (2) is required before a move of type (3) is applied. Hence, part (3) follows by a straightforward induction.

\hfill$\square$

At this point we change our point of view, and instead of the band complex that is associated with the collection of pairs of bases we look at the generalized
equation that is associated with them (see [Ma1], or [Ma2] or definitions 3.1 and 3.2 in [Ca-Ka]). Note that the band complex and the generalized equation that are
associated with a collection of pairs are analogous. i.e., from the band complex it is possible to construct the generalized equation and vice versa.
With the generalized equation Razborov associated a residually free
group. He starts with a group  that is generated by elements that are associated with subintervals between consecutive endpoints of bases, and the relations are
given by the pairs of bases, i.e., each pair of bases indicates that a certain word that is associated with one of the bases is equal to the word that
is associated with its dual base. Then Razborov looks at the residually free quotient of this group (i.e., one replaces the normal closure of the
relations with their radical). The definition appears in detail in definition 3.8 in [Ca-Ka].

Let $\Omega$ be the generalized equation and let $G_{R(\Omega)}$ be the residually free group that is associated with it. By construction the free
group $F$ that is associated with the Levitt component is embedded into $G_{R(\Omega)}$. We continue by performing the operations (3) and (4)
iteratively on the generalized equation $\Omega$, precisely as is done according to the Rips machine, to get a sequence of
generalized equations, $\Omega_1=\Omega$, $\Omega_2, \ldots$.

By construction, the sequence of generalized equations that we obtain factor through one of the (finitely many) Razborov's fundamental sequences
for a generalized equation (see [Ra] and [Ca-Ka] for the construction of Razborov's fundamental sequences).
As in the IET case, we denote the elements that are
associated with the subintervals between consecutive beginnings and  endpoints of bases in the Levitt pseudogroup by
$v_1^t,\ldots,v_{f_t}^t$, where $f_t$ is bounded
by twice the number of bases in the original Levitt pseudogroup.

For large enough $n$ the homomorphisms, $\{h_n:S \to FS_k\}$, extend to homomorphisms, $\{\hat h_n: G_{R(\Omega)} \to F_k\}$. Note that $G_{R(\Omega)}$ is the group that is generated
by the elements $v_1,\ldots,v_f$, the elements that are associated with the original pseudogroup that is associated with the original generalized equation $\Omega$.

For each index $t$, there exists some index $n_t$, such that for all $n>n_t$, $\hat h_n(v^i_t) \in FS_k$.
By the properties of Razborov's fundamental sequences for some index $t_1$ and for every $t>t_1$, the groups $G_{R(\Omega_t)}$ are all isomorphic,
and the sets of generators $v_1^t,\ldots,v_{f_t}^t$ belong to finitely many isomorphism classes of generating sets of these isomorphic groups.

Let $t_1$ be the index for which for every $t>t_1$ only moves (3) and (4) are applied along the iterative process.
Since the process is using only operations (3) and (4), at any step $t>t_1$,
for every $t_3>t_2 > t_1$, the semigroup that is generated by
$v_1^{t_2},\ldots,v_{f_{t_2}}^{t_2}$ is contained in the semigroup that is generated by
$v_1^{t_3},\ldots,v_{f_{t_3}}^{t_3}$.

These properties of the sequence of generalized equations, $\Omega_1,\ldots$ and their associated sets of generators,
 enable  one to modify the argument that was used in the proof of theorem 4.2 (the IET case), to
prove theorem 4.10 (the Levitt case) under the bounded periodicity assumption.

As in the proof of theorem 4.2 in the IET case, while shortening the lengths of the elements $v^t_1,\ldots,v^t_{f_t}$ using
the Rips machine (moves (3)-(4)), the ratios between the lengths of the elements $v^t_1,\ldots,v^t_{f_t}$ may not be bounded.
Hence,
we  need to replace the elements $v^t_1,\ldots,v^t_{f_t}$ by elements that are not longer than them, and for which the ratios
between their lengths are universally bounded. As in the IET case (proposition 4.4), in obtaining such new sets of generating elements,
there is a difference if there is a bound on the periodicity or not. This is why the assumption on the periodicity is added in the statement of
theorem 4.10.

\vglue 1pc
\proclaim{Proposition 4.12 (cf. proposition 4.4)} Suppose that there exists an integer $c_p$, such that the periodicity of the elements
$h_n(s_1),\ldots,h_n(s_r)$ is bounded by $c_p$ for all integers $n$.

There exists a subsequence of indices $t$,
such that for every index $t$ (from the subsequence),
the elements
$v^t_1,\ldots,v^t_{f_t}$ can be replaced by elements $u^t_1,\ldots,u^t_{g}$, with the following properties:
\roster
\item"{(1)}" for every index $t$ (from the subsequence), the semigroup that is generated by $u^t_1,\ldots,u^t_{g}$
contains the semigroup that
is generated by $v^t_1,\ldots,v^t_{f_t}$.

The homomorphism $h_n:F \to F_k$  extends to a homomorphism $\hat h_n$ from the residually
free group $RF_t$ that is generated by $u^t_1,\ldots,u^t_g$ to $F_k$.

\item"{(2)}" for each index $t$, and for large enough index $n$, $\hat h_n(u^t_i) \in FS_k$ for every index $i$, $1 \leq i \leq g$.

\item"{(3)}" for every $t$ it is possible to associate a presentation of the residually free group $RF_t$ with the
generalized equation that is associated with  the generators:
$u^t_1,\ldots,u^t_{g}$.

The original free group  $F$ is a subgroup of the group that is generated by:
$u^t_1,\ldots,u^t_{g}$.

The presentations of the residually free groups $RF_t$
and the presentations of the free group $F$
that are associated with the
elements $u^t_1,\ldots,u^t_g$ are similar for all $t$. i.e., for $t_4,t_5$
the presentations are identical if we replace
$u^{t_4}_i$ with $u^{t_5}_i$, for $i=1,\ldots,g$.

\item"{(4)}"
As in the IET case, each of the generators $s_j$, $j=1,\ldots,r$,
can be written as a word, $s_j=w^{t_1}_j(v^{t_1}_1,\ldots,v^{t_1}_{f_{t_1}})$. Note that although the elements $s_j$ are mapped to the positive cone
in $F_k$ by the homomorphisms $\{h_n\}$, and so are the elements $v^{t_1}_i$, the words $w^{t_1}_j$ are not necessarily positive words in the elements,
$v^{t_1}_1,\ldots,v^{t_1}_{f_{t_1}}$.

By part (1) each of the elements $v^t_i$ can be presented as a positive word in the new generators:
$u^t_1,\ldots,u^t_{g}$. When we substitute these positive words in the words $w^t_j$, we get the words:
$s_j=z^t_j(u^t_1,\ldots,u^t_{g})$. Note that with this notation, the elements $s_j$ are elements of the free group $F$, but the elements
$u^t_1,\ldots,u^t_g$ belong to the residually free group $RF_t$ that is generated by:
$v^{t_1}_1,\ldots,v^{t_1}_f$ or rather in a limit quotient of it.

From the sequence of homomorphisms $\{\hat h_n:<v_1,\ldots,v_f> \to F_k\}$ it is possible to extract a subsequence that converges into an
action of a limit group $L$ on a real tree $\hat Y$. $F<L$ and the tree $\hat Y$ contains the limit tree $Y$ (on which $F$ acts with
indecomposable Levitt dynamics).

As in the IET case (part (5) in proposition 4.4), with each word, $w^t_j$ and $z^t_j$, we can associate a path in the tree $\hat Y$,
$p_{w^t_j}$ and $p_{z^t_j}$. For a fixed $j$, $1 \leq j \leq r$,
$p_{w^t_j}$ and $p_{z^t_j}$  are identical, and they are identical paths for all $t$.

\item"{(5)}" there exist positive constants $d_1,d_2$, such that for every $t$,
there exist a sequence of indices:  $1 \leq i_1 < i_2< \ldots < i_{b} \leq g$, such that
for every $1 \leq  m_1 < m_2 \leq  b$:
$$d_1 \cdot length (u^t_{i_{m_1}})  \, \leq \, length (u^t_{i_{m_2}}) \, \leq \,
 d_2 \cdot length (u^t_{i_{m_1}}).$$
Furthermore,  for every index $i$, $1 \leq i \leq g$ for which $i \neq i_m$, $m=1,\ldots,b$:
$$ 10g \cdot d_2 \cdot length (u^t_i)  \, \leq \, length (u^t_{i_1}) $$

\item"{(6)}"
For each index $t$,
in the words $z^t_j$, $1 \leq j \leq r$, in a  distance bounded by some constant $c$,
either before or after the occurrence of an
element $u^t_i$, for $i \neq i_1,\ldots,i_{b}$, appears one of the elements $u^t_{i_m}$, $1 \leq m \leq b$.
\endroster
\endproclaim

\nfp
As in the IET case, to get the set of new generators $u^t_1,\ldots,u^t_{g}$ we need to modify the
Rips machine or the Makanin combinatorial algorithm. The Rips machine runs two processes. The first is
a process that erases subsegments that are covered exactly once, and if that process terminates, the machine
runs a second process that applies a sequence of entire transformations until there is a subsegment that is covered exactly
once, or until every segment is covered exactly twice.

The pseudogroups that are defined by the finite set of pairs of bases that we analyze, are all Levitt (indecomposable) pseudogroups. Hence,
if one applies the Rips machine to them, the first
process doesn't terminate.
Since our aim is to construct generators with comparable lengths, we need to separate
between long bases  and short ones. When we start with such separation, and modify accordingly the
Rips machine, it may happen that even though we start with a Levitt pseudogroup, the first process terminates after
finitely many steps (or doesn't apply at all), and then one  applies the second process finitely many times.

Recall that by lemma 4.11, for $t>t_0$  no non-trivial word in the generators,  fix pointwise
a non-degenerate subinterval of the interval $I$. Furthermore, in applying the first process in the Rips machine
to the pseudogroup that is defined by the pairs of bases, only moves (3) and (4) of this process are applied.

As in the IET case, we start by dividing the pairs of bases into long and short ones. First, we divide the elements,
$v^t_1,\ldots,v^t_{f_t}$, into finitely many sets according
to their length, from the longest to the shortest.  We place a separator between consecutive sets of elements
whenever there is a pair of consecutive elements (ordered according to length) that satisfy:
$$length(v^t_i) \geq c_1(f_t,c_p) \cdot length(v^t_{i+1})$$
where $c_1(f_t,c_p)=4f_tc_p$.

If there are no separators, we set
$g_t=f_t$, and $u^t_i=v^t_i$, $i=1,\ldots,f_t$. Suppose that there is a separator. In that case
we call the elements
that are in the last (shortest) set $short$ and the elements in all the other (longer) sets $long$.
We further call a pair of bases $long$ if at least one (hence, both) of the bases can be presented as a product (in the semigroup that is generated
by the elements) that contains a long element. Otherwise it is $short$.

The construction of the new set of generating elements and bases starts with (possibly none) applications of a modification of the first
process in the Rips machine.
By lemma 4.11, for every $t>t_0$ the process that is used to construct the pseudogroups in step $t$,  applies only moves
of type (4) and (3). In each step of the modified process we first apply moves of type (4) as long as possible. Since each such
move erases a pair of bases, one can apply only finitely many such moves.
We continue a step by applying a modified move of type (3) if such modified move is possible.

Let $b_1$ be a base with a paired base $b_2$, and suppose that $J_1$ is a maximal subinterval of $b_1$ ($J_1$ is
contained in the interior of $b_1$) that
is covered only once. $J_1$ is one of the elements $v^t_1,\ldots,v^t_{f_t}$. If it is not a long element, we do not
perform a modified move of type (3) along $J_1$. If $J_1$ is long  we perform move (3).

By the proof of lemma 4.11 every point in $J_2$ is covered by $b_2$ and at least one additional
base. As in the proof of lemma 4.11, if there is a point in $J_2$ which is the endpoint of a base, and this endpoint
 is covered only once after erasing $J_2$ from $b_2$, the Euler characteristic of the graph that is associated with the
new pseudogroup has increased by at least 1. In this case,
we replace the system of generators $v^t_1,\ldots,v^t_{f_t}$ and bases, by the new system of generators and bases
that is obtained after the modified move (3), and start the whole process again. The graphs that are
associated with the new pseudogroups that we obtain have strictly bigger Euler characteristics than the previous ones.
Since their Euler characteristic is strictly bigger, and
it has to be negative, such a
replacement can occur only finitely many times. Hence, in the sequel we may assume that such an increase in the Euler characteristic
does not happen along the modified process that we present.

After performing (the modified) move (3), we get a new set of bases. The endpoints of all the new bases are covered at least
twice (so no move of type (2) is required as we argued in the proof of lemma 4.11). We cut a long pair of bases, took out a subinterval, and obtained two new
pairs of bases. As in the IET case, we say that a pair of bases is $seconadry$ $short$ it its length is less that $\frac {1} {c_1(f_t,c_p)}$ times the
length of the
previous shortest (short) base. It is said to be $short$ if it is not secondary short and not longer than $c_1(f_t,c_p)$ times the maximal length of a previous
short base. A (new) pair of bases is $long$ if it is not short nor secondary short.

If the set of short bases was changed, it may be that some previous long bases become short, i.e., that there are no separators between them and the new
set of short bases. In that case we declare previously long bases that their length is bounded above by $c_1(f_t,c_p)$ times the maximal length of the
new short bases to be short. If the set of short bases was changed again, we repeat this process of possibly changing previously long bases to be short. There
were finitely many (previously) long bases, so this iterative process terminates after a finite time.

We repeat these modified steps in the first part of the Rips machine iteratively. We first perform all the possible moves of type (4),
and then perform a modified step of
type (3). i.e., we perform move (3) only along a long subinterval i.e., along a subinterval with length that is bounded below by $c_1(f_t,c_p)$ times
the maximal length of a short base.

As in lemma 4.11, since the endpoints of all the new bases are covered at least twice, no move of type
(2) is required along the modified process.
Afterwards we set some of the new bases  to be secondary short, short and long as we did in the first step. Note that short and secondary short bases
are not touched along the process, apart from being identified with other bases when performing a move a type (4).

The modified process terminates after finitely many steps, since at each step we erase a long subinterval, i.e., a subinterval of length
of length
bounded below by $c_1(f_t,c_p)$ times the maximal length of a short base.
The process terminates
when there are no long subintervals that are covered exactly once. Long bases may be cut along the process,
and at least one of the
 new bases that are obtained after the new cuts are either long or short.

\medskip
Although after the first (modified) process there are no long subintervals that are covered exactly once,
there may still remain long bases
after that process. Hence, although we analyze pseudogroups of Levitt type,
what is left to analyze are long bases with long subintervals that are covered
at least twice, and only short or secondary short subintervals that are covered once.
To analyze these long bases and shorten them we use a variation of the process
that is used in the proof of proposition 4.4.

Each subinterval that is covered only once is either short or secondary short. Over each such subinterval we add a new (formal) base, with no paired base, that
covers precisely the subinterval that is covered only once. With these new formal bases, that have no paired bases, every subinterval is covered at least twice,
so we can run the modified procedure that was used in the IET case, in the proof of proposition 4.4. Just that in analyzing IET pseudogroups every subinterval was
covered twice, and now we have to analyze pseudogroups in which every subinterval is covered at least twice.

Let $U_t$ be the pseudogroup (defined by the pairs of bases) that we ended the first step of the procedure with.
$U_t$ is supported on a finite union of oriented intervals. Every point in these intervals is covered at least twice.

Let $I$ be one of the intervals on which the pseudogroup $U_t$ is supported. $I$ is an oriented interval.
Let $pv$ be the vertex at the positive end of $I$. We perform the following operations on the bases of the pseudogroup $U_t$.
\roster
\item"{(1)}" Suppose that $pv$ does not belong to a subinterval that supports a short
or a secondary short base. Since no non-trivial word in the generators of the pseudogroup fix a non-degenerate subinterval, the set of long bases
that are adjacent to $pv$ belong to distinct pairs.

Let $b_1$ be the longest (long) base that is adjacent to $pv$. Let $b_2$ be the base that contains the other endpoint of $b_1$, such that
$b_2$ has the closest endpoint to $pv$. We transfer all the bases that their support is contained in the support of $b_1$ to the the dual base of $b_1$.
We cut $b_1$ at the endpoint of $b_2$ that is closer to $pv$. We further erase the part of $b_1$ that contains $pv$ and the corresponding part from its dual base.

We redefine the set of short, long and secondary short bases, given the new remaining part of $b_1$ as we did in the first procedure. If the remaining part of $b_1$ is a long base, the number of long, short and secondary short bases did not change. If the remaining part of $b_1$ is short or secondary short, the number of
pairs of long bases is reduced by 1.

Furthermore, if what was cut from $b_1$ is  only a short or a secondary short piece, and what is left
from $b_1$ is still long, then only long bases are adjacent to the new endpoint of the interval that supports the remaining of $b_1$ after we performed part (1).
These long bases have to be different than the long bases that were adjacent to $pv$, except the remaining of $b_1$.

\item"{(2)}" Suppose that $pv$ is the endpoint of both long and short bases. We act as in part (1). We set $b_1$ to be the longest long base with an
endpoint in $pv$. We transfer all the bases that are supported on the interval that supports $b_1$. We define $b_2$ as in (1), cut $b_1$ accordingly, and erase
the part of $b_1$ that contains $pv$, and the corresponding part of its dual base.

We redefine the set of short, long and secondary short bases, given the new remaining part of $b_1$ as we did in the first procedure. If the remaining part of $b_1$ is a long base, the number of long, short and secondary short bases did not change. If the remaining part of $b_1$ is short or secondary short, the number of
pairs of long bases is reduced by 1.
Furthermore, if what was cut from $b_1$ is  only a short or a secondary short piece, and what is left
from $b_1$ is still long, then only long bases are adjacent to the new endpoint of the interval that supports the remaining of $b_1$ after we performed part (2),
and these long bases are distinct than the long bases that are adjacent to $pv$ except the remaining part of $b_1$.

\item"{(3)}" Suppose  that $pv$  is the endpoint of only short and secondary short bases (and no long bases). In this case we set $b_1$ to be the longest long
 base, with
the closest endpoint to $pv$ among all the long bases.

We cut all the short and secondary short bases that have the endpoint of $b_1$ that is close to $pv$ in their interior.
We
 cut their dual bases accordingly.
Hence, we may increase (at most doubled) the number of pairs of short and secondary short bases.

By performing the cut of short and secondary short bases at the endpoint of $b_1$,
no base has the endpoint of $b_1$ that is close to $pv$ in its interior.

We continue from the endpoint of $b_1$ (that is closer to $pv$). Since $b_1$ is a long base and is adjacent to this base, we can act as in part (1) or (2).
We do not move the bases that are between $pv$ and the endpoint of $b_1$.

We further update the list of short, long and secondary short bases. Whenever the endpoint of the active subinterval
is adjacent to a long base we apply parts (1) or (2) of the procedure. Indeed, as long as no long base is changed to a short or secondary short one,
there will be a long base that is adjacent to the endpoint of the active subinterval. Hence, overall we possibly doubled the number of short and secondary short
bases,
and after finite number of steps we decrease the number of long bases by at least 1.
\endroster

Parts (1) and (2) do not increase the number of bases, and a finite number of iterations of actions with them lead to a reduction in the number of pairs of
long bases. Part (3) may at most double the number of short and secondary short bases.
An
action of part (3) must be followed by a finite iteration of actions of type (1) and (2) until the number of pairs of long bases is reduced.

Hence, after a finite number of steps all the remaining bases become short or secondary short. The number of times that part (3) is applied is bounded by
the original number of pairs of long bases. Hence, the number of (short and secondary short) bases that the procedure ends with is bounded by
${f_t}2^{f_t}$, and the number of elements that they define is bounded by $2{f_t}2^{f_t}$.

The final pseudogroup that is obtained when the second procedure terminates, contains no long elements, and the total
number of elements is bounded in terms of the number of elements in the original pseudogroup. By passing to a
subsequence we may assume that it is fixed.

The procedure that we used (that is composed from two procedures) is a modification of the Makanin-Razborov procedure, and uses its moves.
Parts (1) and (2) in proposition 4.12 are straightforward from the Makanin-Razborov moves that we use.
Hence, it is
possible to associate with it a finite graph with finitely many fundamental sequences, as in Razborov's work (see [Ra] and [Ca-Ka]). This finite graph
implies that the set of generators $u^t_1,\ldots,u^t_g$ is associated with finitely many generalized equations, hence, the sets of generators are associated
with finitely many presentations of residually free groups. After passing to a subsequence we may assume that they are all associated with the same presentation
of the residually free group $RF$ and we get part (3).

To prove the part (3) of the proposition, i.e., to prove that there are only finitely many presentations of the free group $F$ that are
associated
with the various generalized equations along the process, we start by analyzing the first part of the procedure,
the one that erases a long segment that is covered only once at each step.


By proposition 4.8 in [Be-Fe1] there are only finitely many closed leaves in a band complex up to homotopy.  The  moves that are used in the
first process can only decrease the length of the closed leaves in the associated band complexes. Since the number of generators
in these band complexes are bounded and so are the lengths of the relations (which are the lengths of the closed leaves in the
built from the pairs of bases that are the output of the first part of the procedure, i.e., that are associated with the tuples
$v^t_1,\ldots,v^t_{f_t}$. Hence, after passing to a subsequence, we can assume that they are all
the same presentations of the free group $F$.

In the second process we perform a sequence of entire transformations where the carrying bases are always long bases.
When the procedure terminates the outcome is a collection
of boundedly many short bases. Makanin in his seminal work defined the excess of a $generalized$ equation (see definition 4.8 in [Ca-Ka]).
Makanin's excess is set to be
the sum of the lengths of the bases in the generalized equation minus twice the length of the segments that support the collection of bases.

Since every point in the
supporting segments in our second procedure is covered at least twice, the excess of the generalized equations in the second procedure is always non-negative.
Entire transformations
on generalized equations in which every point is covered at least twice preserve the excess of the equations. Since we assume that the periodicity is bounded,
no abelian parts appear along the process, i.e., no long base overlaps with its dual with a short or secondary short shift between them.
By part (3) in lemma 4.11 no non-trivial word acts
trivially on a non-degenerate segment. Hence,  no
base overlaps exactly with its dual base, so no pair of bases can be eliminated along the second process.

Therefore, the excess of the generalized equations is preserved along the second procedure, and since the procedure terminates with a bounded
collection of short bases, the
excess has to be bounded by a bounded from above by twice the sum of the lengths of the (final) short bases, which is a bounded multiple of the
sum of the lengths of the  original short elements.

Hence, in the case
of bounded periodicity, the second procedure that we use is essentially similar to the procedure that was used in the IET case with an additional
bounded periodicity assumption.
i.e., the long elements on which the entire transformations along the process are performed cover the segments that support them exactly twice,
except for some short corrections
at their ends. Since these type of IET pseudogroups correspond to only finitely many presentations of the ambient free group $F$,
we can pass to a further subsequence for
which the presentations of $F$ that is associated with the  generalized equations along our second process are the same and we get part (3).

\medskip
Properties (4) and (5)  of the subsequence of generators $u^t_1,\ldots,u^t_g$ that appear in the statement of proposition 4.12, follow by the same arguments
that were used in the IET case, i.e., in the proof of proposition 4.4. To prove part (6) of the proposition we need to modify the argument that was
used to prove part (7) of proposition 4.4.

The procedure that we used to construct the elements $\{u^t_i\}$ is built from two procedures. The first procedure cuts long elements that are
covered only once (and cut a corresponding segment from the dual base), and the second procedure uses iterations of modified elementary transformations.
Hence, to prove part (6) of the proposition we can prove that it holds for each of the two procedures,

We start with the first procedure and denote the elements that it terminates with: $\hat v^t_1,ldots,\hat v^t_{m_t}$. Note that $m_t$ is universally bounded since
the number of bases along the first (and the second) procedure is universally bounded.

If the first procedure terminates with no secondary short elements, the claim of part (6) is empty for this procedure. In particular,
if we start the procedure with no long element that is covered only once, we don't need to run the first procedure, and part (6) holds for its output.

In each step of the first procedure a long element that is covered only once is cut from a base and from its dual base. Such an operation does not effect
the elements that are supported by the original base, and it may cut elements that are supported by the dual base. It may cut a single long element in two
places, where in between the two cuts there is a new long element, or it cuts two distinct elements, each in one cut.

Suppose first that along the first process only operations of type (3) are performed and no operation of type (4). i.e., we only cut long elements and do
not transfer pairs of bases with the same support where no other base supports part of these bases. In that case each operation of type (3) increases the
number of connected components. Since the Euler characteristics of the graphs that are associated with the generalized equations along the
process are fixed, and the valence of each vertex in this graph is at least 3, the numbers of connected components, which are the number of the vertices of the
associated graphs, are universally bounded. Hence, it is possible to apply operation (3) only bounded number of times. Therefore, a long element can be cut a
bounded number of times along the process and part (6) of the proposition follows in this case.

If the first process includes operations of type (3) and operations of type (4) only along long bases, then the same argument implies part (6) of the
proposition. Hence, we need to consider operations of type (4) along short or secondary short bases. Note that a short or secondary short bases are not touched
along the process, and there is a global bound on the number of pairs of short bases at each step of the procedure, so there is a global bound on the equivalence classes
(i.e., the orbits) of short and secondary short bases on which operations of type (4) are performed along the process.

Furthermore, short and secondary short elements do not change their place along the process, and can be cut in boundedly many places,    unless they are part of a pair
of (short or secondary short) bases that are moved using operation of type (4). There are only boundedly many pairs of short bases that the procedure starts with,
so apart from boundedly many operations of type (4), operations of type (4) are performed on pairs of short or secondary short bases that were cut from long bases
along the process.

A short or secondary short pair of bases on which an operation of type (4) is performed, and the pair is  cut from long bases along the process,
 must have long elements that were cut from the long  bases that contained the pair of bases, at least from one of the sides of each of the short or secondary short bases.  In fact
apart from boundedly many such pairs, long elements are cut from the long bases that contained  the pair of short or secondary short bases from the two sides of these bases.

Hence, each original long element $v^t_i$ can be written as a positive word  $v^t_i=w(\hat u^t_1,\ldots,\hat u^t_g)$, where the word $w$ is a product of  positive words
in the elements, $\hat u^t_1,\ldots,\hat u^t_g$,
$w_1,\ldots,w_{e(t,i)}$, each of the words,
$w_1,\ldots,w_{e(t,i)}$, is either short or long, and each of them is a positive word in which some of the short and secondary short elements from the set
$\hat u^t_1,\ldots,\hat u^t_g$  appears at most once, and in addition short and secondary short bases that were cut from long bases along the process, and on which
operations of type (4) were performed may appear, and between the
appearances of these bases there are long or short elements from the set
$\hat u^t_1,\ldots,\hat u^t_g$.

This presentation of a long element $v^t_i$ as a product of words,
$w_1,\ldots,w_{e(t)}$, that are either short or long, and each of them contains  either boundedly many of the elements
$\hat u^t_1,\ldots,\hat u^t_g$, so at least one of these elements  has to be short or long, or it contains in addition short or secondary short bases
on which operations of type (4) was performed and they were cut from long bases, in which case there are long or short elements between the appearances of these
short or secondary short bases. In both cases there is some bound $c$ such that every subword of length $c$ in the word
$v^t_i=w(\hat u^t_1,\ldots,\hat u^t_g)$ contains a short or a long element from the set of elements:
$\hat u^t_1,\ldots,\hat u^t_g$, and part (6) holds for the collection of elements that was produced after the first procedure.

\smallskip
So far we proved that the elements that the first procedure terminates with satisfy part (6) of the proposition. It is left to prove that part (6) holds for the
output of the second procedure.

If all the indices of the elements $\{u^t_i\}$ belong to the set: $i_1,\ldots,i_{b}$,
 the claim in part (6) is empty.
A short (and not secondary short) element, $\hat u$, that is created along the iterative procedure to construct the final elements $\{u^t_i\}$,
can be cut to at most the number of final elements,
so at least one of these new elements has to be short (and not secondary short), so part (6) follows for all the elements that appear in the
presentation of $\hat u$ as a product of final elements from the set $\{u^t_i\}$.

Let $\hat u$ be a long element that is created along the iterative process to construct the elements $\{u^t_i\}$. If $\hat u$ is replaced in the next step
in the process by a product of short and secondary short elements, then at least one of these elements has to be short, and their number is bounded by
the final number of elements, so part (6) follows for all the elements that appear in the presentation of $\hat u$ as a product in the final set of elements
$\{u^t_i\}$.

Otherwise in the first step that $\hat u$ is changed,
$\hat u$ is replaced by a product that contains a new long element. If
$\hat u$ is replaces by a product that starts and ends with a long element, then the number of short and secondary short elements between
two long elements in the product is bounded by the total number of elements $\{u^t_i\}$, which is universally bounded,
so part (6) follows by an iterative application of the claim for long elements that are created along the process.

Suppose that $\hat u$ is replaced by a product that contains a long element, and starts with a short or a secondary short element. The number of
times that such a replacement can occur for a consecutive sequence of replacements of the starting long elements, in the product that
represent $\hat u$, is bounded by the
total number of short and secondary short elements,
before the number of long pairs of bases is reduced by at least 1. Since the total number of elements is universally bounded, and so is the total number
of pairs of long bases, there is a universal bound on the number of times that such a replacement can occur in a product that
represent a long element that is created along
the iterative process in the next steps of the process.

Suppose that $\hat u$ is replaced by a product that terminates with a long element followed by a product of secondary short elements. Let
$te$ be the maximal number of terminating elements $\{u^t_i\}$, and let $ba$ be the maximal number of starting long bases that are associated
with the elements $\hat v^t_1,\ldots,\hat v^t_{\hat f_t}$. After the first replacement of $\hat u$, the number of secondary short elements that may appear after
the terminal long element in
a product that represent $\hat u$ is bounded by $te$, the number of terminal elements in the whole process. Hence, our goal is to argue that such a replacement
can occur a bounded number of times.

Recall that in the iterative process that produces the elements $\{u^t_i\}$, step (3) can occur at most once before the number of long
pairs is reduced by at least 1. In a part of the process in which consecutively a long base is adjacent to the terminal point (that we denoted $pv$),
usually after step (3) was
applied, and until the number of pairs of long elements is reduced, the number of times that secondary short bases are adjacent to the terminal point
is bounded by ${ba}^2$, and the number of times that the terminal part of a long base is replaced by a product that terminates
with a secondary short base is bounded by the number of pairs of long bases, that is bounded by $ba$.

Therefore, in each step of the procedure along such a part of the procedure,
a long element $\hat u$ is replaced by a product, and if the product contains a long
element, then this long element is followed by a consecutive sequence of secondary short elements before a short or a long element,
and the length of such consecutive sequence of secondary short elements is bounded by $4{ba}^2+te$.

This  bound is for a part of the process in which the number of long pairs of bases is reduced by 1, and the number of such
parts is bounded by the initial number of long pairs of bases.
Hence, altogether, in each step of the procedure, the long element $\hat u$ is replaced by a product in which a consecutive sequence of secondary short elements
has length that is bounded by  ${ba} \cdot (4{ba}^2+te)$.
This proves part (6) of the proposition.

\hfill$\square$

As in the IET case, proposition 4.12 enables  the proof of the key claim  in the proof of theorem 4.10 - the combinatorial
bounded cancellation along the process that constructs the new sets of generators:
$u^t_1,\ldots,u^t_{g}$.

\vglue 1pc
\proclaim{Proposition 4.13}
Suppose that there exists an integer $c_p$, such that the periodicity of the elements
$h_n(s_1),\ldots,h_n(s_r)$ is bounded by $c_p$ for all integers $n$.

With the notation of proposition 4.12, let $L$ be a limit group which is the limit of a (convergent) subsequence of the homomorphisms:
$\{\hat h_n:<u^{t_1}_1,\ldots,u^{t_1}_g> \to F_k\}$. Note that $L$ is a quotient of the residually free group $RF$ that is associated
with the generalized equation that is associated with the final elements,
$u^{t_1}_1,\ldots,u^{t_1}_g$.
And the free group $F$, that acts on a real tree as an indecomposable Levitt action,
is a subgroup of $L$.

There exists a constant $C>0$,
so that for a subsequence of the indices $t$, that for brevity we
still denote $t$,
the words $z^t_j$, $1 \leq j \leq r$, can be replaced by words $\hat z^t_j$ with the following properties:
\roster
\item"{(1)}"
 As elements in the limit group $L$:
$\hat z^t_j(u_1^{t},\ldots,u^{t}_{g})=z^{t}_j(u_1^{t},\ldots,u^{t}_{g})$. Note that the element $z^t_j$ represents the element $s_j$ which is an element in
$F<L$.

\item"{(2)}" $\hat z^t_j$ is obtained from $z^{t}_j$ by eliminating distinct pairs of subwords. Each pair of  eliminated
subwords corresponds to two subpaths of the path
$p_{z^{t}_j}$ that lie over the same segment in the (finite) tree $T_{w_j}$, where the two subpaths have opposite orientations.

\item"{(3)}" With the word $\hat z^t_j(u^t_1,\ldots,u^t_{g})$ we can naturally associate a path in a real
tree $\hat Y$, that contains the real tree $Y$ on which $F$ acts, that we denote, $p_{\hat z^t_j}$. The path $p_{\hat z^t_j}$ can be naturally divided into
subsegments according to the appearances of the subwords $u^t_i$ in the word $\hat z^t_j$.

Let $DB_{\hat z^t_j}$ be the number of such subsegments that are associated with subwords $u^t_i$ in $p_{\hat z^t_j}$, that at least
part of them is covered more than once by the path $p_{\hat z^t_j}$. Then for every $t > t_0$ and every $j$, $1 \leq j \leq r$,
$DB_{\hat z^t_j} \leq C$.
\endroster
\endproclaim

\nfp
By the argument that was used to prove proposition 4.5 from proposition 4.4 in the IET case, under the bounded periodicity assumption, there are global bounds on the dual positions of each pair of short
elements from the generating sets $u^t_1,\ldots,u^t_g$ along cancelling paths.
Finally, in case of bounded periodicity, the rest of proposition 4.13 follows from proposition 4.12 by exactly the same argument
that proposition 4.5 follows from proposition 4.4 in the IET case.

\hfill$\square$

Proposition 4.13 is the analogue in the Levitt case of  proposition 4.5 in the IET case. In the IET case  proposition 4.5 was used to prove
theorem 4.2.
By the same argument that proposition 4.5 was used to prove theorem 4.2 in the IET case, theorem 4.10 follows from proposition 4.13 in the Levitt case
(note that the argument is valid under the bounded periodicity assumption).

\hfill$\square$

\medskip
In the IET case we first proved theorem 4.2 under the bounded periodicity assumption. i.e. we constructed a sequence of shortening
automorphisms for the tails of the given sequence of homomorphisms $\{h_n\}$ of a pair $(S,Q)$ where $Q$ is a surface group  assuming the images
of the homomorphisms $\{h_n\}$ have bounded periodicity. Then we modified the procedures and the arguments that we used to prove the same
result omitting the bounded periodicity assumption.

In the Levitt case we didn't manage to get a similar modification. Hence, to generalize theorem 4.10 and omit the bounded periodicity assumption we use a
different strategy. We use the proof of Bulitko's lemma (see lemma 1.3 in [Ma1]) to replace each homomorphism from the given sequence,
$\{h_n:(S,F) \to (FS_k,F_k)\}$, by a possibly shorter homomorphism
$\hat h_n:(S,F) \to  (FS_k,F_k)$, so that the new sequence of homomorphisms $\{\hat h_n\}$ have uniformly bounded periodicity, where the bound on the periodicity
depends only on the pair $(S,F)$.

\vglue 1pc
\proclaim{Proposition 4.14}
Let $F$ be a finitely generated free group, and let $S=<s_1,\ldots,s_r>$ be a f.g.\ subsemigroup of $F$ that generates $F$ as a group.
There exists a positive integer $c_p$ that depends only on the pair $(S,F)$, so that for any pair  pair homomorphism,
$h_n:(S,F) \to (FS_k,F_k)$, there exists a pair homomorphism,
$\hat h_n:(S,F) \to (FS_k,F_k)$, with the following properties:

\roster
\item"{(1)}" $\hat h_n$ is obtained from the homomorphism $h_n$ by iteratively replacing periodic subwords $\alpha^m$, $m>c_p$, in $h_n(s_1),\ldots,h_n(s_r)$, by
periodic subwords $\alpha^t$ where $t \leq c_p$.

\item"{(2)}" as a result of (1) the lengths of the elements in the tuple:
$\hat h_n(s_1),\ldots,\hat h_n(s_r)$ are bounded by the corresponding lengths in the tuple:
$h_n(s_1),\ldots,h_n(s_r)$.

\item"{(3)}" the periodicity of the images,
$\hat h_n(s_1),\ldots,\hat h_n(s_r)$ is bounded by $c_p$.
\endroster
\endproclaim

\nfp Since the free group $F$ is f.p.\ and $S$ generates $F$, $F$ admits a finite presentation: $<s_1,\ldots,s_r \, | \, r_1,\ldots,r_q>$. Hence,
for any given pair homomorphism $h_n:(S,F) \to (FS_k,F_k)$,
$h_n(s_1),\ldots,h_n(s_r)$ that are all elements in $FS_k$ satisfy the finite system of equations $\Sigma$ that corresponds to this finite presentation of
$F$.

By the proof of Bulitko's lemma (lemma 1.3 in [Ma1]), there exists a constant $c_p$ that depends only on the system of equations $\Sigma$,
such that it is possible to iteratively replace periodic subwords in the tuple:
$h_n(s_1),\ldots,h_n(s_r)$ by periodic words of smaller exponent, to obtain a tuple of elements that we denote:
$\hat h_n(s_1),\ldots,\hat h_n(s_r)$, that are all elements in $FS_k$. The elements in the obtained tuple  have periodicity bounded by $c_p$, they shortenings of the
corresponding elements in the tuple:
$h_n(s_1),\ldots,h_n(s_r)$,
the obtained tuple of elements satisfy the system of equations $\Sigma$.

Therefore, the obtained tuple of elements can be extended to a pair homomorphism: $\hat h_n: (S,F) \to (FS_k,f_k)$ that satisfies the conclusions of the lemma.

\hfill$\square$

Lemma 4.14 enables one to shorten a general sequence of pair homomorphisms: $\{h_n:(S,F) \to (FS_k,F_k)\}$ that converges into an indecomposable Levitt action of $F$ on a real tree.
Indeed, it possible to start by replacing each of the homomorphisms $h_n$ by a homomorphism $\hat h_n$ that satisfies the conclusions of lemma 4.14. In particular,
the length of $\hat h_n(S_j)$ is bounded by the length of $h_n(s_j)$, $1 \leq j \leq r$,  and the periodicities of the elements $\hat h_n(s_j)$ are bounded by the universal bound $c_p$.

In general, the sequence $\{h_n\}$ may not converge to an action of $F$ on a real tree, but a subsequence converges. However, on the limit tree that is obtained from the convergent subsequence
$F$ may not act faithfully (but rather a limit quotient $L$ of $F$), and it may be that the limit action is not an indecomposable Levitt action, but rather a more general action
with which it is possible to associate a graph of groups of the limit group $L$.

At this point we have all the tools to shorten the convergent subsequence of homomorphisms $\{\hat h_n\}$. However, this requires the notions and techniques that are used in the last section in this book,
the section that constructs the Makanin-Razborov diagram, so we postpone the shortenings of general sequences of homomorphisms that converge into indecomposable Levitt components
(the generalization of theorem 4.10 to general sequences of homomorphisms with no periodicity assumptions)
to the last section (section 9)  in this book.

\vglue 1.5pc
\chapter{The canonical JSJ decomposition of (some) pairs}
\fullwidthtext
\medskip

In theorem 3.1 we have shown that with any given f.g.\ semigroup, $S$, and its set of
homomorphisms into the free semigroup $FS_k$, $Hom(S,FS_k)$, there is an associated (canonical)
finite collection of pairs, $(S_1,L_1),\ldots,(S_m,L_m)$, where each pair consists of
a limit group, $L_i$, and a semigroup, $S_i$, that is embedded in the limit group $L_i$ (as a subsemigroup),
and generates $L_i$ as a group. By the construction of the pairs, $(S_i,L_i)$, each of them is obtained as
a (maximal) limit from
a sequence of homomorphisms from the set, $Hom(S_i,FS_k)$.

Once we associated the canonical set of maximal pairs, $(S_1,L_1),\ldots,(S_m,L_m)$, with a f.g.\
semigroup $S$, to analyze the structure of the set of homomorphisms, $Hom(S,FS_k)$, it is sufficient
(and equivalent) to analyze the set of homomorphisms of each of the limit groups, $L_i$, into the
free group, $F_k$, that restrict to homomorphisms of the subsemigroups, $S_i$, into the free semigroup,
$FS_k$. We denote this set of homomorphisms (of pairs), $Hom((S_i,L_i),(FS_k,F_k))$.

In the case of free groups [Se1], to analyze the set of homomorphisms, $Hom(L,F_k)$, where $L$ is a limit group,
we used Grushko free decomposition to factor $L$ into a free product, and then associated the canonical JSJ decomposition
with each of the factors. With the JSJ decomposition we used its associated modular group to twist (or shorten)
homomorphisms, that allowed us to associate finitely many (maximal) shortening quotients with the  limit group.
Repeating this procedure iteratively we obtained the Makanin-Razborov diagram, in which every path
(called a resolution) terminates in a free group, and the set of homomorphisms from a free group into
the coefficient group $F_k$, can be naturally presented as $(F_k)^s$, where $s$ is the rank of the free group.

In the case of a semigroup, the geometric tools that are needed in order to analyze the set of homomorphisms,
$Hom((S_i,L_i),(FS_k,F_k))$, are based on the analogous tools for groups, but they need to be further refined,
as the modular automorphisms that can be used to modify (shorten) automorphisms are required to ensure
that the image of the subsemigroup $S_i$ remains a subsemigroup of the standard free semigroup $FS_k$. To ensure that we analyzed
axial, IET and Levitt actions of groups (or rather pairs) on oriented trees in the previous section.

As in analyzing homomorphisms into a free group, the basic object that we need to associate with a pair,
$(S,L)$, where $L$ is a freely indecomposable limit group and $S$ is a subsemigroup of $L$ that generates $L$,
is an analogue of a JSJ decomposition.

Unfortunately, we manage to construct a direct analogue of the JSJ
decomposition for groups only under further restrictions on a pair $(S,L)$.  For general pairs $(S,L)$ in
which the limit group $L$ is freely indecomposable, we replace
the JSJ decomposition with  finitely many sequences of decompositions, i.e., with finitely many $resolutions$ or $towers$.

For pairs for which we construct a JSJ decomposition, the JSJ decomposition is (at least partly) canonical,
but unlike the group analogue it
is not unique, i.e., we associate a (canonical) finite collection of decompositions with a
pair $(S,L)$ and not a single decomposition.

\medskip
Let $(S,L)$ be a pair consisting of a freely-indecomposable limit group, $L$, and a subsemigroup $S$ of $L$
that generates $L$ as a group. We look at all the sequences of homomorphisms, $\{h_n:L \to F_k\}$, that
restrict to semigroup homomorphisms of $S$ into $FS_k$, and converge into the pair $(S,L)$. By a theorem of
F. Paulin [Pa], each such sequence subconverges into an action of the limit group $L$ on a real tree $Y$, and
this action is stable (and even super stable in the sense of Guirardel [Gu]) by lemma 1.3 in [Se1]. By works of
M. Bestvina and M. Feighn [Be-Fe1], [Se3], and V. Guirardel [Gu], with a superstable action of the limit group
$L$ on the real tree $Y$ it is possible to associate (canonically) a graph of groups decomposition.

Therefore, with the pair $(S,L)$ we can associate a collection of graphs of groups decompositions of $L$,
i.e., those graphs of groups that are associated with actions of $L$ on real trees,  where these actions
 are obtained as
a limit from convergent sequences of homomorphisms (that algebraically converge into the limit group $L$).

Note that since all of these graphs of groups are
abelian decompositions, they can be all obtained from the abelian JSJ decomposition of the limit group
$L$. They are obtained from the JSJ decomposition of $L$ by cutting $QH$ vertex groups along some s.c.c.\ and then possibly collapse and fold
some parts of the obtained
graphs of groups.

Also, note that if a cyclic subgroup $C$ of $L$ stabilizes a non-degenerate segment in
a real tree that is obtained as a limit of actions of $L$ that are associated with homomorphisms of $L$ into
$F_k$, then the unique maximal cyclic subgroup of the limit group $L$ that contains $C$ stabilizes
this segment as well. Similarly, if an abelian subgroup $A<L$ stabilizes a non-degenerate segment in
such a real tree, then the unique direct summand in $L$ that contains $A$ as a subgroup of finite index,
 in the unique maximal abelian subgroup
of $L$ that contains $A$, stabilizes this non-degenerate segment as well.

On these graphs of groups we can naturally define
a partial order. We say that given two graphs of groups, $\Lambda_1$ and $\Lambda_2$, $\Lambda_1 > \Lambda_2$,
if $\Lambda_1$ is a proper refinement of $\Lambda_2$, or alternatively, $\Lambda_2$ is obtained from $\Lambda_1$
by (possibly) cutting some QH vertex groups along a finite collection of s.c.c.\ and then (possibly)
 performing some collapses and foldings.

\vglue 1pc
\proclaim{Proposition 5.1} Let $(S,L)$ be a pair of a freely indecomposable limit group $L$,
and its subsemigroup $S$ that generates $L$. Then there exist maximal abelian decompositions of the pair $(S,L)$.
Every strictly increasing sequence of abelian decompositions that are associated with $(S,L)$,
$\Lambda_1 \, < \, \Lambda_2  \, < \, \Lambda_3 \, < \, \ldots$, terminates after finitely many steps.
\endproclaim

\nfp All the abelian decompositions, $\Lambda_i$, are obtained from the abelian JSJ decomposition of
the freely-indecomposable limit group $L$, by cutting QH vertex groups along a finite (possibly empty)
collection of s.c.c.\ and then (possibly)
collapse and fold the obtained abelian decomposition. Since there is a bound
on the size of a collection of  disjoint non-homotopic non null homotopic s.c.c.\ on the surfaces that are
associated with the QH vertex groups in the JSJ decomposition of $L$, and the abelian edge groups of the
JSJ decomposition are all finitely generated, the proposition follows.

\hfill$\square$

Proposition 5.1 proves the existence of maximal elements in the set of abelian decompositions that are
associated with the pair, $(S,L)$, with their natural partial order. To construct an analogue of a
JSJ decomposition for the pair, $(S,L)$, we further prove that under further assumptions on the pair $(S,L)$,
there are only finitely many (equivalence
classes) of such maximal elements.

With a decomposition (graph of groups) of a group we associate its Bass-Serre tree. The group naturally acts on that tree.
We say that two decompositions of a group are in the
same equivalence class if their associated Bass-Serre trees are equivariantly isometric.

\vglue 1pc
\proclaim{Theorem 5.2} Suppose that $L$ is freely indecomposable, and that
all the maximal abelian decompositions, $\{\Lambda_i\}$, that are associated with the pair $(S,L)$,
correspond to simplicial actions of $(S,L)$ on real trees.
Then
there exist only finitely many (equivalence classes of)
maximal abelian decompositions of the pair $(S,L)$.
\endproclaim

\nfp Suppose that there are infinitely many (equivalence classes of) maximal abelian decompositions of
a pair $(S,L)$. i.e., infinitely many inequivalent maximal abelian decompositions of $L$ that are associated with actions of $L$ on real trees that
were obtained from convergent sequences of pair homomorphisms of the pair $(S,L)$. Let $\{\Lambda_i\}_{i=1}^{\infty}$ be the collection of these  (inequivalent)
maximal decompositions. Under the assumptions of the theorem, they are all simplicial.

With the maximal pair $(S,L)$ we associate a sequence of homomorphisms, $\{h_n:(S,L) \to (FS_k,F_k)\}$, that converges into
$(S,L)$.
To define the sequence of homomorphisms, we fix a (symmetric) generating set of the limit group $L$,
$\{g_1,\ldots,g_m\}$, that contains a
generating set, $s_1,\ldots,s_r$, of the subsemigroup $S$. With the given generating set of the pair, $(S,L)$,
we naturally associate its Cayley graph, that we denote $X$. For every positive integer $n$, we denote
the ball of radius $n$ in the Cayley graph $X$, $B_n$.

Each maximal decomposition, $\Lambda_i$, that is associated with the pair $(S,L)$, is obtained from a sequence of
homomorphisms, $\{f_i(j)\}_{j=1}^{\infty}$, from $(S,L)$ into $(FS_k,F_k)$. For each index $i$,
the sequence of actions of $(S,L)$ on the Cayley graph  of $(FS_k,F_k)$, that are associated with the sequence of
homomorphisms,
$\{f_i(j)\}_{j=1}^{\infty}$, converges in the Gromov-Hausdorff topology,
after
rescaling the metric so that the maximal length of the image of a generator (under $f_i(j)$) is 1, to an action
of the pair $(S,L)$ on a real tree, that we denote $T_i$.

For each index $n$, we choose a homomorphism $h_n: (S,L) \to (FS_k,F_k)$ that satisfies the following conditions:
\roster
\item"{(1)}" $h_n$ is one of the homomorphisms, $\{f_n(j)\}_{j=1}^{\infty}$.

\item"{(2)}" let $g$ be a word of length at most $n$ in the fixed set of generators, $g_1,\ldots,g_m$, of
the limit group $L$. Then $h_n(g)=1$ if and only if $g=1$ in $L$.

\item"{(3)}" Let $Y_n$ be the Cayley graph of $(FS_k,F_k)$ after rescaling the
metric so that the maximal length of the image of a generator (under $h_n$) is 1, and let $t_n$
be the base point in $T_n$. Then for every $n$ and every
$\ell_1,\ell_2 \in B_n$:
$$ d_{T_n}(\ell_1(t_n),\ell_2(t_n)) - \frac {1}{n} \, \leq \, d_{Y_n}(h_n(\ell_1),h_n(\ell_2)) \, \leq \,
   d_{T_n}(\ell_1(t_n),\ell_2(t_n)) + \frac {1}{n} .$$
\endroster

With each homomorphism $h_n$, there is naturally an associated action of the pair $(S,L)$ on the simplicial
tree $Y_n$ that is obtained from the Cayley graph of the coefficient free group $F_k$ by rescaling the metric
so that the maximal length of the image of a generator (under $h_n$) is 1.
From the sequence of homomorphisms, $\{h_n\}$, it is possible to extract a subsequence, so that the sequence
of actions of the pair $(S,L)$ on the trees $\{Y_n\}$, converges into a faithful action of the pair
$(S,L)$ on a limit real tree $Y$. Note that this action is precisely the limit of the corresponding subsequence of actions of
the limit group $L$ on the limit real trees $T_n$.

The action of $(S,L)$ on the limit tree $Y$ is non-trivial, has abelian stabilizers of non-trivial
segments,  and is super-stable in the sense of [Gu]. Hence,
with this action it is possible to associate a (graph of groups) decomposition $\Delta_1$ with trivial and
abelian edge stabilizers. Since $L$ was assumed to be freely indecomposable, all the edge stabilizers in
$\Delta_1$ are non-trivial abelian. Furthermore, under the assumptions of the theorem,
the abelian decompositions, $\{\Lambda_n\}$ and $\Delta_1$, must be simplicial, and they are all dominated by the JSJ
decomposition of the (freely indecomposable) limit group $L$.

\medskip
We  divide the edges in the abelian decomposition $\Delta_1$, and the vertices that are
adjacent to these edges into families. The stabilizer of every edge group in $\Delta_1$ is an abelian
subgroup of $L$, and each vertex group in $\Delta_1$ is finitely presented. We fix a finite generating set
for each of the vertex groups in $\Delta_1$.

Also, in the sequel we pass to the subsequence of the abelian decompositions $\{\Lambda_n\}$
that was used for the convergence into the simplicial tree with the abelian decomposition $\Delta_1$. We still denote the abelian
decompositions in this subsequence, $\{\Lambda_n\}$.

\vglue 1pc
\proclaim{Definition 5.3} Let $E$ be an edge in $\Delta_1$, and let $A$ be its edge group.
We say that an edge group $A$ in $\Delta_1$ is $stably$ $elliptic$, if every element
$a \in A$,  is elliptic in almost all the abelian decompositions $\Lambda_n$ (i.e., in all but at most finitely
many decompositions in the subsequence that converges into $\Delta_1$). Otherwise we say that $A$ is
$stably$ $hyperbolic$. In a similar way we define  stably elliptic and stably hyperbolic vertex groups.

Let $E$ be an edge in $\Delta_1$ with an edge group $A$, and suppose that $A$ is a stably hyperbolic edge group. Let $V$ be a vertex group that
stabilizes a vertex that is
adjacent to the edge $E$ in $\Delta_1$. We fix a generating set $v_1,\ldots,v_{\ell}$ of $V$.
For an element $f \in F_n$ we denote
the length of the conjugacy class of $f$ by $|f|$. We say that the edge group $A$ is $periodic$
 in the vertex group $V$, if for every element $a \in A$, there exists a positive constant $c_a>0$,
such that for almost every
index $n$ (i.e., for all except finitely many values of $n$),  every
point in the Cayley graph $T$ of the coefficient group $F_k$ is moved by at least one of the
elements, $h_n(v_1),\ldots,h_n(v_{\ell})$,
a distance of at least $c_a \cdot n \cdot |h_n(a)|$.

We say that a hyperbolic edge group $A$ is $non$-$periodic$ in $V$,
if there exists a  non-trivial element $a \in A$,
and a positive constant $c_a>0$, so that for all but finitely many indices $n$,
there exists a point in the Cayley graph $T$ of the coefficient group $F_k$ that is moved by each of the
elements, $h_n(v_1),\ldots,h_n(v_{\ell})$,
a distance that is bounded above by  $c_a \cdot |h_n(a)|$.
\endproclaim

By passing to a subsequence of the maximal abelian decompositions, $\{\Lambda_n\}$, we may assume that
every edge group in $\Delta_1$ is either stably elliptic or stably hyperbolic. By passing to a further subsequence, we
may assume that every stably hyperbolic edge group in $\Delta_1$ is either periodic or non-periodic in the one or two
vertices that are adjacent to it. We continue with such a subsequence  of the abelian decompositions
$\{\Lambda_n\}$.

At this point we gradually refine the abelian decomposition $\Delta_1$. Let $V$ be a vertex group in $\Delta_1$,
that is connected only to either stably elliptic or to periodic edge groups. If $V$ is not itself stably elliptic in a subsequence of the
abelian decompositions,$\{\Lambda_n\}$, we pass to this subsequence, and analyze the actions of the f.g.\ subgroup
$V$ on the Cayley graph $T$ of the coefficient group $F_k$ via the homomorphisms $\{h_n\}$. By passing to a further
subsequence, these actions of $V$ do converge  into a non-trivial action of $V$ on a real tree. Since all the edge
groups that are connected to $V$ in $\Delta_1$ are either stably elliptic or periodic, the abelian decomposition that is
associated with the limit action of $V$ can be further extended to an abelian decomposition of the pair $(S,L)$, that
strictly refines the abelian decomposition $\Delta_1$.

We repeat this refinement procedure for every non stably elliptic vertex group in the obtained (refined) abelian decomposition,
that is connected only to stably elliptic or periodic edge groups,
and is a point stabilizer in the corresponding action on a real tree. i.e., that is not a  vertex group in $\Delta_1$ or
a refined abelian decomposition that
is associated with an IET or an axial component.
By the accessibility for small splittings of f.p.\ groups [Be-Fe1]
(or alternatively,
by acylindrical accessibility ([Se],[De],[We])), this refinement procedure terminates after finitely (in fact
boundedly) many steps, and we obtain an abelian decomposition that we denote $\Delta_2$.

In $\Delta_2$ every vertex group $V$ that is
connected only to stably elliptic and periodic edge groups, is either associated with an IET or an axial component,
or it is stably elliptic itself, in which case  all the edge groups that are
connected to it are stably elliptic as well.

If all the vertex groups in $\Delta_2$ are stably elliptic, then  for almost all the indices $n$, the abelian
decomposition $\Lambda_n$ is dominated by the abelian decomposition $\Delta_2$ (i.e., $\Delta_2$ is a
(possibly trivial) refinement of $\Lambda_n$ for almost every index $n$). But in this case $\Delta_2$
dominates only finitely many non-equivalent abelian decompositions,
a contradiction to the existence of the infinite sequence of
non-equivalent maximal abelian decompositions $\{\Lambda_n\}$ (from which $\Delta_2$ was obtained).

Suppose that $\Delta_2$ contains a vertex group that is associated with an axial or an IET component. Note that such a vertex group
is connected to
only stably elliptic and periodic edge groups. Then there exists a sequence
of homomorphisms of pairs that converges into a non-simplicial
faithful action of $L$ on a real tree (and contradicts the assumption of theorem 5.2).

\vglue 1pc
\proclaim{Proposition 5.4} Suppose that $\Delta_2$ contains a vertex group that is associated with either an axial or an IET
component.
Then there exists a sequence of homomorphisms of pairs:  $\{\nu_n:(S,L) \to (FS_k,F_k)\}$ that converges into
a faithful action of $L$ on a real tree, and this action contains either an axial or an IET component.
\endproclaim

\nfp Suppose that $\Delta_2$ contains a vertex group $A$ that is associated with an axial component. In that case
$L=V*_{A_0} \, A$. By construction, the decomposition $\Delta_2$ is obtained from a sequence of homomorphisms $\{h_n\}$ using a
finite process of iterative refinements.

By theorem 4.1 the vertex
group $A$ can be written as $A=A_0+<a_1,\ldots,a_{\ell}>$, where $A_0$ is the point stabilizer of the axial
component, $\ell \geq 2$, and there exists some index $n_0$ so that for every $n>n_0$, $h_n(a_i) \in FS_k$, and for each $j$, $1 \leq j \leq r$,
$s_j$ can be written as a word:
 $$s_j \, = \, v^j_1w^j_1v^j_2w^j_2 \ldots v^j_{b_j}w^j_{b_j}$$
where each of the elements $w^j_i$ is a positive word in the basis elements $a_1,\ldots,a_{\ell}$, $v^j_i \in V$,
and for every index $n>n_0$, $h_n(v^j_i) \in FS_k$.

In that case, for every $n>n_0$,  we can modify the homomorphisms $h_n$, by preserving the images of elements in the vertex group $V$, and
modifying the images of the elements $a_1,\ldots,a_{\ell}$. By the properties of the homomorphisms $h_n$, for every
$n>n_0$ we can set (the modified images) $\nu_n(a_i)$ to be arbitrary elements in $FS_k$ that commute with $h_n(A_0)$. Clearly, we can choose these images
of $\nu_n(a_1),\ldots,\nu_n(a_{\ell})$, so that in the limit they will be independent over the rationals, and so that in the limit
the whole vertex group $V$ will stabilize a point. Therefore, the limit of the constructed homomorphisms $\{\nu_n\}$ will
contain a single axial component with an associated group $A$. In particular the limit action is faithful and contains
an axial component.

Suppose that $\Delta_2$ contains an IET component. Let $Q$ be a (hyperbolic) punctured surface group that is a QH vertex group in $\Delta_2$.
$\Delta_2$ is an abelian decomposition of
the limit group $L$.
To apply the techniques and the results that we obtained in the previous section (for shortenings homomorphisms), it is better to embed
$L$ and its QH subgroup $Q$ in a bigger limit group $\hat L$. If the punctured surface that is associated with the QH vertex group $Q$
has only one boundary component, we do not change $L$ and $Q$. If $Q$ has more than one boundary component, we first single out
a vertex that is connected to the QH vertex group $Q$ in $\Delta_2$. Then we  add Bass-Serre generators to
the limit group $L$ that conjugate the vertex groups that are connected to all the other boundary components of $Q$ (as is done in the
construction of a completion of a resolution in section 1 in [Se2]). We denote the
limit group that is generated by $L$ and the additional Bass-Serre generators $\hat L$.

$Q$ and $L$ embed naturally in $\hat L$ and the abelian decomposition $\Delta_2$ extends naturally to an abelian decomposition $\hat \Delta_2$ of
$\hat L$. The semigroup $S$ was assumed to generate $L$. We add to $S$ Bass-Serre generators that we added to $L$ that
are mapped by the homomorphisms $h_n$ to
$FS_k$ for large enough $n$. We get a f.g.\ semigroup $\hat S$ that embeds into $\hat L$ and generates it as a group. In particular, the homomorphisms,
$h_n:(S,L) \to (FS_k,F_k)$, extend to homomorphisms: $\hat h_n:(\hat S,\hat L) \to (FS_k,F_k)$.

We denote the subgroup of $\hat L$ that is generated by the punctured surface group $Q$ and the Bass-Serre generators that we added to $L$
by $\hat Q$. The finite set of Bass-Serre generators that we added to $S$ to get $\hat S$, together with a finite set of generators of
$S$ generate $\hat L$. We denote this generating set of $\hat S$, $s_1,\ldots,s_r$. We write each of the generators in a normal form w.r.t.\ the
abelian decomposition $\hat \Delta_2$ of $\hat L$.

These normal forms contain finitely many elements $q_1,\ldots,q_f \in \hat Q$ that together generate the subgroup $\hat Q$ of $\hat L$, and such that
for some index $n_1$, for every $n>n_1$, $\hat h_n(q_j) \in FS_k$, $j=1,\ldots,f$. We denote by $U$ the subsemigroup of $\hat Q$
that is generated by $q_1,\ldots,q_f$.

Theorem 4.2 is stated
for shortening homomorphisms of subsemigroups of a closed (hyperbolic) surface, but the same arguments apply to shorten homomorphisms
of subsemigroups of a punctured surface group together with Bass-Serre elements that conjugate its boundary elements. i.e., they apply
to the restrictions of the homomorphisms $\{\hat h_n\}$ to the pair $(U,\hat Q)$.

In order to prove proposition 5.4 we use the same
arguments that were used in the proof of theorem 442 to shorten homomorphisms in order to make the restrictions of the
homomorphisms, $\{\hat h_n\}$ to the subgroup $\hat Q$ longer.
By proposition 4.7  it is possible to find a sequence of tuples of elements $\{u^t_1,\ldots,u^t_g\}$ such that:
\roster
\item"{(1)}" the tuples
generate a limit group  $L_1$ that contains $\hat Q$. They
satisfy the properties that are listed in proposition 4.7.

\item"{(2)}" the restrictions of the homomorphisms $\hat h_n$ to the subgroup $\hat Q$ extend to homomorphisms $\tilde h_n:L_1 \to F_k$.
For a fixed $t$, for large enough index $n$, $\tilde h_n(u^t_i) \in FS_k$, $1 \leq i \leq g$.


\item"{(3)}" for each $t$, and each index $j$, $1 \leq j \leq f$, there exists a word $\hat z^t_j$, such that
$q_j=\hat z^t_j(u^t_1,\ldots,u^t_j)$.

\item"{(4)}" For every pair of indices $t_1,t_2$, and every $j$, $1 \leq j \leq f$, the elements
$\tilde h_n(\hat z^{t_1}_j(u^{t_2}_1,\ldots,u^{t_2}_j)) \in FS_k$ for large enough index $n$.
\endroster

Given an index $t>1$, we define an automorphism of the group $\hat Q$, $\psi_t$,
by sending the  generators $g_j$,
$1 \leq j \leq f$ to:
$\psi_t(q_j)=\hat z^t_j(u^1_1,\ldots,u^1_g)$ (note that
$\hat z^t_j(u^1_1,\ldots,u^1_g)$  is an element of $\hat Q$, and by construction this assignment defines an automorphism opf $\hat Q$r).

By part (4), for large enough $n$, $\tilde h_n \circ \psi_t(q_j) \in FS_k$, $1 \leq j \leq f$. Suppose that for every index $t$ we choose an
index $n_t$, such that $n_t$ grows to $\infty$, and $\tilde h_{n_t} \circ \psi_t(q_j) \in FS_k$, $1 \leq j \leq f$. A subsequence of the
sequence $\tilde h_{n_t} \circ \psi_t$ converges into an action of $\hat Q$ on some real tree. By construction and the properties
of the elements $u^t_1,\ldots,u^t_g$ (that are listed in proposition 4.7), this action is faithful, and the action
has to be bi-Lipschitz equivalent to the limit  action of $\hat Q$ on a real tree that is obtained from the restrictions of
the homomorphisms $\{\hat h_n\}$ to $\hat Q$. In particular, the limit action has to contain  an IET component that contains  (a finite index) IET action of $Q$ on a real tree.

The restrictions of the homomorphisms  $\tilde h_n \circ \psi_t$ to $\hat Q$ extend naturally to pair homomorphisms of the pair  $(\hat S, \hat L)$ (viewed as the fundamental
group of the abelian decomposition $\hat \Delta_2$). By construction, the
sequence of extended homomorphisms of $(\hat S,\hat L)$ restricts to a sequence of pair homomorphisms of $(S,L)$,
and this restricted sequence of homomorphisms has  a subsequence that converges into a faithful action of $L$ on a real
tree, and this action contains an IET component that is stabilized by the QH subgroup $Q$.
This concludes the proof of proposition 5.4.

\hfill$\square$

\medskip
After constructing the refined abelian decomposition $\Delta_2$, we use (part of) its modular group to
shorten the sequence of homomorphisms, $\{h_n\}$, that were used to construct the abelian decompositions
$\Delta_1$ and $\Delta_2$, to obtain a new abelian decomposition of the pair $(S,L)$. After repeating this
shortening procedure iteratively finitely many times, we are able to replace the "machine" that is used in
constructing the JSJ decomposition for groups [Ri-Se], and by the existence of such a machine,  eventually deduce the finiteness
of the maximal decompositions of the given pair $(S,L)$.

In the abelian decomposition $\Delta_2$ no non stably elliptic vertex group is connected to
only stably elliptic and periodic edge
groups, because otherwise we can further refine the abelian decomposition $\Delta_2$.

We fix a finite set of generators for each of the vertex groups in $\Delta_2$.
We order the vertex groups in $\Delta_2$ by the order of magnitude of the displacements of their fixed sets
of generators. Let $v_1,\ldots,v_{\ell}$ be the fixed set of generators of a vertex group $V$.
For each index $n$, we associate with $V$  the minimal displacement of a point in the Cayley graph $T$
of the coefficient group $F_k$, under the action of the tuple of elements $h_n(v_1),\ldots,h_n(v_{\ell})$.
We denote this minimal displacement $disp_n(V)$.

Given the sequence of homomorphisms $\{h_n\}$ we can not say that the displacement function of a vertex group $V_1$ is bigger than
the displacement function of a vertex group $V_2$, even up to a natural equivalence relation on displacement functions.
However, after passing to a subsequence of the homomorphisms $\{h_n\}$, and up to certain equivalence relation on the displacement functions,
the displacement functions define an order on the vertex groups in $\Delta_2$. We say that two displacement
functions, $disp_n(V_1),disp_n(V_2)$, are $comparable$ if there exists positive constants $c_1,c_2$, so
that for every index $n$, $c_1disp_n(V_1) \, < \, disp_n(V_2) \, < \, c_2 disp_n(V_1)$. We say that
$disp_n(V_1)$ dominates $disp_n(V_2)$ if $disp_n(V_2)=o(disp_n(V_1))$.

Since there are only finitely many vertex groups in $\Delta_2$, there is a finite collection of vertex groups
in $\Delta_2$ with comparable displacement functions, that dominate the displacement functions of all
the other vertex groups in $\Delta_2$. Note that none of the vertex groups in the dominating set
can be stably elliptic. By our construction of $\Delta_2$, for each vertex group in this dominating
subset of vertex groups, there exists a non-periodic non stably elliptic   edge group that is connected to it. Furthermore,
a non-periodic non stably elliptic edge group that is connected to a vertex group in the dominating set, must be
connected only to vertex groups in the dominating subset (since this edge group has a displacement function that is comparable with that of
a dominating vertex group).

Let $E$ be an edge in a decomposition (graph of groups) $\Lambda$ of a group $G$, and suppose that $E$ has an abelian edge group $A$.
A Dehn twist along the edge $E$ is obtained by conjugating by an element of $A$ the fundamental group of a connected component of $\Lambda \setminus E$ in
case $E$ is a bridge, and multiplying the Bass-Serre generator by an element of $A$ in case $\Lambda \setminus A$ is connected.

We set $Mod(\Delta_2)$ to be the modular group of the pair $(S,L)$ that is associated with the
abelian decomposition $\Delta_2$. We set $MXMod(\Delta_2)$ to be the subgroup of $Mod(\Delta_2)$,
that is generated by Dehn twists only along non-periodic (stably hyperbolic) edge groups for which their corresponding edges
connect between dominating vertex groups.

\noindent
In the graph of groups $\Delta_2$, the edges with non-periodic stably hyperbolic edge groups that
connect between dominating vertex groups, are grouped in several connected components, that we denote:
$\Gamma_1,\ldots,\Gamma_u$. We fix finite sets  of generators for the fundamental  groups
of each of the connected subgraphs, $\Gamma_1,\ldots,\Gamma_u$.

Recall that the abelian decompositions, $\Delta_1$ and $\Delta_2$ were obtained from a limit of a sequence of pair homomorphisms,
$h_n:(S,L) \to (FS_k,F_k)$. Each homomorphism, $h_n$, was chosen to be one of the homomorphisms $f_n(j)$, where the sequence
of homomorphisms $\{f_n(j)\}_{j=1}^{\infty}$ converges into a faithful  simplicial action of the limit group $L$ on a real tree with
an abelian decomposition $\Lambda_n$.

For each index $n$, we replace the homomorphism, $h_n$, with a homomorphism, $h_n^1$, that has similar properties to those
of $h_n$, and  choose
an automorphism $\varphi_n \in MXMod(\Delta_2)$, so that $h_n^1$, $\varphi_n$, and the twisted homomorphism,
$h_n^2 = h_n^1 \circ \varphi_n$, have the following properties:
\roster
\item"{(1)}" Let $T_n$ be the Bass-Serre tree that is associated with the abelian decomposition $\Lambda_n$.
$\varphi_n$ is chosen so that:

\itemitem{(i)} for large enough $j$, $f_n(j) \circ \varphi_n : L \to F_k$, is a homomorphism of the pair $(S,L)$ into the pair
$(FS_k,F_k)$ (i.e., $f_n(j) \circ \varphi_n$ sends the positive cone in $L$ into the standard positive cone $FS_k$ in $F_k$).

\itemitem{(ii)} for
each connected component, $\Gamma_i$, $i=1,\ldots,u$, there exists a point $t_i(n) \in T_n$, so that the displacement of
$t_i(n)$ under the action of  the fixed finite collection of generators of the fundamental group of $\Gamma_i$ twisted
by $\varphi_n$, is the shortest among all the points in $T_n$ and all the
automorphisms $\varphi \in MXMod(\Delta_2)$ that satisfy property (i).

\item"{(2)}" Having fixed $\varphi_n$, we choose $h^1_n$, to be one of the homomorphisms $f_n(j)$, for which both $h_n^1$ and
$h_n^1 \circ \varphi_n$ satisfy the conditions that $h_n$ was required to satisfy, with respect to the given action of $L$ on
$T_n$, $\lambda_n: L \times T_n \to T_n$,  and with respect to the twisted action $\lambda_n \circ \varphi_n: L \times T_n \to T_n$,
in correspondence. i.e., we require properties (2) and (3) that the homomorphisms $\{h_n\}$ satisfy.

We require
that both $h_n^1$ and $h_n^1 \circ \varphi_n$ will send the elements in the ball of
radius $n$ in the Cayley graph of $L$ monomorphically into $F_k$, and that both $h_n^1$ and $h_n^1 \circ \varphi_n$ will map  the
elements in a ball of radius $n$ in the Cayley graph of $L$ into elements in $F_k$ that approximate (up to dilatation) their images in
the tree $T_n$.
\endroster

The sequence of homomorphisms, $\{h^1_n\}$, converges into a faithful action of the limit group $L$ on the same real tree as the sequence,
$\{h_n\}$, i.e., to the action of $L$ on the real tree $Y$. With the action of $L$ on $Y$ we have associated the abelian decomposition,
$\Delta_2$. We set $h^2_n=h^1_n \circ \varphi_n$.
As the automorphisms $\varphi_n$ were chosen to be from the modular group $MXMod(\Delta_2)$, the sequence $h^2_n$ converges
into a faithful action of $L$ on a real tree $Y_1$.

With the action of $L$ on the real tree $Y_1$ we can naturally associate an abelian decomposition $\Delta_3$. Since we assumed that
all the maximal abelian decompositions that are associated with the pair $(S,L)$ are simplicial, the abelian decomposition $\Delta_3$
has to be simplicial as well. Starting with $\Delta_3$, we can possibly successively refine it and obtain an abelian
decomposition $\Delta_4$, in the same way that we successively refined the abelian decomposition $\Delta_1$ and obtained
the abelian decomposition $\Delta_2$.

\vglue 1pc
\proclaim{Lemma 5.5} $\Delta_4$ has the following properties:
\roster
\item"{(1)}" all the edges in $\Delta_2$ that are in the complement of the connected subgraphs $\Gamma_1,\ldots,\Gamma_u$
remain edges with the same edge groups in $\Delta_4$.

\item"{(2)}" it can be that all the vertex groups in $\Delta_4$ are stably elliptic. $\Delta_4$ can not have QH or abelian vertex groups that are
associated with axial components.

\item"{(3)}" the edge groups in the connected subgraphs $\Gamma_1,\ldots,\Gamma_u$ are not contained in vertex groups in $\Delta_4$.
Hence, there exists at least one edge group in $\Delta_2$ that is not elliptic in $\Delta_4$.
\endroster
\endproclaim

\nfp Since the automorphisms $\{\varphi_n\}$ act as the identity and with Dehn twists along edges that are in $\Delta_2$ and in
the complement of $\Gamma_1,\ldots,\Gamma_u$,  all these edges
remain edges in $\Delta_4$ and have the same edge groups as in $\Delta_2$.

Unlike the abelian decomposition $\Delta_2$, it may be that all the vertex and edge groups in $\Delta_4$ are stably elliptic.
Since we assumed that all the maximal abelian decompositions of the pair $(S,L)$ are simplicial,
proposition 5.4 implies that it can not be that
$\Delta_4$ contains an axial or a QH vertex group.

Since $\Delta_3$ is obtained from $\Delta_2$ by shortening edges that are stabilized by non-periodic non stably elliptic edge groups and
connect between dominating vertex groups, the edge groups in the connected subgraphs, $\Gamma_1,\ldots,\Gamma_u$, are
not contained in vertex groups in $\Delta_4$ (i.e., they must be hyperbolic in $\Delta_4$). Therefore, there exists
at least one edge group in $\Delta_2$ that is not elliptic in $\Delta_4$.

\hfill$\square$

The  abelian decomposition $\Delta_2$ of the limit group $L$ was
obtained from  an action of $L$ on a real tree $Y$, that is itself obtained as a limit of a sequence of
homomorphisms of $L$ into the free group $F_k$. Hence,  every edge group in $\Delta_2$ that is associated with a stabilizer of a
non-degenerate segment in the real tree $Y$, is either the centralizer of itself, or it is a direct summand in its centralizer.

We start by assuming that all the edge groups in $\Delta_2$ that have cyclic centralizers can be conjugated into vertex groups  in $\Delta_4$.
Let $E_1$ be an edge in $\Delta_2$ with a non-periodic (hence, non stably elliptic) abelian edge group $A_1$, so that the edge $E_1$ is contained in one of the
connected subgraphs of
groups $\Gamma_1,\ldots,\Gamma_u$  (so $E_1$ connects between dominating vertex groups, and  $A_1$ is not contained  in a vertex
group in $\Delta_4$).
By our assumption
the centralizer of $A_1$ in $L$ is non-cyclic.

First, we further assume that
$A_1$ is a strict direct summand in its centralizer. Since the action of $L$ on the real tree $Y$ is assumed to be
simplicial, and $A_1$ stabilizes pointwise a line in $Y$, $\Delta_2$ contains a circle so that all the edge groups in this
circle are stabilized by $A_1$, and the Bass-Serre generator that is associated with that circle, that we denote $t_1$, commutes
with $A_1$, and the centralizer of $A_1$ is the direct sum of $A_1$ and $<t_1>$.

Since the abelian decompositions, $\Delta_3$ and $\Delta_4$, were obtained by shortening the homomorphisms, $\{h_n^1\}$,
using the subgroup of modular automorphisms, $MXMod(\Delta_2)$, that shortens only edges with non-periodic edge
groups that connect between dominating vertex groups, the maximal abelian subgroup $A=A_1+<t_1>$, and in particular
its subgroup $A_1$,  do not
fix a vertex in $\Delta_4$. Hence, there is a circle in $\Delta_4$ so that all the edges in this circle are stabilized by
a direct summand $A_2<A$, and $A_1 \cap A_2$ is a direct summand of $A$ that is also a strict summand in both $A_1$ and $A_2$.
Therefore, there exists an element $t_2 < A_1$, such that $A_1= (A_1 \cap A_2) +<t_2>$.

If there are no vertex groups along the circle that is stabilized by $A_1$ in $\Delta_2$, or along the circle that
is stabilized by $A_2$ in $\Delta_4$, then the limit group $L$ has
to be a non-cyclic free abelian group. In that case, by theorem 4.1,  $L$ admits an axial action on a real tree that is obtained as a limit
of a sequence of homomorphisms of pairs from $(S,L)$ into $(FS_k,F_k)$, and such an action contradicts our assumption that
every action of $(S,L)$ that is obtained as a limit of a sequence of homomorphisms from $(S,L)$ into $(FS_k,F_k)$ is
simplicial.
Hence, we may assume that the circle that is stabilized by $A_1$ in $\Delta_2$ and by $A_2$ in $\Delta_4$ contain
non-trivial vertex groups.

Suppose that $A_1 \cap A_2$ is the trivial subgroup of $A$. Let $V$ be  the fundamental groups of a connected components of the
graph of groups that is obtained
from $\Delta_2$ by deleting the circle that is stabilized by $A_1$. $V$ being a subgroup of $L$ inherits an abelian decomposition from $\Delta_4$.
Since $V \, \cap \, A=A_1$, $V$ inherits a free decomposition from $\Delta_4$, and since a generator of $A_1$ can be taken to be a Bass-Serre generator
that corresponds to a circle in $\Delta_4$ that is stabilized by a subgroup of $A$ that intersects $V$ trivially, $V$ inherits a free decomposition
from $\Delta_2$ in which $A_1$ is a (free) factor.

Therefore, from the free decompositions of the fundamental groups of the connected components in the graph of groups that is
obtained from $\Delta_2$ by erasing the edges in the circle that is stabilized by $A_1$, $L$ admits a free decomposition in which
$A$ is a factor. This contradicts our assumption that $L$ is freely indecomposable.
Hence, we may
assume that $A_1 \cap A_2$ is a non-trivial direct summand of $A$.

Let $V_1,\ldots,V_{\ell}$ be the vertex groups in $\Delta_2$, that are placed in the circle that is stabilized by $A_1$ in $\Delta_2$.
Each of these vertex groups inherits a graph of groups decomposition from the abelian decomposition (of the ambient group $L$) $\Delta_4$.
Each of these graphs of groups of the vertex groups, $V_1,\ldots,V_{\ell}$, contains a circle, so that each edge in that circle is
stabilized by the direct summand $A_1 \cap A_2$.

The limit group $L$ is assumed to be freely indecomposable, and the centralizer of the edge group $A_1 \cap A_2$ is the non-cyclic
abelian subgroup $A$. Hence, by the existence of an abelian  JSJ decomposition for limit groups, an edge group $C$ in $\Delta_2$,
which is not conjugate to $A_1$, and is a subgroup of a vertex group $V_i$ (that is
contained in the circle that is stabilized by $A_1$ in $\Delta_2$), is elliptic with respect to a decomposition of $L$ along edges that
are stabilized by $A_1 \cap A_2$. Therefore, $C$, can be conjugated into a connected subgraph of the graph of  groups that
$V_i$ inherits from $\Delta_4$, where this subgraph does not contain any of the edges in the circle that is stabilized by $A_1 \cap A_2$
in the graph of groups that is inherited by $V_i$.

Each of the given set of generators $s_1,\ldots,s_r$ of the subsemigroup $S$ of the limit group $L$, can be written in a normal form with
respect to the graph of groups $\Delta_2$. These normal forms contain elements from the vertex groups $V_1,\ldots,V_{\ell}$ of $\Delta_2$. We
can further write each of the elements of $V_1,\ldots,V_{\ell}$, that appear in the normal forms of $s_1,\ldots,s_r$, in a normal form
with respect to the graph of groups that the vertex groups $V_1,\ldots,V_{\ell}$ inherit from the abelian decomposition
$\Delta_4$. If we substitute each of these last normal forms in the normal forms of  $s_1,\ldots,s_r$ with respect to $\Delta_2$, we
represented each of the elements $s_1,\ldots,s_r$ as (fixed) words in vertex groups and Bass-Serre generators of $\Delta_2$ and of
the vertex groups in the graphs of groups that are inherited by $V_1,\ldots,V_{\ell}$ from  $\Delta_4$. In particular, these fixed words contain powers of
the elements $t_1$ and $t_2$ that are contained in the abelian subgroup $A$, that are  Bass-Serre generators in $\Delta_2$ and $\Delta_4$,
in correspondence.

At this point we can finally modify the sequence of homomorphisms, $\{h^1_n\}$, in order to get a new sequence of
homomorphisms of the pair $(S,L)$ into the standard pair $(FS_k,F_k)$ that converges into a faithful action of the pair
$(S,L)$ on a real tree, and this action contains an axial component, a contradiction to our assumption that every
such limit action is simplicial.

With each of the elements, $s_1,\ldots,s_r$, we have associated a (fixed) word, that was constructed from a normal
form of the element $s_i$ with respect to the abelian decomposition $\Delta_2$, and normal forms of elements from the
vertex groups $V_1,\ldots,V_{\ell}$ with respect to the graphs of groups that these vertex groups inherit from
$\Delta_4$. These words contain powers of the elements $t_1$ and $t_2$.
We set $d$ to be the sum of the absolute values of powers of the element $t_2$ that appear in all the
(fixed) words that we have associated with the elements $s_1,\ldots,s_r$.

We choose $\ell-1$ positive  irrational numbers, $\alpha_1,\ldots,\alpha_{\ell-1}$, so that
$1,\alpha_1,\ldots,\alpha_{\ell-1}$ are independent over the rationals,  and
additional positive real numbers $\alpha_{\ell}$ and $\beta$ so that:
\roster
\item"{(1)}" $\alpha_1+\alpha_2+ \ldots +\alpha_{\ell}=1$.

\item"{(2)}" for every $i$, $1 \leq i \leq \ell$:
$$| \alpha_i - \frac {1} {\ell} | \, < \, \frac {1} {10  \ell}.$$

\item"{(3)}" $\frac {1} {4 d \ell} < \beta < \frac {1} {3 d \ell}$.

\item"{(4)}" $\beta$ is not in the subspace that is spanned by $1,\alpha_1,\ldots,\alpha_{\ell-1}$ over the rationals.
\endroster
Given the positive irrational numbers, $\alpha_1,\ldots,\alpha_{\ell}$ and $\beta$, we modify the sequence of homomorphisms
$\{h^1_n\}$. First, for every index $n$, we precompose the homomorphism $h^1_n$ with a modular automorphism
$\psi_n \in MXMod (\Delta_2)$. The automorphisms $\psi_n$ only conjugate the vertex groups $V_2,\ldots,V_{\ell}$,
and the vertex groups that are connected to them in $\Delta_2$ and are located outside the circle that is stabilized by $A_1$,
by different powers of the element $t_2$. The conjugations that determine the automorphisms $\{\psi_n\}$,
are chosen to guarantee
that in the limit tree that is obtained from the sequence of homomorphisms $\{h^1_n \circ \psi_n \}$, by using
the same rescaling constants as those that were used for the sequence $\{h^1_n\}$, the vertex groups $V_1,\ldots,V_{\ell}$
fix points, and the distances in the limit tree (that we denote) $\hat Y$ satisfy,
$d_{\hat Y}(Fix(V_i),Fix(V_{i+1}))=\alpha_i tr(t_1)$,
for every $i$, $i=1,\ldots,\ell-1$, and where $tr(t_1)$ is the displacement of the element $t_1$ along its axis
(that is stabilized by $A_1$) in $\hat Y$. Since the homomorphisms $\{h^1_n\}$ are homomorphisms of the pair $(S,L)$,
for large enough index $n$, the homomorphisms $h^1_n \circ \psi_n$ are homomorphisms of the pair $(S,L)$ as well, and
in the limit, we get a faithful action of the pair $(S,L)$ on the real tree $\hat Y$ (note that the action of $(S,L)$ on the
real tree $\hat Y$ is similar to the action of $(S,L)$ on the original real tree $Y$, where the difference is only in the lengths of
the segments between the points that are stabilized by the vertex groups $V_1,\ldots,V_{\ell}$).

At this point we further modify the homomorphisms $\{h^1_n \circ \psi_n\}$, by changing only the image of the element $t_2$.
Each of the vertex groups, $V_1,\ldots,V_{\ell}$, inherits an abelian decomposition from the graph of groups $\Delta_4$. In each
of these inherited graphs of groups there is a circle that is stabilized by the subgroup $A_1 \cap A_2$, and with a Bass-Serre
generator $t_2$. Hence, we can modify each of these graphs of groups of the vertex groups, $V_1,\ldots,V_{\ell}$, and replace the circle that
is stabilized by $A_1 \cap A_2$, with a new vertex group that is stabilized by $A_1$, that is connected to the vertex $V_{\ell}$ with an
edge that is stabilized by $A_1 \cap A_2$.

Therefore, each of the vertex groups, $V_1,\ldots,V_{\ell}$, can be written as an amalgamated product, $V_i = W_i \, *_{A_1 \cap A_2} \, A_1$.
Furthermore, all the edge groups in $\Delta_2$ of edges that are connected to the vertex that is stabilized by $V_i$ in $\Delta_2$ and
are not stabilized by $A_1$  can
be conjugated into $W_i$. This implies that the graph of groups $\Delta_2$ can be modified and collapsed to an amalgamated product:
$L=W*_{A_1 \cap A_2} \, A$, where $A$ is the centralizer of $A_1$ and $A_2$, and $A= (A_1 \cap A_2)+<t_2>+<t_1>$.

Now, we finally modify the homomorphisms $\{ h^1_n \circ \psi_n \}$ by changing the images of the element $t_2$. We define a sequence
of homomorphisms $\{h^3_n:(S,L) \to (FS_k,F_k)\}$ as follows. For each index $n$ and every element  $w \in W$ we set:
$h^3_n(w)=h^1_n \circ \psi_n (w)$. We further set $h^3_n(t_1)=h^1_n \circ \psi_n(t_1)$, and $h^3_n(t_2)$ to be an element that commutes
with $h^1_n(t_2)$, and so that in the limit tree $\tilde Y$,
that is obtained from the sequence of homomorphisms $\{h^3_n\}$, by using
the same rescaling constants as those that were used for the sequence $\{h^1_n\}$, $tr(t_2) = \beta tr(t_1)$. By our choice of the constants
 $\alpha_1,\ldots,\alpha_{\ell-1}$ and $\beta$, the homomorphisms $h^3_n$ map the fixed generators of the semigroup $S$, $s_1,\ldots,s_r$, to elements in the
standard free semigroup $FS_k$. Hence, $h^3_n$ are homomorphisms of the pair $(S,L)$ into the pair $(FS_k,F_k)$.

By construction, the subgroup $W$ acts faithfully on the limit tree $\tilde Y$, since it acts faithfully on the real tree $\hat Y$,
and the restrictions of the homomorphisms $h^3_n$ and $h^1_n$ to the subgroup $W$ are identical. Since
$\beta$ was chosen to be irrational, the subgroup $<t_1,t_2>$ acts indiscretely, with a dense orbit, on a line in $\tilde Y$. Since
the real numbers, $1,\alpha_1,\ldots,\alpha_{\ell-1},\beta$, are independent over the rationals, $L$ modulo the kernel of the action
of it on $\tilde Y$ can be written as an amalgamated product:   $W*_{A_1 \cap A_2} \, A$. Hence,
$L$ acts faithfully on the real tree $\tilde Y$. Finally, the real tree $\tilde Y$ was
constructed from a sequence of homomorphisms of the pair $(S,L)$ into $(FS_k,F_k)$, that we denoted $\{h^3_n\}$, and $\tilde Y$
 contains an axial component (the component
that contains the axis of $t_1$ and $t_2$), a contradiction to the assumption (of theorem 5.2) that every such limit action is discrete.

So far we assumed that all the edge groups in $\Delta_2$ that have cyclic centralizers can be conjugated into vertex groups in
$\Delta_4$,
and that there exists
an edge
in $\Delta_2$ with a non-periodic abelian edge group $A_1$, that is contained in one of the connected subgraphs, $\Gamma_1,\ldots,\Gamma_u$
(i.e., that this non-periodic edge group connects between dominating vertex groups),
and so that
$A_1$ is a (strict) direct summand in its centralizer. i.e., that the centralizer of $A_1$ is not contained in a vertex  in $\Delta_2$.

We have already pointed out that since we have started with an infinite sequence of inequivalent graphs of groups of $L$,
not all the edge groups in $\Delta_2$ are stably elliptic.
Suppose that all the edge groups in $\Delta_2$ that have cyclic centralizers, are not in the connected subgraphs of groups,
$\Gamma_1,\ldots,\Gamma_u$, so they  can all be conjugated into vertex groups in $\Delta_4$. Hence, the subgraphs $\Gamma_1,\ldots,\Gamma_u$
contain only edges with edge groups that have non-cyclic centralizers.
Suppose that there is no edge $E_1$ in $\Delta_2$, with an edge  group $A_1$,  for which:
\roster
\item"{(1)}" $E_1$  is contained in one of the subgraphs $\Gamma_1,\ldots,\Gamma_u$ (hence, $A_1$ has a non-cyclic centralizer).

\item"{(2)}" $A_1$ is a proper subgroup (a proper direct summand) in its centralizer.
\endroster
In this case $\Delta_4$ does contain edge groups with non-cyclic stabilizers, so that their centralizers can not be conjugated
into vertex groups in $\Delta_4$.
At this stage we can not apply the argument that we used previously and we continue as follows.

Suppose that all the edge groups in
$\Delta_4$ are stably elliptic. Since by our assumptions all the cyclic edge groups
in $\Delta_2$ are also cyclic edge groups in $\Delta_4$, the other edge groups in both $\Delta_2$ and $\Delta_4$ are abelian groups
with non-cyclic
centralizers. Hence, it is possible to construct a common refinement, $\Delta$, to both $\Delta_2$ and $\Delta_4$, and in $\Delta$
all the edge groups are stably elliptic. Now, since the edge groups in
$\Delta$ are stably elliptic, the common refinement $\Delta$ dominates all the abelian decompositions $\Lambda_n$ for large
enough $n$. Hence, the abelian decompositions $\{\Lambda_n\}$ belong to finitely many equivalence classes,
a contradiction to our assumptions.

Therefore at least one
edge group in $\Delta_4$ is not stably elliptic. Hence, we can start with $\Delta_4$, and apply the construction of $\Delta_3$ and $\Delta_4$,
to get new abelian decompositions of $L$, that we denote $\Delta_5$ and $\Delta_6$.

If there exists an edge $E$ in $\Delta_4$ with  an edge group $A$ that has a cyclic centralizer, so that $A$ is non-periodic and $E$
connects between two dominating vertex groups in $\Delta_4$ and hence not stably elliptic), $A$ can not be conjugated into a vertex group in $\Delta_6$.
In that case we show
in the sequel (proposition 5.6) that it's possible to find a sequence of homomorphisms of $(S,L)$ into $(FS_k,F_k)$
that converges into a faithful action
of $L$ on a real tree, and this action contains an axial or an IET component, a contradiction to the assumptions of theorem 5.2. Hence, we may
assume  that all the edges in $\Delta_4$ that have edge groups with cyclic centralizers, are either periodic or they don't connect
between dominating vertex groups.

If there exists an edge $E$ in $\Delta_4$ with an abelian edge group $A$ that has a non-cyclic centralizer,
the centralizer can not be conjugated into a vertex group
in $\Delta_4$, $A$ is non-periodic and connects between two dominating vertex groups in $\Delta_4$, we can apply the previous argument
to the edge group $A$ and its centralizer, and to the two abelian decompositions, $\Delta_4$ and $\Delta_6$, and obtain an action of $L$
on a real tree with an axial component, a contradiction to the assumptions of theorem 5.2.

Therefore, we may assume that  all the non-periodic edge groups in $\Delta_4$ that connect between dominating vertex groups,
have non-cyclic centralizers, and their centralizers can be conjugated into a vertex group in $\Delta_4$. In that case, there exists at least
one edge in $\Delta_4$ and at least 2 edges in $\Delta_6$, that have edge groups with non-cyclic centralizers, and these centralizers can not
be conjugated into vertex groups in $\Delta_4$ and $\Delta_6$ in correspondence.

By starting with $\Delta_6$ and iteratively constructing abelian decompositions, we either obtain:
\roster
\item"{(1)}"  an abelian decomposition in which all the edge and vertex groups are stably elliptic. In this case it is possible to
construct a common refinement of all the abelian decompositions that are constructed until the abelian decomposition in which all the edge
and vertex groups are stably elliptic. This abelian decomposition dominates all the abelian decompositions $\{\Lambda_n\}$ for large enough $n$. But
these decompositions were assumed to be inequivalent, and the constructed common refinement dominates only finitely many decompositions,
a contradiction.

\item"{(2)}" an edge $E$ with a stabilizer $A$ that connects between two dominating vertex group in $\Delta_{2i}$, such that
$A$ has non-cyclic centralizer, $A$ is non-periodic, and the centralizer of $A$ is not contained in a vertex group in $\Delta_{2i}$.
In that case we apply the construction that we used in case $\Delta_2$ has such an edge, and associate a sequence of
pair homomorphisms with the pair $(S,L)$ that converges into an action of $L$ on a real tree and the action contains an
axial component. A contradiction to the assumptions of theorem 5.2.

\item"{(3)}" an edge group with cyclic centralizer in $\Delta_{2i}$, so that this centralizer can not be conjugated into a vertex group
in $\Delta_{2(i+1)}$ - in that case we prove in the sequel (proposition 5.6)
that there exists a sequence of homomorphisms of the pair $(S,L)$
into $(FS_k,F_k)$ that converges to a faithful action of $(S,L)$ on a real tree, and this limit real tree contains an axial or an IET
component, a contradiction to the assumptions of theorem 5.2.

\item"{(4)}" if cases (1)- (3) do not occur, then the abelian decomposition $\Delta_{2i}$ contains at least $(i-1)$ edges with edge
groups that have non-cyclic centralizers, and these centralizers can not be conjugated into vertex groups in $\Delta_{2i}$.
\endroster
By the accessibility
for small splittings of a f.p.\ group [Be-Fe], or by acylindrical accessibility ([Se],[De],[We]), there is a global
bound on the number of edges in the abelian decompositions $\Delta_{2i}$,
that depends only on $L$ (in fact on the number of generators of $L$ [We]). Hence, if cases (1) - (3) do not occur, there is a global bound on the
number of steps that case (4) can occur.

Therefore, to complete the proof of theorem 5.2, we need the following basic proposition, that motivates our approach to the construction
of the JSJ decompositions for pairs.

\vglue 1pc
\proclaim{Proposition 5.6}
Suppose that for some index $i$, there exists an edge $E$ in $\Delta_{2i}$ with  an edge group $A$ that has a cyclic centralizer,
so that $A$ is non-periodic and $E$
connects between two dominating vertex groups in $\Delta_{2i}$, and $A$ can not be conjugated into a vertex group in $\Delta_{2(i+1)}$.
Then there exists a sequence of homomorphisms of $(S,L)$ into $(FS_k,F_k)$ that converges into a faithful action
of $L$ on a real tree, and this action contains  an IET component.
\endproclaim

\nfp By the properties of the JSJ decomposition of the limit group $L$, the subgroup $A$ corresponds to a s.c.c.\ on some maximal
QH subgroup $MSQ$ in the JSJ decomposition of $L$. Hence, to prove proposition 5.6 it is possible to use geometric properties of
(oriented) curves and arcs on  surfaces (i.e., on the surface that is associated with $MSQ$), together with
properties of  the simplicial actions
of $L$ on the real trees that are associated with the abelian decompositions $\Delta_{2i}$ and $\Delta_{2(i+1)}$.

The abelian decomposition $\Delta_{2i}$ was obtained from a limit of a sequence of homomorphisms of pairs that we denote:
$\{h_n:(S,L) \to (FS_k,F_k)\}$. The abelian decomposition $\Delta_{2(i+1)}$ is obtained from a sequence of
shortened homomorphisms, $\{h_n^s:(S,L) \to (FS_k,F_k)\}$, i.e., from the homomorphisms $\{h_n\}$ that are pre-composed
with automorphisms in the modular group
that is associated with $\Delta_{2i}$ (generated by Dehn twists along the edges in $\Delta_{2i}$), $h_n^s=h_n \circ
\varphi_n$, where $\varphi_n \in MXMod(\Delta_{2i})$.

Since the sequence of  shortened homomorphisms, $\{h_n^s\}$, converges into an action of $L$ on a real tree with an associated
abelian decomposition, $\Delta_{2(i+1)}$, after possibly passing to a subsequence,
each of the shortened homomorphisms, $h^s_n$, can be written as a composition:
$h_n^s=h_n^b \circ \nu_n$, where $h_n^b$ is obtained from $h_n^s$ by shortening it with automorphisms from
$MXMod(\Delta_{2(i+1)})$ and :
\roster
\item"{(1)}"  $\nu_n \in Mod(\Delta_{2(i+1)})$.

\item"{(2)}" $h_n^b$ is a pair homomorphism: $h_n^b:(S,L) \to (FS_k,F_k)$.

\item"{(3)}" for each index $n$: $$\max_{1 \leq j \leq r}  length(h_n^s(s_j)) \, \geq  \,
n \cdot \max_{1 \leq j \leq r}  length(h_n^b(s_j)).$$
\endroster

Let $A_1,\ldots,A_t$ be the edge groups in $\Delta_{2i}$ that are non-periodic, have cyclic stabilizers, and can not be
conjugated into a vertex group in $\Delta_{2i+1}$. Let $a_1,\ldots,a_t$ be the positive generators of $A_1,\ldots,A_t$, i.e.,
the generators of the cyclic groups $A_1,\ldots,A_t$ that have conjugates that are mapped to $FS_k$, by each of the homomorphisms
$\{h_n\}$.

With the limit group $L$ there is an associated JSJ decomposition. Since the subgroups $A_1,\ldots,A_{\ell}$ have cyclic centralizers, and they can not be conjugated into vertex groups in $\Delta_{2(i+1)}$, their positive generators $a_1,\ldots,a_t$ correspond to
oriented s.c.c.\  $C_1,\ldots,C_t$ on surfaces (that are associated with $MQH$ subgroups) in the JSJ decomposition of $L$.

Each of edge subgroups $A_1,\ldots,A_t$ can not be conjugated into a vertex group in $\Delta_{2(i+1)}$. Hence, with each
of the edge group $A_1,\ldots,A_t$ we can associate a minimal, non-degenerate subgraph of $\Delta_{2(i+1)}$, that
contains a conjugate of it.
Let $B_1,\ldots,B_{\ell}$ be the edge groups in the union of these minimal non-degenerate  subgraphs. By the properties of the
JSJ decomposition of the limit group $L$, each of the subgroups $B_1,\ldots,B_{\ell}$ is cyclic, has a cyclic centralizer, and
is associated with (the subgroup generated by) a s.c.c.\ on a surface in the JSJ decomposition of $L$ (as a limit group).

Let $b_1,\ldots,b_{\ell}$ be the positive generators of the cyclic groups $B_1,\ldots,B_{\ell}$. These positive
generators correspond to
oriented s.c.c.\  $c_1,\ldots,c_{\ell}$ on surfaces in the JSJ decomposition of the limit group $L$. Furthermore, since
the generators $a_1,\ldots,a_t$ and $b_1,\ldots,b_{\ell}$ are positive, if two oriented curves $C_i$ and $c_j$ intersect, each
of their intersection points is positively oriented.

With each of the given generators $s_1,\ldots,s_r$, being elements in the limit groups $L$ that can be represented in a normal form with
respect to the JSJ decomposition of $L$, we can naturally associate a (possibly empty) finite
collection of  oriented arcs and curves on the surfaces in the JSJ decomposition of the limit group $L$. Since the
graph of groups $\Delta_{2i}$ and $\Delta_{2(i+1)}$ were constructed from pair homomorphisms, the elements
$s_1,\ldots,s_r$ act positively on the real trees that are associated with both $\Delta_{2i}$ and $\Delta_{2(i+1)}$. Hence, if an oriented arc or curve
from the finite collections that
are associated with each of the elements $s_1,\ldots,s_r$ intersect non-trivially with any of the oriented s.c.c.\ $C_1,\ldots,C_t$ or
$c_1,\ldots,c_{\ell}$, then every such intersection point is positively oriented.

\smallskip
With each
homomorphism from a (punctured) surface group into a free group it is possible to associate a map from the
(punctured) surface into a bouquet of circles. With a pair homomorphism the image of such a map is a bouquet of
oriented circles. By homotoping such a map to be transversal at the midpoints of the (oriented) circles, we may assume that
the preimages of the midpoints of the oriented circles are finite collections of oriented arcs and s.c.c.\ on the (punctured)
surface.

The pair homomorphisms $h_n^s$ are compositions $h_n^s=h_n^b \circ \nu_n$, where $\nu_n \in Mod(\Delta_{2(i+1)})$,  and the
homomorphisms $h_n^b$ are much shorter than the homomorphisms $h_n^s$. With each of the homomorphisms $h_n^b$ there is an associated
map from each of the (punctured) surfaces in the JSJ decomposition of $L$ into a bouquet of oriented circles. As we indicated, we
may assume that the preimages of the midpoints of these oriented circles are finite collections of oriented s.c.c.\ and arcs.
Up to conjugacy,
each of the oriented s.c.c.\ $c_j$ (that are associated with the positive generators $b_1,\ldots,b_{\ell}$ of edge groups in
$\Delta_{2(i+1)}$) are mapped by $h_n^b$ into $FS_k$. Hence, we may homotope the map from the surface to the bouquet of
oriented circles
so that the preimages of the midpoints of the oriented circles are oriented s.c.c.\ and arcs that intersect the s.c.c.\
$c_1,\ldots,c_{\ell}$ only in positive orientation.

We set the (modular) automorphism $\psi \in MXMod(\Delta_{2i})$ to be the composition of positive Dehn twists along the edges
that are associated with the edge groups $A_1,\ldots,A_t$ in $\Delta_{2i}$, where positive Dehn twists refers to Dehn twists along the oriented s.c.c.\
$C|_1,\ldots,C_t$ that are associated with $A_1,\ldots,A_t$. We set the (modular) automorphism
$\varphi \in MXMod(\Delta_{2(i+1)})$ to be the composition of positive Dehn twists along the edges that are
associated with the edge groups $B_1,\ldots,B_{\ell}$.

We look at homomorphisms of the form $h_n^b \circ \varphi^{m}$ for large $m$. Since the homomorphisms $h_n^s$ map the
elements $s_1,\ldots,s_r$ and $C_1,\ldots,C_t$ into $FS_k$, so do the homomorphisms $h_n^b \circ \varphi^m$. We look
at the map from the surfaces that are associated with the JSJ decomposition of $L$ into the bouquet of
oriented circles that is associated with homomorphisms $h_n^b \circ \varphi^m$. The preimages of the midpoints of the oriented
circles under this map, are obtained from the preimages of the midpoints of these circles under the maps that are associated with
the homomorphisms $h_n^b$, by performing $m$ powers of Dehn twists along the s.c.c.\ $c_1,\ldots,c_{\ell}$.

Since the homomorphisms $h_n^s$ were obtained as shortenings, we may assume that $h_n^s \circ \psi^v (C_i) \in FS_k$
and $h_n^s \circ \psi^v(s_j) \in FS_k$, for
$1 \leq i \leq t$, $1 \leq j \leq r$, and $v=0,1,\ldots,v_0$ for some fixed (previously chosen) positive integer $v_0$. Hence,
for every every large $m$,
$h_n^b \circ \varphi^m \circ \psi^v (C_i) \in FS_k$
and $h_n^s \circ \varphi^m \circ \psi^v(s_j) \in FS_k$, for
$1 \leq i \leq t$, $1 \leq j \leq r$, and every positive integer $v$.

Because  the  intersection points between the s.c.c.\ $c_i$ and $C_j$, $1 \leq i \leq t$, $1 \leq j \leq \ell$, are all
positively oriented, and because the high power $\varphi^m$ of positive Dehn twists along the curves $\{c_j\}$,
surrounds each such s.c.c.\ with long positive (periodic) words, it follows that
for every large positive pair $m_1,m_2$,  and for every large $m$ that is much bigger than both $m_1$ and $m_2$:
$h_n^b \circ \varphi^m \circ \psi^{m_2} \circ \varphi^{m_1}(c_j) \in FS_k$, $1 \leq j \leq \ell$, and
$h_n^b \circ \varphi^m \circ \psi^{m_2} \circ \varphi^{m_1}(C_i) \in FS_k$, $1 \leq i \leq t$.

Now, with each of the elements $s_1,\ldots,s_r$ we can associate oriented arcs and curves on the surfaces that are associated
with the JSJ decomposition of $L$. Whenever these curves  intersect  the curves $C_1,\ldots,C_{\ell}$ and $c_1,\ldots,c_t$, they
intersect them positively. Since $h_n^b(s_j) \in FS_k$, and $h_n^s(s_j) \in FS_k$, for $1 \leq j \leq r$, for every large
positive pair $m_1,m_2$,  and for every large $m$ that is much bigger than both $m_1$ and $m_2$:
$h_n^b \circ \varphi^m \circ \psi^{m_2} \circ \varphi^{m_1}(s_j) \in FS_k$, $1 \leq j \leq r$.

Furthermore, the positive intersection numbers, and the high power of positive Dehn twists along the curves $\{c_j\}$,
 imply that for every positive integer $p$, and every
large tuple of positive integers: $e_1,f_1,\ldots,e_p,f_p$,
and for every large $m$ that is much bigger than the sum of these positive integers:
$$h_n^b \circ \varphi^m \circ \psi^{f_p} \circ \varphi^{e_p} \circ
\psi^{f_1} \circ \varphi^{e_1}(C_i) \in FS_k \ , \ 1 \leq i \leq t$$
$$h_n^b \circ \varphi^m \circ \psi^{f_p} \circ \varphi^{e_p} \circ
\psi^{f_1} \circ \varphi^{e_1}(c_j) \in FS_k \ , \ 1 \leq j \leq \ell$$
$$h_n^b \circ \varphi^m \circ \psi^{f_p} \circ \varphi^{e_p} \circ
\psi^{f_1} \circ \varphi^{e_1}(s_j) \in FS_k \ , \ 1 \leq j \leq r.$$

Therefore, the sequence of powers that is used in constructing the JSJ decomposition for groups (theorem 4.5 in [Ri-Se]),
when taken to consist of
only positive powers (so that the twisted homomorphisms are indeed pair homomorphisms),
can be used to construct a sequence of homomorphisms that converges into a faithful action of
the limit group $L$ on a real tree, and this limit action contains an IET component, as proposition 5.6 claims.

\hfill$\square$

Proposition 5.6 completes the proof of theorem 5.2.

\hfill$\square$

\medskip
Theorem 5.2 proves that if every faithful action of a
freely indecomposable  limit group $L$ on a real tree $Y$, that is obtained  as a limit from a sequence of homomorphisms of the
pair $(S,L)$ into $(FS_k,F_k)$, is discrete (or simplicial), then the pair $(S,L)$ has only finitely many maximal abelian
decompositions, that we view as its (finite collection of) canonical JSJ decompositions. A similar statement is valid in case
every such action contains only simplicial and axial components.

\vglue 1pc
\proclaim{Theorem 5.7} Suppose that $L$ is freely indecomposable, and that
all the maximal abelian decompositions, $\{\Lambda_i\}$, that are associated with the pair $(S,L)$,
correspond to faithful actions of $L$ on a real tree, where these actions contain only simplicial and axial components.

Then
there exist only finitely many (equivalence classes of)
maximal abelian decompositions of the pair $(S,L)$.
\endproclaim

\nfp Suppose that there are infinitely many (equivalence classes of) maximal abelian decompositions of
a pair $(S,L)$. Let $\{\Lambda_i\}_{i=1}^{\infty}$ be the collection of these  (inequivalent)
maximal decompositions.

All the abelian decompositions, $\Lambda_i$, are dominated by the JSJ decomposition of the freely indecomposable
limit group $L$. If there exists an infinite subsequence of the decompositions, $\{\Lambda_i\}$, that is
dominated by an abelian decomposition, $\Theta_1$,
 that is strictly dominated by the JSJ decomposition of $L$, we pass
to that subsequence. We continue iteratively. If the infinite subsequence of the decompositions, $\{\Lambda_i\}$,
that is dominated by $\Theta_1$, has a further infinite subsequence that is dominated by an abelian decomposition,
$\Theta_2$, that is strictly dominated by $\Theta_1$, we pass to that subsequence.
Since every strictly decreasing sequence of abelian decompositions of the limit
group $L$ terminates after finitely many steps, the sequence of maximal abelian decompositions of the pair,
$(S,L)$, contains an infinite subsequence, that is dominated by an abelian decomposition, $\Theta$, and
no infinite subsequence of that infinite subsequence is dominated by an abelian decomposition, $\Theta'$, that
is strictly dominated by $\Theta$. We (still) denote this infinite subsequence of maximal abelian decompositions,
$\{\Lambda_i\}$, and their dominating abelian decomposition of $L$, $\Theta$.

Our goal is to show that the abelian decomposition, $\Theta$, that strictly dominates the entire sequence of maximal
abelian decompositions of the pair $(S,L)$, $\{\Lambda_i\}$, is itself a maximal abelian decomposition of $(S,L)$,
a contradiction to the maximality of each of the abelian decompositions, $\Lambda_i$. Hence, a contradiction to the
existence of an infinite sequence of maximal abelian decompositions that are associated with the pair $(S,L)$, and the
theorem follows.

We start with the infinite sequence of inequivalent maximal abelian decompositions $\{\Lambda_i\}$, that are all dominated by the
abelian decomposition $\Theta$.
With the sequence of abelian decompositions, $\{\Lambda_i\}$, we have associated a sequence
of homomorphisms $\{h_n\}$ that satisfy properties (1)-(3) that are listed in the beginning of the proof of theorem 5.2. By possibly
passing to a subsequence of the homomorphisms $\{h_n\}$ (still denoted $\{h_n\}$), we obtain a convergent sequence, that
converges into a superstable action of the limit group $L$ on a real tree $Y$, and with this action there is an
associated abelian decomposition of $L$ that we denoted $\Delta_1$.

Starting with $\Delta_1$ we (possibly) refine it to get an abelian decomposition $\Delta_2$, precisely as we did in case all the
actions on real trees are simplicial, i.e., precisely as we did in the course of proving theorem 5.2. Since all
the abelian decompositions $\Lambda_i$ are dominated by $\Theta$, $\Delta_2$ is dominated by $\Theta$ as well.

\proclaim{Proposition 5.8} If $\Delta_2$ is equivalent to the abelian decomposition $\Theta$, then there exists a sequence
of homomorphisms $u_n:(S,L) \to (FS_k,F_k)$ that converges into a superstable
 action of the limit group $L$ on a real tree $\hat Y$, and the abelian decomposition that is dual to this action is $\Theta$
itself.
\endproclaim

\nfp
Note that $\Delta_2$ was constructed after  iterative refinements  from the sequence of homomorphisms $\{h_n\}$, and it is
suppose to be obtained from the limit of the convergent sequence $\{u_n\}$ without any refinements.

$\Theta$ can not be simplicial since it dominates an infinite sequence of inequivalent maximal abelian decompositions
of $L$. Hence, $\Theta$ contains vertex groups that are associated with
either axial or IET components. $\Delta_2$ that is assumed to be equivalent to $\Theta$ is obtained from a sequence of pair
homomorphisms: $\{h_n:(S,L) \to (FS_k,F_k)\}$. If $\Theta$ does not contain vertex groups that are associated with IET components,
we vary the homomorphisms $\{h_n\}$ by precomposing them with Dehn twists along edge groups (that connect between non-axial
vertex groups), and vary the values of a prefered basis of generators of vertex groups that are associated with axial
components using theorem 4.1, as we did in the proof of proposition 5.4 in case $\Delta_2$ contains an axial component.
Such modifications gives a new sequence of homomorphisms $\{\hat h_n\}$ that converges into a faithful action of the limit
group $L$ on a real tree, where the abelian decomposition that is associated with the limit action is $\Delta_2$ that is
equivalent to $\Theta$.

If $\Delta_2$ contains vertex groups that are associated with IET components, then we modify the homomorphisms $\{h_n\}$
by performing Dehn twists along edges that connect between non-axial non-QH vertex groups, vary the values of prefered
basis of axial vertex groups using theorem 4.1, and precomposing the obtained homomorphisms with automorphisms that are
extensions of automorphisms of the QH
vertex groups - precisely as we did in proving proposition 5.4 (in case $\Delta_2$ contains a QH vertex group). Again, the
obtained sequence of homomorphisms converges into a faithful action of $L$  on a real tree, where the abelian decomposition that
is associated with this action is $\Delta_2$ that is assumed to be equivalent to $\Theta$.

\hfill$\square$

Since $\Theta$ dominates all the maximal abelian decompositions, $\{\Lambda_i\}$, proposition 5.8 implies that
if $\Delta_2$ is equivalent to $\Theta$ we
obtained a contradiction, and theorem 5.2 follows. Hence, we may assume that $\Theta$ strictly dominates $\Delta_2$.

We continue in a similar way to what we did in proving theorem  5.2.
If all the edge groups in $\Delta_2$ are stably elliptic, and the only non stably elliptic vertex groups are  axial,
then  for almost all the indices $i$, the abelian
decomposition $\Lambda_i$ is dominated by the abelian decomposition $\Delta_2$ (i.e., $\Delta_2$ is a
(possibly trivial) refinement of $\Lambda_i$ for almost every index $i$). Hence, there is a subsequence of the maximal
abelian decompositions, $\{\Lambda_i\}$, that is dominated by the abelian decomposition $\Delta_2$, that is strictly dominated
by the abelian decomposition $\Theta$, a contradiction to our choice of $\Theta$. Therefore, $\Delta_2$ must contain
a stably hyperbolic edge group.

Suppose that  $\Delta_2$ contains an axial vertex group $A=A_1+A_2$, so that $\Delta_2$
collapses to the abelian decomposition $L=L_1*_{A_1} \, A$. If in addition $A_1$ is stably elliptic, then it is elliptic in almost
all the abelian decomposition $\{\Lambda_i\}$, so it must be elliptic in the dominating abelian decomposition $\Theta$, and
therefore, $A$ is an abelian vertex group in $\Theta$ as well. If $A_1$ is not stably elliptic, it can not be elliptic in $\Theta$, and
therefore, $rk(A_1) \geq 2$. In both cases we continue by essentially analyzing the subgraph that is obtained from $\Delta_2$
by taking out the vertex group $A$, and continue with the fundamental group of the obtained graph of groups that we denote $L_1$.

Conversely, if $A$ is an axial vertex group in $\Delta_2$, $A=A_1+A_2$, so that $\Delta_2$ collapses to an abelian decomposition
$L=L_1*_{A_1} \, A$, and $A$ is conjugate to an abelian vertex group in $\Theta$,
that is connected to other vertex groups in $\Theta$ by subgroups that generate a conjugate to $A_1$, then $A_1$ is elliptic in all
the abelian decompositions $\{\Lambda_i\}$.

From the graph of groups $\Delta_2$ we take out the edges with stably elliptic edge groups,
and the axial vertex groups. What left are several connected subgraphs of $\Delta_2$,
that we denote: $C_1,\ldots,C_u$.
For each of these subgraphs we fix a finite generating set of its fundamental group.

As in the proof of theorem 5.2, we shorten the homomorphisms $\{h_n\}$ (using Dehn twists along dominating edge groups),
and pass to a subsequence of them
for which the shortenings converge into a (new) faithful action of the limit group $L$ on some real tree.
We further restrict each of the shortened homomorphisms  to the fundamental groups
of the connected components (that are subgraphs in $\Delta_2$), $C_i$, $i=1,\ldots,u$. After possibly passing to a subsequence,
the restrictions of the shortened homomorphisms converge to a faithful action
of that group on a real tree, that we denote $Y_i$. We denote the abelian decomposition that is associated with these actions,
$\Delta_3^i$, $i=1,\ldots,u$. Given each of the abelian decompositions $\Delta_3^i$ we further refine it in the same way in which
we refined the abelian decomposition $\Delta_1$ to obtain the abelian decomposition $\Delta_2$. We denote the obtained refined
abelian decompositions, $\Delta_4^i$, $1 \leq i \leq u$.

The edges that were taken out from $\Delta_2$ to obtain the connected components, $C_i$, $i=1,\ldots,u$, are either stably elliptic
or contained in abelian subgroups of rank at least 2 in the fundamental groups of the connected
components $C_i$. Hence, these edge groups can be assumed to be either elliptic in the
abelian decomposition, $\Delta_4^i$, $i=1,\ldots,u$, or the abelian decompositions $\Delta_4^i$ can be modified to guarantee
that all the non-cyclic abelian subgroups are elliptic, so that the edge
groups that were taken out from $\Delta_2$ are elliptic.
Therefore, from the abelian decomposition $\Delta_4^i$, $1 \leq i \leq u$,
possibly after a modification that guarantees that all the non-cyclic abelian subgroups are elliptic,
it is possible to obtain an abelian decomposition of the ambient limit group $L$, by adding the elliptic edge groups, and
the QH and (axial) abelian vertex groups that were taken out from $\Delta_2$. We denote the obtained abelian
decomposition $\Delta_4$.

Since the abelian decomposition $\Theta$ strictly dominates $\Delta_2$, the construction of  $\Delta_4$ from (almost)
shortest homomorphisms
of the connected components $C_i$ of $\Delta_2$ implies that at least one of the following possibilities must hold:
\roster
\item"{(1)}" there exists a QH vertex group in $\Delta_4$.

\item"{(2)}" there exists an edge group with a cyclic centralizer
in $\Delta_4$ that can not be conjugated into an edge group nor into a vertex group in $\Delta_2$.
Equivalently, there is an edge group with cyclic centralizer in $\Delta_2$ that
can not be conjugated into a  vertex group or into an edge group in $\Delta_4$.

\item"{(3)}" there exists a non-trivial abelian subgroup $A_1$ of an abelian subgroup $A$ of rank at least 2 in $L$,
so that $A_1$ was
contained in a non-QH, non-axial vertex group in $\Delta_2$ and $A$ was not contained in such a vertex group in
$\Delta_2$, and $A_1$ can not be conjugated into a non-QH, non-axial vertex group in any of the abelian decompositions $\Delta_4^i$.

\item"{(4)}" there exists a non-trivial abelian subgroup $A$ of rank at least 2 in $L$, so that $A$ was
contained in a non-QH, non-axial vertex group in $\Delta_2$, and $A$
can not be conjugated into a  vertex group in any of the abelian decompositions $\Delta_4^i$.
\endroster

If case (1) occurs we get a contradiction to the assumptions of theorem 5.7 according to
proposition 5.4. If cases (3) or (4) occur we
apply the argument that was used in the proof of theorem 5.2, and obtain a proper refinement of $\Delta_2$. If case (2) occurs
it is possible to apply the argument that was used to prove proposition 5.6 and construct a sequence of homomorphisms of the pair $(S,L)$ that converge into a faithful action of
$L$ on a real tree, where this real tree contains an IET component, a contradiction to the assumptions of theorem 5.7.

As long as the obtained abelian decomposition is strictly dominated by $\Theta$ we can continue refining the obtained abelian
decomposition iteratively. Therefore, after finitely many refinements we get the abelian decomposition $\Theta$, and
a sequence of homomorphisms of the pair $(S,L)$ that converge into a faithful action of $L$ on a real tree with an associated
abelian decomposition $\Theta$. Since $\Theta$ dominates the sequence of maximal abelian decompositions $\{\Lambda_i\}$, we
obtained a contradiction to their maximality, and theorem 5.7 follows.

\hfill$\square$

\smallskip
In case $L$ is freely indecomposable, and all the faithful actions of $L$ on a real tree that are obtained as limit of sequences
of homomorphisms of the pair $(S,L)$ contain only simplicial and axial components, we set each of the finitely many maximal
abelian decompositions that are associated with such a pair
$(S,L)$, to be a $JSJ$ $decomposition$ of the pair $(S,L)$. Hence,
with such a pair $(S,L)$ we have canonically assigned finitely many (abelian) JSJ decompositions, that
are all obtained from the abelian JSJ decomposition of the limit group $L$, by possibly cutting some
of the QH vertex groups along finitely many s.c.c.\ and further collapsing and folding.

In case  of limit pairs $(S,L)$ in which the the limit group $L$ is freely indecomposable, and the abelian decompositions that are
associated with limit actions of $L$ that are obtained from convergent sequences of pair homomorphisms of $(S,L)$ contain QH
vertex groups the situation is more involved.

In that case we don't associate finitely many abelian decompositions with the pair $(S,L)$ that can be viewed as its JSJ
decomposition. Using the techniques that appear in later sections in this work and in [Se8] it is possible to associate
with $(S,L)$ a finite collection of finite hierarchical sequences of decompositions of $L$. These finite collection of hierarchical
sequences is canonical and can be viewed as the JSJ decomposition of a pair $(S,L)$ in which $L$ is freely indecomposable.



\vglue 1.5pc
\chapter{Limit sets}
\fullwidthtext
\medskip
One of the fundamental objects that is associated with a Kleinian group is its limit set. For such a
(non-elementary) group the
limit set is a non-empty closed subset of $S^2$, that often has a fractal structure. In this section we define
a natural limit set that is associated with a pair $(S,L)$, where $L$ is a freely indecomposable limit group,
 and one of its associated (canonical) JSJ decompositions.

Given a (Gromov) hyperbolic  group or a limit group one can study all the small stable faithful
actions of such a group on a real
tree. For a freely indecomposable torsion-free hyperbolic group, the properties of its canonical JSJ decompositions,
together with
a theorem of Skora [Sk], and Thurston's compactification of the Teichmuller space, imply that the set of
small stable faithful
actions of such a group on a real tree up to dilatation,  can be naturally viewed as a topological space, that is
homeomorphic to a finite union of (possibly trivial) products of finitely many spheres and a (possibly trivial) Euclidean factor. The
properties of the abelian decomposition of a limit group, implies the same conclusion for the structure
of the set of small stable
faithful actions of a freely indecomposable limit group on a real tree up to dilatation.

\vglue 1.5pc
\proclaim{Definition 6.1} Let
$(S,L)$ be a pair in which $L$ is a freely indecomposable limit group.
Look at all the  sequences of pair homomorphisms $\{h_n:(S,L) \to (FS_k,F_k)\}$ that converge into a faithful action of
$L$ on a real tree. Note that the action must be small and (super) stable.

The set of all
these limit actions up do dilatation, is a subset of all the faithful small stable actions of $L$ on real trees up to dilatation,
and can naturally be equipped with the induced topology. We call this topological space the $limit$ $set$ of the pair
$(S,L)$.
\endproclaim

It is natural to ask if such a limit set is homeomorphic to a finite simplicial complex. If so it may be described more
explicitly, perhaps even computed. In the sequel, for the purposes of encoding the structure of all pair homomorphisms
of a given pair, we don't need a detailed understanding of this limit set, apart from some of its global properties.

If it is homeomorphic to a simplicial complex, when is freely indecomposable, one can ask the same question for general
limit groups $L$, omitting the indecomposability requirement.

\vglue 1.5pc
\chapter{Enlarging the positive cone}
\fullwidthtext
\medskip

The JSJ decomposition of a freely indecomposable limit group plays an essential role in
the geometric construction of the Makanin-Razborov diagram that encodes the set of solutions to a system
of equations over a free group ([Se1]). In section 5 we managed to associate canonically finitely
many JSJ decompositions with pairs $(S,L)$ only under additional assumptions. Therefore, to construct a Makanin-Razborov
diagram that encodes the solutions to a system of equations over a free semigroup, we need to modify the construction
that is used over groups.

The main object that we use in the case of semigroups is a $resolution$. A resolution has finitely many steps,
it starts  with the given pair $(S,L)$, and continues with quotients of it that are not necessarily proper quotients.
Later on we show that there exist finitely many such
resolutions that encode all the pair homomorphisms of a given pair $(S,L)$.

To construct resolutions we start with the $extended$ $cone$ that we associate with an abelian decomposition, that was
already used in the shortenings in section 4. We start with a strengthening of theorem 4.1 which is not needed in the sequel,
but is of independent interest.

\vglue 1pc
\proclaim{Theorem 7.1} Let $A$ be a f.g.\ free abelian group, let $rk(A)=\ell$, and let
$S$ be a f.g.\ subsemigroup of $A$ that generates
$A$ as a group. Let $s_1,\ldots,s_r$ be a fixed generating set of $S$.

Then there exist finitely many $positive$ collections of (free) bases of $A$, $a^i_1,\ldots,a^i_{\ell}$, $1 \leq i \leq c$,  such that
for any sequence of homomorphisms $h_n:(S,A) \to (FS_k,F_k)$, that converges into a free action of $A$ on
a real tree,
there exists
an index $i$, $1 \leq i \leq c$, and an index $n_0$, for which for every index $n>n_0$:
\roster
\item"{(1)}" $h_n(a^i_j) \in FS_k$  for
$1 \leq j \leq \ell$.

\item"{(2)}" for each of the (fixed set of) generators $s_1,\ldots,s_r$ of the semigroup $S$,
there are fixed words, that depend only on the index
$i$, $1 \leq i \leq c$,  so that for every $m$, $1 \leq m \leq r$:
$$s_m \, = \, w_{i,m}(a^i_1,\ldots,a^i_{\ell})$$
where the words $w_{i,m}$ are (fixed) positive words in the elements $a^i_1,\ldots,a^i_{\ell}$.
\endroster
\endproclaim

\nfp  The theorem is immediate if $\ell=1$, hence, we may assume that $\ell>1$. Given a sequence of homomorphisms
$\{h_n:(S,A) \to (FS_k,F_k)\}$, that converges into a free action of $A$ on a real tree, theorem 4.1 implies that
there exists a free basis of $A$, $a_1,\ldots,a_{\ell}$,
and an index $n_0$, so that for every  $n>n_0$, the homomorphisms, $\{h_n\}$,  satisfy properties (1) and (2)
in the statement of the theorem.

To prove that there exists a finite collection of positive words that suffices
for all the convergent  sequences we use a compactness argument.
We have already shown that given a sequence of homomorphisms $\{h_n:(S,A) \to (FS_k,F_k)\}$ that converges into a free action
of $A$ on a real tree $Y$, there exists a collection of positive words and a basis of $A$ for which the conclusions
of the theorem hold.

Suppose that finitely many such collections of free bases together with positive words do not suffice. There exist only
countably many possible
collections of positive words and free bases, so we order the infinite set of the collections of free bases and positive words that are associated
with convergent sequences of homomorphisms that satisfy the assumptions of the theorem.

For each positive integer $t$, let $\{h_n^t:(S,A) \to (FS_k,F_k)\}$ be a sequence of homomorphisms that satisfies the assumption
of the theorem (i.e., it converges to a free action of $A$ on a real tree).
Suppose further that each of the homomorphisms in the sequence $\{h_n^t\}$ does not satisfy
the conclusions of the theorem with respect to the first $t$ collections of free bases and positive words  from the ordered list.

For each positive integer $t$, the sequence $\{h_n^t\}$ converges into a free action of the (free) abelian group
$A$ on an oriented  line $Y_t$,  $\lambda_t:A \times Y_t \to Y_t$. From the sequence of actions $\{\lambda_t\}$ it is a possible
to extract a convergent subsequence, that converge into a non-trivial, not necessarily faithful, action of the abelian group $A$
on a directed real line, $\alpha_0:A \times R_0 \to R_0$.

If the action $\alpha_0$ is not free, let $A_1$ be the kernel of the action $\alpha_0$.
$A_1$ is a direct summand in $A$. If $A_1$ is non-trivial, then we can pass to a
further subsequence, for which the sequence of actions $\lambda_t$ restricted to the direct summand $A_1$ converges into
a non-trivial, not necessarily faithful action  of $A_1$ on a directed line:
$\alpha_1:A_1 \times R_1 \to R_1$.

If the action $\alpha_1$ is not free, we continue iteratively by restricting the actions $\lambda_t$ to the kernel of the
previous action (which is a direct summand in the previous summand), and pass to a convergent subsequence. Since $A$ is a f.g.\ free
abelian group this iterative procedure terminates after finitely many steps.

Let $\{\lambda_{t_d}\}$ be the final convergent subsequence that converges to an action of the abelian group $A$ on a directed line $R$. Hence,
there exists a convergent sequence of pair homomorphisms
$h_{t_d}:(S,A) \to (FS_k,F_k)$ that converges to the same limit action as the limit action of the sequence $\{\lambda_{t_d}\}$, i.e. to an
action of $A$ on the directed line $R$,
such that the homomorphism $h_{t_d}$ does not satisfy the conclusion of the theorem with respect to the first $t_d$ collections of free bases and
positive words from the ordered list.

By the argument that
was used to prove theorem 4.1, the conclusions of theorem 4.1 hold for the convergent sequence $\{h_{t_d}\}$. Hence, there exists a
fixed free basis of $A$ that act positively on $R$ for large enough $d$, and a fixed set of words $\{w_m\}$, $1 \leq m \leq r$,
that satisfy part (2) in the statement of the theorem.

This free basis of $A$, together with the words $\{w_m\}$, appear in the ordered list of bases and words. Hence, for large
enough index $d$, the  homomorphisms $\{h_{t_d}\}$ satisfy part (2) in the conclusion of the theorem with respect to a fixed free basis and
a collection of words $\{w_m\}$ that appear in the ordered list. But the choice of the homomorphisms $\{h_{t_d}\}$ assumes that
the homomorphism $h_{t_d}$ does not satisfy part (2) in the conclusions of the theorem for all the first collections of free bases
and words $\{w_m\}$ from the ordered list, and we get a contradiction.

Therefore, finitely many
free bases together with a finite collection of words suffice for all the convergent sequences that satisfy the assumption of the
theorem, and the conclusion of the theorem follows.

\hfill$\square$

For a pair $(S,Q)$ in which $Q$ is a (closed) surface group, $S$ is a f.g.\ subsemigroup that generates $Q$ as a group, and
$(S,Q)$ is obtained as a limit of a sequence of pair homomorphisms into $(FS_k,F_k)$, we state a weaker statement that associates
a $standard$ $cone$ with such a pair.

\vglue 1pc
\proclaim{Lemma 7.2} Let $(S,Q)$ be a pair of a closed (hyperbolic) surface  group $Q$ and a subsemigroup $S$
that generates $Q$ as a group. Let $s_1,\ldots,s_r$ be a fixed generating set of $S$.

There exist finitely many residually free groups, $RF_1,\ldots,RF_m$, with embeddings $\nu_i: Q \to RF_i$, $1 \leq i \leq m$, and with
standard IET presentations, where the
generators of $RF_i$ are $q^i_1,\ldots,q^i_{\ell_i}$.

For every sequence of pair homomorphisms:
$\{h_n:(S,Q) \to (FS_k,F_k)\}$  that converges into a free action of
$Q$ on a real tree, there exists an index $i$, $1 \leq i \leq m$, such that for large enough index $n$, the homomorphisms $\{h_n\}$ extend to homomorphisms
$\{ \hat h_n:RF_i \to F_k \}$, where $h_n= \hat h_n \circ \nu_i$, and $\hat h_n (q^i_b) \in FS_k$, $ 1 \leq b \leq \ell_i$.
\endproclaim

\nfp Since the sequence of homomorphisms $\{h_n\}$ converges into a minimal IET action on a real tree,
by the Rips machine, or the Makanin procedure, there exists a standard set of generators of one of finitely many residually free groups $RF_i$,
with embeddings $\nu_i: Q \to RF_i$ that satisfies the
conclusion of the lemma (cf. theorem 4.2).

\hfill$\square$

Theorem 7.1 and lemma 7.2 generalize to pairs $(S,L)$ in which the ambient limit group $L$ is freely indecomposable. The
generalization is crucial in constructing the Makanin-Razborov diagram  of a pair.

\vglue 1pc
\proclaim{Theorem 7.3} Let $(S,L)$ be a pair of a freely indecomposable limit group $L$ and a
subsemigroup $S$ that generates $L$ as a group. Let $s_1,\ldots,s_r$ be a fixed generating set of $S$. Let
$\{h_n:(S,L) \to (FS_k,F_k)\}$ be a sequence of pair homomorphisms that converge into a faithful action of $L$ on a real tree $Y$
with an associated abelian decomposition $\Lambda$. Suppose that $A$ is the stabilizer of an axial component in the real tree $Y$
(and a vertex group in $\Lambda$), and $Q$ is a QH vertex group in $\Lambda$ that is associated with an IET component.

Let $A_0<A$ be the direct summand of $A$ that contains the subgroup that is generated by the edge
groups that are connected to the vertex that is stabilized by $A$ in $\Lambda$ as a subgroup of finite index. Suppose that:
$rk(A)-rk(A_0)=\ell$. Let $\hat g$ be the number of vertices in $\Lambda$, that are not stabilized by the abelian group $A$, for
which their vertex groups have conjugates that stabilize points in the axis of $A$ (these are all the vertices that are adjacent
to the vertex that is stabilized by $A$ in $\Lambda$). Let $g=max(\hat g -1,0)$. Let $EA=A+<e_1,\ldots,e_g>$ be a free abelian
group of rank $rk(A)+g$.

With the QH subgroup $Q$ we associate a natural extension $\hat Q$.  $\hat Q$
is generated by $Q$ and finitely many elements that map a point with non-trivial stabilizer in the IET component in $Y$  that is
stabilized by $Q$, to another point with a non-trivial stabilizer in that subtree that is not in the same orbit under the
action of $Q$. By adding these elements all the points with non-trivial stabilizers in the IET component in $Y$ that is stabilized by
$Q$ are in the same orbit under the natural extension of $Q$.

We replace the graph of groups $\Lambda$ with a graph of groups $\hat \Lambda$, by enlarging the vertex group $A$ to be $EA$ and
adding edges that are stabilized by boundary subgroups of the vertex group $Q$ and Bass-Serre generators, that make sure that
all the point stabilizers in the IET component that is stabilized by $Q$ in the real tree $Y$ are in the same orbit. We denote the
fundamental group of the graph of groups $\hat \Lambda$ by $\hat L$.

Each of the given generators of the semigroup $S$, $s_1,\ldots,s_r$,  that generate $L$ as a group, can be written as a word (in a normal form)
w.r.t.\ the graph of groups $\hat \Lambda$. These normal forms contains
 elements that lie in $\hat Q$, the natural extensions of the surface group $Q$, that is generated by $Q$ and the Bass-Serre generators that we added
to $\hat \Lambda$, and are associated
with the subpaths of the paths that are associated with $s_1,\ldots,s_r$ that are contained in subtrees that are stabilized
by conjugates of $Q$ in the real tree $Y$, and elements that lie outside these extensions.
Let $t_1,\ldots,t_d$ be the elements from $\hat Q$  that appear as subwords in the normal
forms  of the elements $s_1,\ldots,s_r$ w.r.t.\ the graph of groups $\hat \Lambda$..

Then there exist:
\roster
\item"{(i)}" $\ell+g$ elements that are part from a (free) basis of a free abelian group $EA$, $a_1,\ldots,a_{\ell+g}$,
such that $EA=<a_1,\ldots,a_{\ell+g}>+A_0$. The elements $a_1,\ldots,a_{\ell}$ are part of a free basis of $A$, $A=A_0+<a_1,\ldots,a_{\ell}$.
The elements $a_{\ell+1},\ldots,a_{\ell+g}$ are elements that are added to $EA$ so that all the non-trivial point stabilizers on the axis of
$A$ in $L$ will belong to the same orbit.

\item"{(ii)}" a residually free group $RF$ with an embedding $\nu: \hat Q \to RF$ and a standard IET presentation and generators,  $q_1,\ldots,q_m$.
The subgroups $\hat Q$ that contains the QH vertex group $Q$,  and the residually free group $RF$, satisfy properties that are analogous to the ones that are listed in lemma 7.2.
 The elements $t_1,\ldots.t_d$ are equal to words in the generators $q_1,\ldots,q_m$.

For large index $n$, the restrictions of the homomorphisms $h_n$ to the QH vertex group $Q$ extend to homomorphisms $\hat h_n: RF \to F_k$.

\item"{(iii)}" each of the generators $s_j$ can be written as a word in terms of elements $\{t_d\}$, that are elements in $RF$,
elements $\{u^j_e\}$ that are positive words in
$a_1,\ldots,a_{ell+g}$, and elements in the other vertex and edge groups in $\Lambda$.
\endroster

Such that
there exists an index $n_0$ for which for every $n>n_0$ the
following properties hold:
\roster
\item"{(1)}"
$h_n(a_b) \in FS_k$  for
$1 \leq b \leq \ell+g$.

\item"{(2)}" $\hat h_n(q_b) \in FS_k$, $1 \leq b \leq m$,  and
$\hat h_n(t_d) \in FS_k$.
\endroster
\endproclaim

\nfp Follows by the same arguments that were used to prove theorem 7.1 and lemma 7.2.

\hfill$\square$

\vglue 1.5pc
\chapter{A Makanin--Razborov diagram -- the freely indecomposable case}
\fullwidthtext
\medskip

In the previous section we associated a standard set  of generators with an abelian decomposition. In this section we use
shortenings that were constructed in section 2, together with the machine for the construction of the JSJ
decomposition for groups, to associate finitely many resolutions with a given pair $(S,L)$ in which the limit group
$L$ is freely indecomposable. These resolutions enable one to encode the set of all the pair homomorphisms of such a pair
using the set of pair homomorphisms of finitely many proper quotient pairs.

For presentation purposes we will start by proving the main theorem in case the limit group $L$ is freely indecomposable and
contains no no
non-cyclic abelian subgroups, and then combine it with arguments that were used in the construction of the JSJ
decompositions of pairs (in special cases) in section 3 to omit the assumption on abelian subgroups.
In the sequel we generalize the construction of
resolutions to all pairs, omitting the freely indecomposable assumption.
These resolutions are the building block of the Makanin-Razborov diagram that encodes all the pair homomorphisms of a given
pair, or alternatively, all the solutions to a system of equations over a free semigroup.

\vglue 1pc
\proclaim{Theorem 6.1}
Let $(S,L)$ be a pair, where $L$ is a freely indecomposable limit group, and let $s_1,\ldots,s_r$ be a fixed generating
set of the semigroup $S$. Suppose that the limit group $L$ contains no non-cyclic abelian subgroups.
Let $\{h_n:(S,L) \to (FS_k,F_k)\}$ be a sequence of pair homomorphisms that converges into
a faithful action of $L$ on a real tree $Y$.

Then there exists a $resolution$:
$$(S_1,L_1) \to  (S_2,L_2) \to \ldots \to (S_m,L_m) \to (S_f,L_f)$$
that satisfies the following properties:
\roster
\item"{(1)}" $(S_1,L_1)=(S,L)$, and  $\eta_i:(S_i,L_i) \to (S_{i+1},L_{i+1})$ is an isomorphism for $i=1,\ldots,m-1$ and
$\eta_m:(S_m,L_m) \to (S_f,L_f)$ is a
proper quotient map.

\item"{(2)}" with each of the pairs $(S_i,L_i)$, $1 \leq i \leq m$, there is an associated abelian decomposition that we denote
$\Lambda_i$.

\item"{(3)}" there exists a subsequence of the homomorphisms $\{h_n\}$ that factors through the resolution. i.e., each homomorphism
$h_{n_r}$ from the subsequence, can be written in the form:
$$h_{n_r} \ = \ \hat h_r \, \circ \, \varphi^m_r \, \circ \, \ldots \, \circ \, \varphi^1_r$$
where $\hat h_r:(S_f,L_f) \to (FS_k,F_k)$, each of the automorphisms $\varphi^i_r \in Mod(\Lambda_i)$, and
each of the homomorphisms:
$$h^i_{n_r} \ = \ \hat h_r \, \circ \, \varphi^m_r \, \circ \, \ldots \, \circ \, \varphi^i_r$$
is a pair homomorphism $h^i_{n_r}: (S_i,L_i) \to (FS_k,F_k)$.
\endroster
\endproclaim

\nfp If the abelian decomposition $\Lambda$ that is associated with the action of $L$ on the limit tree $Y$ is equivalent to
the abelian JSJ decomposition of the limit group $L$ (as a group), then the theorem  follows from the possibility to shorten
using $Mod(\Lambda)$ that was proved in section 2, and claim 5.3 in [Se1]. Note that in theorem 2.2 we used an embedding of a QH vertex group
into a residually free group and then used automorphisms of the bigger residually free group to shorten homomorphisms. However, the
automorphisms of the residually free group that we used preserve the QH subgroup, so they reduce to automorphisms of
the QH vertex group.

In the case in which $\Lambda$ is the JSJ of the freely indecomposable limit group $L$, by theorem 5.3 in [Se1]  every shortening quotient is
a proper quotient, so a subsequence of shortenings of the homomorphisms $\{h_n\}$ converge into a proper quotient of the
pair $(S,L)$, and the conclusion of the theorem follows with a resolution of length 1: $\eta_1:(S,L) \to (S_f,L_f)$.

Hence, we may assume that the abelian decomposition $\Lambda$, that is associated with the action of $L$ on the limit tree
$Y$, is strictly dominated by the abelian JSJ decomposition of $L$. In particular, $\Lambda$ does not correspond to a
closed surface, so it contains edges.


Assume that the sequence $\{h_n\}$ converges into a faithful action of the limit group $L$ on some real tree $Y$
with an associated abelian decomposition $\Lambda$,
$L$ is freely indecomposable, and contains no non-cyclic abelian subgroups.

We start by refining the abelian decomposition $\Lambda$ precisely as we did in proving theorem 3.2. Using definition 3.3
we divide the edges that are adjacent to a given non-QH vertex group in $\Lambda$ to periodic and non-periodic. If there
is a non-QH vertex group in $\Lambda$ that is adjacent only to periodic edge groups, we pass to a subsequence of the homomorphisms
$\{h_n\}$ for which their restrictions to this vertex group converges to a non-trivial faithful action on a real tree. The abelian
decomposition that is associated with this action (in which all the edge groups that are connected to this vertex group in $\Lambda$
are periodic), enable us to further refine $\Lambda$. By Bestvina-Feighn accessibility [Be-Fe1], or alternatively by
acylindrical accessibility [Se3],
this refinement process terminates after finitely many steps, and we obtain an abelian decomposition $\Lambda_1$. Note that in
$\Lambda_1$ every non-QH vertex group is adjacent to at least one non-periodic edge group.

In section 3 we divided the edges and vertices in $\Lambda_1$ into several equivalence classes,
according to the growth of the translation
lengths of fixed sets of generators of the corresponding edge groups. We fix a finite set of generators for each of the
vertex groups and edge groups
in $\Lambda_1$.
Recall that we say that two displacement
functions of two vertex groups or edge groups,
$disp_n(V_1),disp_n(V_2)$, are $comparable$ if there exists positive constants $c_1,c_2$, so
that for every index $n$, $c_1 \cdot disp_n(V_1) \, < \, disp_n(V_2) \, < \, c_2 \cdot disp_n(V_1)$. We say that
$disp_n(V_1)$ $dominates$ $disp_n(V_2)$ if $disp_n(V_2)=o(disp_n(V_1))$.

After passing to a subsequence  of the homomorphisms $\{h_n\}$ this defines an order on the equivalence classes of the edge
groups and vertex groups
in $\Lambda_1$. In particular, there exists a collection of edge groups and vertex groups
with comparable displacement functions that
dominate all the other (displacement functions of) vertex groups and edge groups in $\Lambda_1$.

Suppose first that $\Lambda_1$ does not contain QH vertex groups. In that case
we set $Mod(\Lambda_1)$ to be the modular group of the pair $(S,L)$ that is associated with the
abelian decomposition $\Lambda_1$, and  $MXMod(\Lambda_1)$ to be the subgroup of $Mod(\Lambda_1)$,
that is generated by Dehn twists only along dominating edge groups.

We start by using the full modular group $Mod(\Lambda_1)$. For each index $n$, we set the
pair homomorphism: $h_n^1:(S,L) \to (FS_k,F_k)$ and
$h_n^1=h_n \circ \varphi_n$ where $\varphi_n \in Mod(\Lambda_1)$, to be the shortest pair homomorphism
that is obtained from $h_n$ by precomposing it with a modular automorphism from $Mod(\Lambda_1)$. If
there exists a subsequence of the homomorphisms $\{h_n^1\}$ that converge into a proper quotient of the pair
$(S,L)$, we set the limit of this subsequence to be $(S_f,L_f)$, and the conclusion of the theorem follows
with a resolution of length 1.

Therefore, we may assume that every convergent subsequence of the homomorphisms $\{h_n^1\}$
converges into a faithful action of the limit group $L$ on some real tree. In that case we use only Dehn twists along
dominant edge groups. For each index $n$, we set $h_n^1=h_n \circ \varphi_n$, where $\varphi_n \in MXMod(\Lambda_1)$,
and $h_n^1$ is one of the
shortest homomorphisms that is obtained by precomposing $h_n$ with a modular automorphism from $MXMod(\Lambda_1)$. We
pass to a subsequence of the homomorphisms $\{h_n^1\}$ that converges into (necessarily faithful) action of $L$ on some
real tree with an associated abelian decomposition $\Delta_2$. We further refine $\Delta_2$ by analyzing actions of non-QH vertex
groups in $\Delta_2$ that are connected only to periodic edge groups and obtain a (possibly) refinement of $\Delta_2$ that
we denote $\Lambda_2$.

For presentation purposes we assume that $\Lambda_2$ contains no QH vertex groups. Note that by construction every
dominant edge group in $\Lambda_1$ is not elliptic in $\Lambda_2$.
With $\Lambda_2$ we associate its modular group $Mod(\Lambda_2)$, and the modular group that is associated
only with dominant edge groups that we denote $MXMod(\Lambda_2)$.

We continue iteratively. First we shorten using the full modular group, $Mod(\Lambda_i)$, and check if there a subsequence of
shortened homomorphisms that converge into a proper quotient of the pair $(S,L)$. If there is such a subsequence we obtained a
finite resolution. If not we shorten only along edges with dominant edge groups, i.e., using automorphisms from $MXMod(\Lambda_i)$,
pass to a convergent subsequence and
further refine the obtained abelian decomposition. For presentation purposes we assume that all the obtained abelian decompositions,
$\{\Lambda_i\}$ do not contain QH vertex groups.

\vglue 1pc
\proclaim{Definition 6.2} Let $\Lambda_1,\ldots,\Lambda_i,\ldots$ be the infinite sequence of  abelian decompositions of the
limit groups $L$ that are constructed along the iterative procedure. Note that we assume that all these abelian decompositions
do not contain QH vertex groups, and that the limit group $L$ in the pair $(S,L)$ is freely indecomposable and contains no
non-cyclic abelian subgroups. By construction, a dominant edge group in $\Lambda_i$ is not elliptic in $\Lambda_{i+1}$.

All the abelian decompositions $\Lambda_i$ are dominated by the abelian JSJ decomposition of $L$. For each index $i$, we set
$\Theta_i$ to be the minimal abelian decomposition (w.r.t. to the natural partial order that we defined on abelian decompositions
of $L$) that dominates all the abelian decompositions, $\Lambda_i,\Lambda_{i+1},\ldots$.

The sequence of abelian
decompositions  $\{\Theta_i\}$ is non-increasing, and every strictly decreasing sequence of abelian
decompositions of $L$ has to terminate,
hence,
there exists some index $i_0$, such that for every $i>i_0$, $\Theta_{i_0}=\Theta_i$. We call $\Theta_{i_0}$ the $stable$ $dominant$
abelian decomposition of the sequence $\Lambda_1,\Lambda_2,\ldots$.
\endproclaim

The abelian decompositions
$\Lambda_i$ were assumed to have no  QH vertex groups.
Since every dominant edge group in $\Lambda_i$ is hyperbolic in $\Lambda_{i+1}$,
each of the dominant abelian decompositions, $\Theta_i$, dominates a pair of hyperbolic-hyperbolic splittings of
the limit group $L$. Hence, it must contain a QH vertex group.
In particular, the stable
dominant abelian decomposition, $\Theta_{i_0}$, must contain a QH vertex group.

The next proposition is crucial in our approach to proving theorem 6.1.
It enables the substitution of the entire suffix of the sequence of
abelian decompositions $\{\Lambda_i\}$,
$\Lambda_{i_0},\Lambda_{i_0+1},\ldots$, with a single abelian decomposition - the stable dominant decomposition, $\Theta_{i_0}$.

\vglue 1pc
\proclaim{Proposition 6.3} There exists a subsequence of shortened pair homomorphisms $\{h^i_{n_i}\}$, $i \geq i_0$, and
a sequence of automorphisms: $\psi^i, \nu^i \in Mod (\Theta_{i_0})$, $i \geq  i_0$, with the following properties:
\roster
\item"{(1)}" $h^i_{n_i}$ is obtained from the pair homomorphism $h_{n_i}$, by shortening $h_{n_i}$ using a sequence of elements from
 the dominant modular groups:
$MXMod(\Lambda_1),\ldots,MXMod(\Lambda_{i})$. i.e.,
$$h_{n_i}\ = \ h^i_{n_i} \, \circ  \, \varphi^{i}_{n_i} \, \circ \, \ldots \, \circ \, \varphi^1_{n_i}$$
where $\varphi^j_{n_i} \in MXMod(\Lambda_j)$, $1 \leq j \leq i$  and for each $j$, $1 \leq j \leq i-1$, the homomorphism:
$$h^j_{n_i}\ = \ h^i_{n_i} \, \circ  \, \varphi^{i}_{n_i} \, \circ \, \ldots \, \circ \, \varphi^{j+1}_{n_i}$$
is a pair homomorphism $h^j_{n_i}:(S,L) \to (FS_k,F_k)$.

\item"{(2)}" For each index $i$ we define a pair homomorphism:
$$f_i\ = \ h^i_{n_i} \, \circ  \, \nu^{i+1} \, \circ \psi^i \, \circ \, \ldots \, \circ \, \psi^{i_0}.$$
The sequence of pair homomorphisms, $\{f_i\}$, converges into a faithful action of the limit group
$L$ on a real tree with an associated abelian decomposition $\Theta_{i_0}$. In particular, the limit action contains an IET
component.
\endroster
\endproclaim

\nfp To prove the proposition we basically imitate the construction of the JSJ decomposition as it appears in
[Se1] and [Ri-Se]. We fix a generating set $<s_1,\ldots,s_r>$ of the semigroup $S$, that also generate the limit group
$L$ (as a group).
We start with the subsequence of pair homomorphisms that are obtained using a (finite) iterative
sequence of shortenings from a subsequence of the given sequence of homomorphisms $\{h_n\}$, using the dominant
modular groups $MXMod$ that are associated with the abelian decompositions $\Lambda_1,\ldots,\Lambda_{i_0-1}$. By our assumptions
the sequence $\{h^{i_0-1}_{n_t}\}$ converges into a faithful action of $L$ on a real tree, with an associated abelian decomposition
that can be further refined (by restricting the homomorphisms to non-QH vertex groups that are connected only to
periodic edge groups) to the abelian decomposition $\Lambda_{i_0}$.

We set the automorphisms $\nu^{i_0}_{n_t} \in Mod (\Lambda_{i_0})$ to be automorphisms
that guarantee that the compositions $h^{i_0-1}_{n_t} \circ \nu^{i_0}_{n_t}$ converge into a faithful action of $L$ on a real
tree $Y_{i_0}$, in which:
\roster
\item"{(1)}" the abelian decomposition that is associated with the action of $L$ on the new real tree $Y_{i_0}$
is $\Lambda_{i_0}$. In particular, all
the edges in $\Lambda_{i_0}$ correspond to non-degenerate segments in $Y_{i_0}$.

\item"{(2)}" for some pre-chosen (and arbitrary) positive constant $\epsilon_{i_0}$, and after rescalings,
the lengths of the segments $[y_0,s_j(y_0)]$, $1 \leq j \leq r$, differ by at most
$\epsilon_{i_0}$ from the length of that segment in the original real tree, which is the limit of the sequence of homomorphisms
$h^{i_0-1}_{n_t}$. Furthermore, the length of each segment in the quotient graph of groups that is associated with the
action of $L$ on $Y_{i_0}$ differ by at most $\epsilon_{i_0}$ from the length of the corresponding segment in the graph
of groups that is associated with the action of $L$ on the real tree that is obtained as a limit of the sequence
$h^{i_0-1}_{n_t}$.

\item"{(3)}" the segments $[y_0,s_j(y_0)]$ are degenerate or positively oriented, and for each $t$,
$h^{i_0-1}_{n_t} \circ \nu^{i_0}_{n_t}(s_j) \in FS_k$.
\endroster
By the properties of the abelian decomposition $\Lambda_{i_0}$, such automorphisms $\nu^{i_0}_{n_t}$ can be constructed, by choosing appropriate Dehn twists
along the edges in $\Lambda_{i_0}$.

We continue iteratively. In analyzing the next level, the one in which the sequence $\{h^{i_0}_{n_t}\}$ converges
into a limit action on a real tree that is associated with $\Lambda_{i_0+1}$,
we set the automorphisms $\psi^{i_0} \in MXMod(\Lambda_{i_0})$
and $\nu^{i_0+1}_{n_t} \in Mod(\Lambda_{i_0+1})$,
to be automorphisms
that guarantee that the compositions $h^{i_0}_{n_t} \circ \nu^{i_0+1}_{n_t}$
converge into a faithful action of $L$ on a real
tree $Y_{i_0+1}$, in which:
\roster
\item"{(1)}" the abelian decomposition that is associated with the action of $L$ on the new real tree $Y_{i_0+1}$
is $\Lambda_{i_0+1}$. In particular, all
the edges in $\Lambda_{i_0+1}$ correspond to non-degenerate segments in $Y_{i_0+1}$.

\item"{(2)}" we fix finite generating sets for all the vertex groups in $\Lambda_{i_0}$ that are connected
to dominant edge groups, and for all the dominant edge groups in $\Lambda_{i_0}$.

For some pre-chosen (and arbitrary) positive constant $\epsilon_{i_0+1}$, after rescalings,
the lengths and the translation lengths of all the elements in these fixed finite generating sets when acting on
$Y_{i_0+1}$, differ by
at most $\epsilon_{i_0+1}$
from the corresponding lengths and translation lengths in the limit tree that is obtained from the sequence of
homomorphisms $\{h^{i_0}_{n_t}\}$.

\item"{(3)}" after rescaling,
the lengths of the segments $[y_0,\psi^{i_0}(s_j)(y_0)]$, $1 \leq j \leq r$, in $Y_{i_0+1}$, differ by at most
$\epsilon_{i_0+1}$ from the lengths of the corresponding segments, $[y_0,s_j(y_0)]$ in the limit tree $Y_{i_0}$.

\item"{(4)}" the segments $[y_0,s_j(y_0)]$ are degenerate or positively oriented in $Y_{i_0+1}$, and for each $t$,
$h^{i_0}_{n_t} \, \circ \, \nu^{i_0+1}_{n_t}(s_j) \in FS_k$ and
$h^{i_0}_{n_t} \, \circ \, \nu^{i_0+1}_{n_t} \, \circ \, \psi^{i_0}(s_j) \in FS_k$.
\endroster

In analyzing the next level, the one in which the sequence $\{h^{i_0+1}_{n_t}\}$ converges
into a limit action on a real tree that is associated with $\Lambda_{i_0+2}$,
we consider all the elements in the ball of radius 2, $B_2$, in the Cayley graph
of the limit group $L$ w.r.t. the generating set $s_1,\ldots,s_r$.

We set the automorphisms $\psi^{i_0+1} \in MXMod(\Lambda_{i_0+1})$
and $\nu^{i_0+2}_{n_t} \in Mod(\Lambda_{i_0+2})$
to be automorphisms
that guarantee that the compositions $h^{i_0+1}_{n_t} \circ \nu^{i_0+2}_{n_t}$
converge into a faithful action of $L$ on a real
tree $Y_{i_0+2}$, in which:
\roster
\item"{(1)}" the abelian decomposition that is associated with the action of $L$ on the new real tree $Y_{i_0+2}$
is $\Lambda_{i_0+2}$. In particular, all
the edges in $\Lambda_{i_0+2}$ correspond to non-degenerate segments in $Y_{i_0+2}$.

\item"{(2)}" we fix finite generating sets for all the vertex groups in $\Lambda_{i_0+1}$ that are connected
to dominant edge groups, and for all the dominant edge groups in $\Lambda_{i_0+1}$.

For some pre-chosen (and arbitrary) positive constant $\epsilon_{i_0+2}$, after rescalings,
the lengths and the translation lengths of all the elements in these fixed finite generating sets when acting on
$Y_{i_0+2}$, differ by
at most $\epsilon_{i_0+2}$
from the corresponding lengths and translation lengths in the limit tree that is obtained from the sequence of
homomorphisms $\{h^{i_0+1}_{n_t}\}$.

\item"{(3)}" after rescaling, for all the elements $u \in B_2$ (w.r.t.\ the fixed generating set $s_1,\ldots,s_r$ of the semigroup $S$ that generate $L$),
the lengths of the segments $[y_0,\psi^{i_0+1} \circ \psi^{i_0}(u)(y_0)]$ in $Y_{i_0+2}$, differ by at most
$\epsilon_{i_0+2}$ from the lengths of the corresponding segments, $[y_0,\psi^{i_0}(u)(y_0)]$ in the limit tree $Y_{i_0+1}$.

\item"{(4)}" the segments $[y_0,s_j(y_0)]$ are degenerate or positively oriented in $Y_{i_0+2}$, and for each $t$,
$h^{i_0+1}_{n_t} \, \circ \, \nu^{i_0+2}_{n_t}(s_j) \in FS_k$ and
$h^{i_0+1}_{n_t} \, \circ \, \nu^{i_0+2}_{n_t} \, \circ \, \psi^{i_0+1} \, \circ \, \psi^{i_0}(s_j) \in FS_k$.
\endroster

We continue iteratively, constructing at each level $i$  the automorphisms $\{\nu^{i+1}_{n_t}\}$, and
$\psi^i$. Note that in the iterative construction the automorphisms $\psi^i$ are fixed at step $i$, and the
sequence $\{\nu^{i+1}_{n_t}\}$ is constructed at step $i$. All these automorphisms are from the modular group
$Mod(\Theta_{i_0})$.

We set the homomorphisms $f_i$, the sequence of homomorphisms that appear in part (2) of the statement of the
proposition, to be:
$$f_i\ = \ h^i_{n_i} \, \circ  \, \nu^{i+1} \, \circ \psi^i \, \circ \, \ldots \, \circ \, \psi^{i_0}$$
where $\nu^{i+1}=\nu^{i+1}_{n_i}$, and for a suitable strictly increasing sequence of indices: $\{n_i\}$, such that the sequence $\{f_i\}$ converges into
an action of the limit group $L$ on a real tree $Y$. By construction, every non-trivial element of $L$ is mapped to a
non-trivial element of $F_k$ for large enough $i$, since the homomorphisms: $h^i_{n_t} \circ \nu^{i+1}_{n_t}$
converge into a faithful action of $L$ on a real tree. Hence, the action of $L$ on the limit tree $Y$
is faithful. Furthermore, by the construction of the automorphisms
$\psi^i$ and $\nu^{i+1}_{n_t}$, every element $g \in L$ that is not elliptic in the abelian decomposition
$\Theta_{i_0}$ acts hyperbolically on the real tree $Y$. Clearly, every element that is elliptic in $\Theta_{i_0}$ fixes
a point in $Y$, since it is elliptic in all the abelian decompositions $\Lambda_i$, $i \geq i_0$, and the autmorphisms
$\psi^i$ map it to a conjugate. Therefore, the abelian decomposition that is associated with the (faithful) action of
$L$ on the limit tree $Y$ is the stable dominant abelian decomposition $\Theta_{i_0}$.

\hfill$\square$

Proposition 6.3 enables one to replace a suffix of the sequence of abelian decompositions: $\Lambda_1,\ldots$ that was constructed
by iteratively shortening the given sequence of homomorphisms $\{h_n\}$ (and further pass to subsequences), with a single abelian
decomposition that minimally dominates the suffix.

Proposition 6.3 assumes that all the constructed abelian decompositions
$\{\Lambda_i\}$ do not contain QH vertex groups. To prove theorem 6.1 in case the limit group $L$ contains no non-cyclic abelian
vertex groups, we need to generalize the construction of the abelian decompositions $\{\Lambda_i\}$, and the conclusion
of proposition 6.3, in case the abelian decompositions $\Lambda_i$ do contain QH vertex groups.

\medskip
Suppose that $L$ is freely indecomposable with no non-cyclic abelian subgroups. Recall that  $\Lambda$ is the abelian decomposition
that $L$ inherits from its action on the real tree $Y$, that it is obtained as a limit from the convergent sequence of homomorphisms
$\{h_n\}$. We further refine $\Lambda$. If a non-QH vertex group is connected only to periodic edge groups,
we  restrict the homomorphisms $\{h_n\}$ to such a vertex group and obtain a non-trivial splitting of
it in which all the previous (periodic) edge groups are elliptic. Hence, the obtained abelian decomposition
of the non-QH vertex group can be used to refine the abelian decomposition $\Lambda$. Repeating this refinement
procedure iteratively, we
get an abelian decomposition that we denote $\Lambda_1$.

We fix finite generating sets of all the edge groups and all the non-QH vertex groups in $\Lambda_1$.  We divide the edges and
non-QH vertex
groups in $\Lambda_1$ into
finitely many equivalence classes of their growth rates as we did
in case $\Lambda_1$ contained no QH vertex groups. Two non-QH vertex groups or edge groups (or an edge and a non-QH vertex group)
 are said to be in the same equivalence
class if the maximal length of the image of their finite set of generators have comparable lengths, i.e., the length of a maximal
(translation) length image of one is bounded by a (global) constant times the maximal length of an image of the other, and vice versa. After passing
to a subsequence of the homomorphisms $\{h_n\}$ (that we still denote $\{h_n\}$), the non-QH and edge groups in $\Lambda_1$ are divided
into finitely many equivalence classes. We say that one class dominates another, if the maximal length of an image of a generator
of a group from the first class dominates the maximal length of an image of the second, but not vice versa. By definition, the
classes are linearly ordered (possibly after passing to a further subsequence). There exists a class that dominates
all the other classes, that
we call the $dominant$ class, that includes (possibly) dominant edge groups and (possibly) dominant non-QH vertex groups.

As in the case in which $\Lambda_1$ contains no QH vertex groups, we denote by $Mod(\Lambda_1)$ the modular group that is associated
with $\Lambda_1$. We set $MXMod(\Lambda_1)$ to be the $dominant$ subgroup that is generated by Dehn twists along
dominant edge groups,
and modular groups of those QH vertex groups that contain at least one dominant boundary element (i.e. a QH vertex group that is connected to at least one
dominant edge group in $\Lambda_1$).
We call these QH vertex groups,
$dominant$ QH vertex groups.

As we did in the simplicial case, we start by using the full modular group $Mod(\Lambda_1)$.
For each index $n$, we set the
pair homomorphism: $h_n^1:(S,L) \to (FS_k,F_k)$,
$h_n^1=h_n \circ \varphi_n$, where $\varphi_n \in Mod(\Lambda_1)$, to be a shortest pair homomorphism
that is obtained from $h_n$ by precomposing it with a modular automorphism from $Mod(\Lambda_1)$.
If
there exists a subsequence of the homomorphisms $\{h_n^1\}$ that converges into a proper quotient of the pair
$(S,L)$, we set the limit of this subsequence to be $(S_f,L_f)$, and the conclusion of the theorem follows
with a resolution of length 1.

Therefore, we may assume that every convergent subsequence of the homomorphisms $\{h_n^1\}$
converges into a faithful action of the limit group $L$ on some real tree. In that case we use only elements from the dominant
modular group, $MXMod(\Lambda_1)$.

First, we shorten the action of each of the dominant QH vertex groups
using the
procedure that is used in the proof of theorem 2.2.
%
Since an IET action of a QH vertex
group $Q$ on a real tree is indecomposable in the sense of $[Gu]$, if we fix a finite set of (non-elliptic) generators of the QH vertex group $Q$,
that we assume to be a
finite set of positively oriented generators, then
finitely many (fixed) translates of  the  path in the tree that is associated with one of the fixed set of generators,
cover the  path that is associated with any of the other fixed generators.
Keeping this covering along the shortening procedure
guarantees that the ratios between the lengths of the paths that are associated with the different generators
remain globally bounded along the shortening procedure.




At this point we shorten the homomorphisms $\{h_n\}$ using the dominant modular group $MXMod(\Lambda_1)$. For each homomorphism
$h_n$, we pick the shortest homomorphism after precomposing with an element from $MXMod(\Lambda_1)$ that keeps the positivity
of the given set of the images of the generators $s_1,\ldots,s_r$, and the covering of paths that are associated with fixed sets of generators of QH vertex groups
by translates of paths that are associated with other paths of any of the other generators.
We further keep the length of the fixed set of generators of a dominant QH vertex group
to be at least the maximal length of the image of a generator of a dominant
edge or non-QH vertex group (the generators are chosen from the fixed finite sets of generators of each of the vertex and edge groups).
We further require that if a path that is associated with one of the generators,
$s_1,\ldots,s_r$, passes through an edge with a dominant edge group in $\Lambda_1$, then after shortening the path that
is associated with such a generator contains  a subpath that is associated with the dominant edge group. This guarantees
that the length of such a generator remains bigger than the length of the dominant edge group along the entire procedure.
We (still) denote the obtained (shortened) homomorphisms $\{h_n^1\}$.

By the shortening arguments that are used to prove theorem 2.2, the lengths of the images, under the shortened
homomorphisms $\{h^1_n\}$, of the fixed set of generators
$s_1,\ldots,s_r$, and the fixed set of generators of the dominant QH vertex groups,
are bounded by some  constant $c_1$ (that is independent of $n$) times the maximal length of the images of
the fixed generators of the dominant
vertex and edge groups.

We pass to a subsequence of the homomorphisms $\{h_n^1\}$ that converges into an action of $L$ on some
real tree with an associated abelian decomposition $\Delta_2$. If the action of $L$ is not faithful, the conclusions of
theorem 6.1 follow. Hence, we may assume that the action of $L$ is faithful.
We further refine $\Delta_2$, by restricting (a convergent subsequence of)
the homomorphisms to non-QH vertex groups that are connected only
to periodic edge groups,  as we did with $\Lambda$, and
construct an abelian decomposition that we denote
$\Lambda_2$.


By construction, a dominant edge group, a dominant boundary component of a QH vertex group,
  a dominant non-QH vertex group, and a dominant QH vertex group
in $\Lambda_1$
can not be elliptic in $\Lambda_2$ (i.e., they can not be contained in a non-QH vertex group or an edge group in
$\Lambda_2$).
Furthermore, a  QH vertex group that has a dominant boundary component in $\Lambda_1$, is not elliptic nor a QH vertex group
in $\Lambda_2$.

With $\Lambda_2$ we associate its modular group $Mod(\Lambda_2)$. We further associate with
$\Lambda_2$ its dominant edge groups and dominant QH and non-QH vertex groups, and the modular group that is associated
only with dominant edge groups, and with dominant QH vertex groups,
that we denote $MXMod(\Lambda_2)$.

We continue in a similar way to what we did with the sequence of homomorphisms $\{h_n\}$ and with $\Lambda_1$.
We first shorten using automorphisms from the ambient modular group $Mod(\Lambda_2)$. If there exists a subsequence of
shortened homomorphisms that converges into a proper quotient of $L$, the conclusion of theorem 6.1 follows. If
there is no such subsequence, we restrict the shortenings to automorphisms
from the dominant modular group $MXMod(\Lambda_2)$, and shorten the homomorphisms that were used to construct the abelian decomposition
$\Lambda_2$ in the same way that we shortened the homomorphisms  $\{h_n\}$ (from which we obtained the abelian decomposition
$\Lambda_1$).

We associate sets of  generators with each of the new dominant QH vertex groups in $\Lambda_2$, i.e., those QH vertex
groups in $\Lambda_2$ that are not QH vertex groups in $\Lambda_1$. We further choose finitely many translates of the
paths that are associated with each generator that cover the paths that are associated with  the other generators.
These translates guarantee that the ratios between the lengths of the paths that are associated with these generators
along the iterative (shortening) procedure that we use
remain globally bounded.

As we did in the first step, we shorten the homomorphisms $\{h^1_n\}$ using the dominant modular group $MXMod(\Lambda_2)$.
For each homomorphism
$h^1_n$, we pick the shortest homomorphism after precomposing it with an element from $MXMod(\Lambda_2)$ that keeps the positivity
of the given set of the images of the generators $s_1,\ldots,s_r$,
and positivity of the generators of the dominant edge groups in $\Lambda_1$ and $\Lambda_2$. We
further keep the covering of paths that are associated with the fixed set of generators of QH vertex groups in $\Lambda_1$ and $\Lambda_2$ by
finite number of translates the paths
that are associated with the other generators of these QH vertex groups. We also keep  the positivity and negativity of the positively and negatively oriented
subpaths in the paths that are associated
with elements in the ball of radius 2 in the Cayley graph of $L$ w.r.t. $s_1,\ldots,s_r$.
We require the shortened homomorphisms to  keep the lengths of the fixed set of generators of the QH vertex groups in $\Lambda_2$
to be at least the maximal length of the image of a generator of a dominant
edge or vertex group in $\Lambda_2$
(the generators are chosen from the fixed finite sets of generators of each of the vertex and edge groups).

We further require that after the shortening,
if a path in the limit tree $Y_2$ that is associated with $\Lambda_2$ and with one of the following:
\roster
\item"{(1)}" a generator
of a QH vertex group in $\Lambda_1$ or a  generator of the dominant edge groups in $\Lambda_1$.

\item"{(2)}" a positively or
negatively oriented subpath
in the paths that are associated
with elements in the ball of radius 2 in the Cayley graph of $L$ w.r.t. $s_1,\ldots,s_r$.
\endroster
passes through an edge with a dominant edge group in $\Lambda_2$, then after shortening the path that
is associated with such a generator contains at least a subpath that is associated with the dominant edge group. This guarantees
that the length of the corresponding element  remains bigger than the length of the dominant edge group along the entire procedure.
We  denote the obtained (shortened) homomorphisms $\{h_n^2\}$.


By the shortening arguments that are proved along the proof of theorem 2.2,  lengths of the images, under the shortened
homomorphisms $\{h^2_n\}$, of the fixed set of generators of the dominant QH vertex groups
are bounded by some  constant $c_2$ (that is independent of $n$) times the maximal length of the images of
the fixed generators of the dominant
vertex and edge groups.

We pass to a subsequence of the homomorphisms $\{h_n^2\}$ that converges into an action of $L$ on some
real tree with an associated abelian decomposition $\Delta_3$. If the action of $L$ is not faithful, the conclusions
of theorem 6.1 follow. Hence, we may assume that the action of $L$ is faithful.
We refine $\Delta_3$ to an abelian decomposition $\Lambda_3$,
precisely as we refined $\Delta_2$ to obtain $\Lambda_2$.

With $\Lambda_3$ we associate its modular group $Mod(\Lambda_3)$ and dominant modular group,
$MXMod(\Lambda_3)$. We first shorten using $Mod(\Lambda_3)$, and if every convergent shortened subsequence converges
into a faithful action that we assume to be a
finite set of positively oriented generators,  of $L$, we further use the dominant modular group $MXMod(\Lambda_3)$.

In shortening using $MXMod(\Lambda_3)$, we keep the positivity of all the fixed (positive) sets of generators of the
QH vertex groups in $\Lambda_1$,  $\Lambda_2$ and $\Lambda_3$, and the generators of the dominant edge groups in
$\Lambda_1$ and $\Lambda_2$. We also keep the positivity and negativity of all the finitely many positive and negative
subpaths of the paths that are associated with the elements in the ball of radius 3 in the Cayley graph of
$L$ w.r.t. the generating set $s_1,\ldots,s_r$.


We denote the obtained (shortened) homomorphisms $h_n^3$, and continue iteratively.
If in all steps the obtained actions are faithful, we get an infinite sequence of
abelian decompositions, $\Lambda_1,\Lambda_2,\ldots$.
Given the infinite sequence of abelian decompositions, we define the stable dominant abelian decomposition,
$\Theta_{i_0}$, in the same way as we did in the simplicial case (definition 6.2). From the convergent
sequences $\{h^i_n\}_{n=1}^{\infty}$, we choose a subsequence $\{f_i\}$ that has a subsequence that converges into a
faithful action of the limit groups $L$, with an associated abelian decomposition $\Theta_{i_0}$.

\vglue 1pc
\proclaim{Proposition 6.4} Suppose that the limit group $L$ is freely indecomposable, and contains no non-cyclic abelian subgroups.
It is possible to choose a sequence of pair homomorphisms, $f_i:(S,L) \to (FS_k,F_k)$, such that:
\roster
\item"{(1)}" for each $i$, $f_i$ is a homomorphism from the sequence $\{h^i_n\}$.

\item"{(2)}" the sequence $f_i$ has a subsequence that converges into a faithful action of $L$ on a real tree $Y$.

\item"{(3)}" With the action of $L$ on the limit tree $Y$ there is an associated abelian decomposition $\Delta$. By possibly
refining $\Delta$, using restrictions of the homomorphisms $\{f_i\}$ to some of its non-QH vertex groups, it is possible to
obtain the stable abelian decomposition, $\Theta_{i_0}$.
\endroster
\endproclaim

\nfp For each index $i$, we choose $f_i$ to be a homomorphism $h^i_n$ (from the sequence $\{h^i_n\}$), that maps the fixed positive
generators of all the QH vertex groups and the generators of all
the edge groups in $\Lambda_1,\ldots,\Lambda_{i}$ to $FS_k$. We further require that the ratios between the lengths of the images
of these generators, and the ratios between the lengths of the images of all the elements in the ball of radius $i$
in the Cayley graph of $L$ w.r.t. the generating set $s_1,\ldots,s_r$ will be approximately the ratios between the
lengths of the paths that are associated with these elements in the corresponding limit tree $Y_{i+1}$ (the trees that
are obtained as the limits of the sequences $\{h^i_n\}$). We further require that $f_i$ maps the elements in a ball of radius $i$
in the Cayley graph of $L$ monomorphically into $F_k$.

The sequence $\{f_i\}$ has a subsequence that converges into a faithful action of $L$ on a real tree $Y$. Let $\Delta$
be the abelian decomposition that is associated with the action of $L$ on $Y$. Since every elliptic element in the stable abelian
decomposition $\Theta_{i_0}$ must fix a point in $Y$, $\Delta$ is dominated by $\Theta_{i_0}$.

\vglue 1pc
\proclaim{Lemma 6.5} Let $Q$ be a QH vertex group in $\Theta_{i_0}$ that does not appear in any of the abelian
decompositions $\Lambda_i$, for $i \geq i_0$.

If there is a non-peripheral element in  $Q$ that fixes a point in
$Y$, then the entire QH vertex group $Q$ fixes a point in $Y$.
\endproclaim

\nfp
Suppose that a non-peripheral element $q \in Q$ fixes a point in $Y$. $q$ is contained in some ball $B_m$ in the Cayley graph
of $L$ w.r.t. $s_1,\ldots,s_r$. $\Theta_{i_0}$ is the stable dominant abelian decomposition of the sequence of abelian
decompositions: $\Lambda_1,\ldots$. Hence, there must exist an abelian decomposition $\Lambda_i$, for some $i>max(i_0,m)$, for which
either:
\roster
\item"{(i)}" $q$ is a non-peripheral element in some QH vertex group in $\Lambda_i$.

\item"{(ii)}" $q$ is hyperbolic in the abelian decomposition $\Lambda_i$.
\endroster
According to the procedure that was used to construct the abelian decompositions, $\Lambda_1,\ldots$, if either
(i) or (ii) hold for $\Lambda_i$, then the traces and the lengths of $q$ in its actions on the trees that are
associated with the abelian decompositions $\Lambda_{i+1},\ldots$ (i.e., the limit trees $Y_{i+1},\ldots$),
are bounded below by either a (global) positive constant times the lengths of the fixed set of generators of
the QH vertex group  that contains $q$ in $\Lambda_i$ (in case (i)), or by either  a (global) positive constant times
the lengths of a fixed set of generators in a QH vertex group in $\Lambda_i$, or the length of a generator of an edge group
in $\Lambda_i$ that becomes a dominant edge group in some $\Lambda_{i'}$ for some $i' \geq i$ in case (ii).

By the properties of the stable abelian decomposition $\Theta_{i_0}$, and the structure of the procedure for the construction
of the abelian decompositions $\Lambda_1,\ldots$, for some index $j_0 > i$,  the lengths of either:
\roster
\item"{(1)}"  the fixed set of generators of a QH vertex group in $\Lambda_i$ that contains the element $q$

\item"{(2)}" the fixed set of generators of a QH vertex group $\tilde Q$ in $\Lambda_i$, or a generator of an abelian group $\tilde E$ in $\Lambda_i$,
such that $q$ is hyperbolic w.r.t.\ the splitting that is obtained from $\Lambda_i$ by collapsing all the edges except the abelian edge group $\tilde E$ in $\Lambda_i$
or collapsing all those edge groups except those that
are connected to the QH vertex group $\tilde Q$ in $\Lambda_i$.
\endroster
have to be bigger than fixed multiples of the lengths of  the fixed
sets of generators of all the QH vertex groups and the lengths of generators of all the edge groups in the abelian decompositions
$\Lambda_j$, $j \geq j_0$,  for
edge groups and QH vertex groups that are
contained in the QH vertex group $Q$ in $\Theta_{i_0}$ that contains the element $q$.

Hence, if the element $q$ is elliptic in the limit action on the real tree $Y$, all the QH vertex groups and all the edge groups that
appear in $\Lambda_j$ for $j \geq j_0$, and are contained in $Q$, must fix points in $Y$, which means that
the entire QH vertex group $Q$ in $\Theta_{i_0}$
that contains the non-peripheral element $q$ must fix a point in $Y$.

\hfill$\square$

Edge groups and QH vertex groups in $\Theta_{i_0}$, that appear in some abelian decomposition $\Lambda_j$,
for some $j \geq i_0$, have to  fix points in the limit tree $Y$. From the graph of groups $\Theta_{i_0}$ we take out all
the edge groups and all the QH vertex groups that appear in some abelian decomposition $\Lambda_j$ for some $j \geq i_0$.
We denote the obtained (possibly disconnected) graph of groups $\hat \Theta$. $\hat \Theta$ contains QH vertex groups and non-QH vertex
groups that are connected to them.

We restrict the homomorphisms $\{f_i\}$ to the fundamental groups of each of the connected components in $\hat \Theta$. By lemma
6.5 if a component of $\hat \Theta$ contains a QH vertex group, then the restrictions of the sequence $\{f_i\}$ to the
fundamental group of that component subconverges into an action on a real tree, such that the abelian decomposition
that is associated with that action, contains one or more QH vertex groups from $\hat \Theta$, that are also QH vertex
groups in $\Theta_{i_0}$.

We continue the refinement process by erasing these QH vertex groups from $\hat \Theta$, and restrict the homomorphisms $\{f_i\}$
to the fundamental groups of connected components in the remaining graph of groups. By lemma 6.5, after finitely many such
revisions of the graph of groups $\hat \Theta$ (i.e. erasing the QH vertex groups that were visible), we uncover all the QH
vertex groups in $\Theta_{i_0}$ that do not appear as QH vertex groups in an abelian decomposition $\Lambda_j$ for $j \geq i_0$.

Now, starting with the abelian decomposition $\Delta$, that was read from the original faithful action of the limit group $L$ on
the limit tree $Y$, possibly refining $\Delta$ using the QH vertex groups and the edge groups from $\Theta_{i_0}$ that appear
in some abelian decomposition $\Lambda_j$ for some $j \geq i_0$, and further refining the obtained decomposition using the uncovered
QH vertex groups, we finally obtain the stable dominant abelian decomposition $\Theta_{i_0}$. This concludes the proof of
proposition 6.4.

\hfill$\square$

Let $j_0 \geq i_0$ be an index for which all the QH vertex groups  in $\Theta_{i_0}$ that appear in some
$\Lambda_j$ for some $j \geq i_0$, already appear in $\Lambda_j$ for some $i_0 \leq j \leq j_0$.
Proposition 6.4 enables us to replace the suffix of the sequence of abelian decompositions, $\Lambda_{i_0},\ldots$, with a finite
resolution $\Lambda_{i_0},\ldots,\Lambda_{j_0},\Theta_{i_0}$, and hence the entire sequence $\Lambda_1,\ldots$
with the finite resolution:   $\Lambda_1,\ldots,\Lambda_{j_0}, \Theta_{i_0}$. Furthermore, we may continue with
the sequence of pair homomorphisms $\{f_i\}$, that are obtained from a subsequence of the pair homomorphisms that we
started with $\{h_n\}$, by precomposing them with automorphisms from the modular groups that are associated with the abelian
decompositions that appear along the finite resolution.

We continue to the next step starting with the sequence of homomorphisms $\{f_i\}$, and the abelian decomposition $\Theta_{i_0}$.
In the abelian decomposition that is obtained from a subsequence of homomorphisms that were obtained from the
homomorphisms $\{f_i\}$ by precomposing them with automorphisms from the (dominant) modular group of $\Theta_{i_0}$,
either a dominant edge group or a dominant non-QH vertex group is not elliptic.

At this point we repeat the whole construction of a sequence of abelian decompositions.
If the sequence terminates after a finite number of steps
 the conclusion of theorem 6.1 follows. Suppose it ends up with an infinite sequence of abelian decompositions.

Note that for any abelian decomposition $\Lambda_i$ along this sequence, after finitely many steps there exists an
abelian decomposition $\Lambda_{i'}$ for some $i' \geq i$, with a dominant edge group. This means that the dominant edge
group in $\Lambda_{i'}$ is not elliptic in $\Lambda_{i'+1}$. Furthermore, as long as the boundary elements of a QH
vertex group $Q$ in $\Lambda_i$  remain elliptic in the next abelian decompositions, the QH vertex group $Q$ remains a QH vertex
group in the next abelian decompositions.

Using
proposition 6.4, we can replace the suffix of the sequence with a finite resolution that stably
dominates the entire sequence  (see definition 6.2).

We continue iteratively. At each step we start with a sequence of homomorphisms that are obtained from a subsequence of
the  homomorphisms $\{h_n\}$ by precomposing them with automorphisms from the modular groups
of the abelian decompositions that appear along the previously constructed finite resolution. We continue for a single step,
so that at least one dominant edge group or a dominant vertex group in the previous abelian decomposition is not elliptic.
Then we either associate with the shortened sequence of homomorphisms
a finite resolution, that completes the proof of theorem 6.1, or we associate with it an infinite sequence of
abelian decompositions. By proposition 6.4 a suffix of this last infinite sequence can be replaced by a finite resolution that ends with
a stable dominant abelian
decomposition of it.

If this iterative procedure terminates after finitely many steps, the conclusion of theorem 6.1 follows. Otherwise we obtained
an infinite sequence of abelian decompositions. Note that the infinite sequence of abelian decompositions contains a subsequence
of abelian decompositions that contain  QH vertex groups, and dominant edge groups in the abelian decompositions from the
subsequence are not elliptic in the abelian decomposition that appears afterwards in the sequence.

Now, we apply proposition 6.4 to the sequence of abelian decompositions that we constructed. By proposition 6.4 a suffix of the
sequence can be replaced with a finite resolution that terminates with the abelian decomposition that stably dominates the original suffix
of the sequence. Since the sequence of abelian
decompositions contains a subsequence with QH vertex groups, the stable dominant abelian decomposition must contain a
QH vertex group as well. Since every abelian decomposition in this subsequence contains a (dominant) edge that is not elliptic
in the next abelian decomposition, the stable dominant abelian decomposition contains either 2 QH vertex groups, or a
single QH vertex group with an associated surface group $S$ that is either:
\roster
\item"{(i)}" an orientable surface  with $\chi(S) \leq -2$.

\item"{(ii)}" a non-orientable surface with $genus(S)+bnd(S) \geq 3$, where $bnd(S)$ is the number of boundary
components of the surface $S$.
\endroster

We repeat the whole construction starting with the abelian decomposition that we obtained and the subsequence of homomorphisms
that is associated with it according to the construction that is used in the proof of proposition 6.4. Either the construction
terminates in finitely many steps, or we get a sequence of abelian decompositions that has a subsequence that satisfies the
properties that the previous stably dominant abelian decomposition satisfied, and in each abelian decomposition from the
subsequence there exists a dominant edge group that is not elliptic in the next abelian decomposition.

Once again we apply proposition 6.4 to the constructed sequence of abelian decompositions. The stable dominant abelian
decomposition of the sequence has to contain QH vertex groups with topological complexity that is bounded below by
a larger lower bound. Hence, it either contains at least
3 QH vertex groups, or 2 QH vertex groups so that at least one of them
satisfies properties (i) or (ii), or a single QH vertex group with an associated surface group $S$ that is either orientable
with $\chi(S) \leq -3$ or non-orientable with $genus(S)+bnd(S) \geq 4$.

We continue iteratively. Since after each iteration (of the entire construction) we obtain an abelian decomposition
that contain QH vertex groups with topological complexity that is bounded below by larger and larger bounds, and since
the obtained abelian decompositions are all dominated by the JSJ decomposition of the freely-indecomposable limit
group $L$, the procedure has to terminate after finitely many iterations, and the obtained resolution satisfies the
conclusions of theorem 6.1.

\hfill$\square$

So far we assumed that $L$ is freely indecomposable and contains no non-cyclic abelian groups. To get a conclusion similar
to the one that appears in theorem 6.1 in the presence of non-cyclic abelian subgroups, we need to slightly modify it. Instead of
using only the modular group or the dominant modular group, we need to further
allow generalized Dehn twists in the presence of a non-cyclic abelian group.
i.e., we further allow Dehn twists in roots of the the values of a generator of an edge group
(in case the edge group belongs to a non-cyclic
maximal abelian subgroup), and not just modular automorphisms.

\vglue 1pc
\proclaim{Definition 6.6} Let $(S,L)$ be a pair in which $L$ is a freely indecomposable limit group, and let $\{h_n:(S,L) \to (FS_k,F_k)\}$
be a sequence of homomorphisms that converge into a faithful action of $L$ on a real tree with an associated abelian decomposition $\Lambda$.

With $\Lambda$ we associate its dominant rigid and QH vertex groups, and dominant edge groups. We say that an abelian vertex group $A$ in $\Lambda$ is
dominant, if $A=A_0+<a_1,\ldots,a_{\ell}>$, where $A_0<A$ is the subgroup that is generated by the edge groups that are connected to the
vertex that is stabilized by $A$, and the lengths of the images $\{h_n(a_i)\}$, $1 \leq i \leq \ell$, grows more than a constant times the
lengths of the images of the generators of dominant edge groups in $\Lambda$,  for at least one of the indices $i$, $1 \leq i \leq \ell$.

Let $h:(S,L) \to (FS_k,F_k)$ and $f:(S,L)\to (FS_k,F_k)$ be two pair homomorphisms. We
say that $f$ is obtained from $h$ using $generalized$ $Dehn$ $twists$ if:
\roster
\item"{(i)}" there exists a pair homomorphism $\hat h:(S,L) \to (FS_k,F_k)$ that is obtained from $h$ by precomposing it
with a modular automorphism of $\Lambda$: $\hat h= h \circ \varphi$, $\varphi \in Mod(\Lambda)$.

\item"{(ii)}" suppose that  $A<L$ is an abelian vertex group in $\Lambda$. Suppose that:
$A=A_0+<a_1,\ldots,a_{\ell}>$ as above. Then $f$ is obtained
from $\hat h$ by modifying the images of the (non-elliptic) generators $a_1,\ldots,a_{\ell}$.  For each $i$, $1 \leq i \leq \ell$,
$f(a_i)$ is an element from the cyclic group $\hat h(A)$ in $F_k$.
\endroster

We say that a pair homomorphism $f$ is obtained from a homomorphism $h$ using $dominant$ generalized Dehn twists if:
\roster
\item"{(iii)}" there exists a pair homomorphism $\hat h:(S,L) \to (FS_k,F_k)$ that is obtained from $h$ by precomposing it
with a dominant modular automorphism of $\Lambda$: $\hat h= h \circ \varphi$, $\varphi \in MXMod(\Lambda)$.

\item"{(iv)}"
$f$ is obtained
from $\hat h$ by modifying the images of only dominant abelian vertex groups in $\Lambda$ (where the modification type is specified in part (ii)).
\endroster

Note that in general $f$ is not obtained from $h$ by a precomposition with a modular automorphism.
Also,
note that generalized Dehn twists are used in studying systems of equations with parameters over a free group (see sections
9-10 in [Se1]).
\endproclaim

The addition of generalized Dehn twists
enable us to use our treatment of non-cyclic abelian subgroups in sections 2 (theorem 2.1) and 3 and generalize the statement
of theorem 6.1 to include all pairs $(S,L)$ with $L$ a general freely indecomposable limit group.

\vglue 1pc
\proclaim{Theorem 6.7}
Let $(S,L)$ be a pair, where $L$ is a freely indecomposable limit group, and let $s_1,\ldots,s_r$ be a fixed generating
set of the semigroup $S$.
Let $\{h_n:(S,L) \to (FS_k,F_k)\}$ be a sequence of pair homomorphisms that converges into
a faithful action of $L$ on a real tree $Y$.

Then there exists a $resolution$:
$$(S_1,L_1) \to  (S_2,L_2) \to \ldots \to (S_m,L_m) \to (S_f,L_f)$$
that satisfies the following properties:
\roster
\item"{(1)}" $(S_1,L_1)=(S,L)$, and  $\eta_i:(S_i,L_i) \to (S_{i+1},L_{i+1})$ is an isomorphism for $i=1,\ldots,m-1$ and
$\eta_m:(S_m,L_m) \to (S_f,L_f)$ is a
proper quotient map.

\item"{(2)}" with each of the pairs $(S_i,L_i)$, $1 \leq i \leq m$, there is an associated abelian decomposition that we denote
$\Lambda_i$.

\item"{(3)}" there exists a subsequence of the homomorphisms $\{h_n\}$ that factors through the resolution. i.e., each homomorphism
$h_{n_r}$ from the subsequence, is obtained from a homomorphism of the terminal pair $(S_f,L_f)$ using a composition of
a
modification that uses generalized Dehn twists that are associated with $\Lambda_m$, and modular automorphisms that
are associated with $\Lambda_1,\ldots,\Lambda_{m-1}$.
\endroster
\endproclaim

\nfp Suppose that $L$ is freely indecomposable and   contains a
non-cyclic abelian group. If all the non-cyclic abelian groups in $L$ remain elliptic in all the abelian decompositions that
are constructed along the iterative procedure that was
used in the proof of theorem 6.1, in case there are no  abelian vertex groups, i.e., in the iterative applications of the
construction that is used in the proof of proposition 6.4, the same construction proves the conclusion of theorem 6.7.

Suppose that  at some
step along an application of the iterative procedure that is used in the proof of proposition 6.4,
$\Lambda_i$ contains an abelian vertex group.

In that case  $A$ is either the set
stabilizer of an axial component or it is the set stabilizer of a line in the simplicial part of $Y_i$. Let $A_0<A$ be the point
stabilizer of the axial component with set stabilizer $A$ or of the axis of $A$.

$A_0$ is the stabilizer of an edge  in $\Lambda_i$. If $A_0$ is not dominant, we don't include Dehn twists along elements
of $A$ in the dominant modular group of $\Lambda_i$, and proceed with the procedure that is
used in the proof of theorem 6.1 as long as $A_0$ is not dominant. The abelian decompositions that we
consider in all the steps of the procedure in which $A_0$ is not dominant, are abelian decompositions relative to $A_0$.

If $A_0$ remains not dominant along the entire procedure, the procedure that is used in proving proposition 6.4 and its
conclusions remain valid. Suppose that at some step $j$, $A$ is dominant. $A=A_0+<a_1,\ldots,a_{\ell}>$. By theorem
2.1 in the axial case, and in case $A$ acts simplicially, when $A$ is  dominant we can modify a subsequence of the
homomorphisms $\{h_n\}$ using generalized Dehn twists along $A$ (that do not change $A_0$), so that the images
of $a_1,\ldots,a_{\ell}$ under the modified homomorphisms are identified with fixed elements in $A_0$. Hence, the modified sequence
of homomorphisms converges into a limit group in which the image of $A$ is $A_0$. Therefore, the modified sequence of
homomorphisms converges into a proper quotient of the limit group $L$, and the conclusion of theorem 6.6 follows.

\hfill$\square$

Theorem 6.7 proves that given a pair $(S,L)$ in which $L$ is freely indecomposable, and a
sequence of pair homomorphisms $\{h_n:(S,L) \to (FS_k,F_k)\}$, that converges into a faithful action of $L$ on some real tree,
it is possible to extract a subsequence that factors through a finite resolution of the pair $(S,L)$ that terminates in a
proper quotient of the pair $(S,L)$. Using a compactness argument it is not difficult to apply the conclusion
of theorem 6.7 and get a Makanin-Razborov diagram for a pair $(S,L)$, when $L$ is a restricted limit group with no
freely decomposable restricted limit quotients.

\medskip
Let $FS_k=<a_1,\ldots,a_k>$ be a free semigroup that generates the free group $F_k$,
and let $(S,L)$ be a restricted pair, i.e.,
a pair that contains the subpair $(FS_k,F_k)$. Suppose that the pair $(S,L)$ does not admit
a quotient restricted map $\eta: (S,L) \to (\hat S, \hat L)$ with the following properties:
\roster
\item"{(1)}" $\eta$ maps the subpair $(FS_k,F_k)<(S,L)$ monomorphically onto the corresponding subpair
$(FS_k,F_k)<(\hat S, \hat L)$.

\item"{(2)}" $\hat L$ admits a non-trivial free decomposition in which $\eta(F_k)$ is
contained in a factor.
\endroster

The restricted pair $(S,L)$ is in particular freely indecomposable with respect to the subpair $(FS_k,F_k)$,
and so is every restricted quotient pair of $(S,L)$. Given a sequence of homomorphisms $\{h_n:(S,L) \to (FS_k,F_k)\}$, it
is possible to extract a subsequence that converges into a faithful action of some restricted quotient $(S_1,L_1)$
of $(S,L)$. By theorem 6.7 it is possible to further extract a subsequence that factors through a resolution
that terminates in a proper quotient of $(S_1,L_1)$. Applying theorem 6.7 iteratively, finitely many times, it
is possible to extract a further subsequence that factors through a resolution:
$$(S_1,L_1) \to  (S_2,L_2) \to \ldots \to (S_m,L_m) \to (FS_k,F_k)$$
that terminates in the standard pair $(FS_k,F_k)$, and
in which some of the epimorphisms: $\eta_i:(S_i,L_i) \to (S_{i+1},L_{i+1})$ are isomorphisms and some are proper quotient
maps.
With each of the pairs $(S_i,L_i)$, $1 \leq i \leq m$, there is an associated abelian decomposition that we denote
$\Lambda_i$, with which we naturally associate a modular group.

Although f.g.\ subsemigroups of limit groups need not be f.p.\ in general, a pair $(S,L)$ in which $L$ is a limit group
and $S$ is a f.g.\ subsemigroup of $L$ is naturally a finitely presented object. Hence, a resolution of the form:
$$(S_1,L_1) \to  (S_2,L_2) \to \ldots \to (S_m,L_m) \to (FS_k,F_k)$$
together with the abelian decompositions that are associated with the various restricted pairs $(S_i,L_i)$ can be encoded
using a finite amount of data. In particular, there are only countably many such resolutions (that are associated
with convergent restricted pair homomorphisms of a given pair $(S,L)$), and we can order these resolutions using their
finite encoding.

\vglue 1pc
\proclaim{Theorem 6.8} Let $(S,L)$ be a restricted pair, that has no restricted quotient $(\hat S, \hat L)$ in which $\hat L$
admits a free decomposition relative to (the embedding of the pair) $(FS_k,F_k)$. Then there are finitely many resolutions of the form that
is constructed in theorem 6.7:
$$(S_1,L_1) \to  (S_2,L_2) \to \ldots \to (S_m,L_m) \to (FS_k,F_k)$$
where $(S_1,L_1)$ is restricted quotient pair of $(S,L)$ , such that:
\roster
\item"{(1)}" every restricted pair homomorphism, $h:(S,L) \to (FS_k,F_k)$, factors through at least one of these
finitely many resolutions.

\item"{(2)}" for each of the resolutions in the collection, there exists a sequence of homomorphisms:
$\{h_n:(S,L) \to (FS_k,F_k)\}$, that converges into a faithful action of the initial pair $(S_1,L_1)$ on
a real tree with an associated abelian decomposition (after the refinement that is used in the proof of proposition
6.4) $\Lambda_1$.

Furthermore, the sequence of homomorphisms $\{h_n\}$ can be modified using modular automorphisms and generalized
Dehn twists that are associated with the abelian decompositions: $\Lambda_1,\ldots,\Lambda_m$, to get sequences
of pair homomorphisms $\{h_n^1\},\ldots,\{h_n^m\}$. Each of these modified sequences of homomorphisms
$\{h_n^i\}$ converges into a faithful action of the pair $(S_i,L_i)$ on a real tree.
\endroster
\endproclaim

\nfp Let $\{h \, | \, h:(S,L) \to (FS_k,F_k)\}$ be the collection of all the restricted pair homomorphisms of
the pair $(S,L)$. We look at all the possible subsequences of such homomorphisms $\{h_n\}$, that converge into
quotients of $(S,L)$, with which we can associate a resolution:
$$(S_1,L_1) \to  (S_2,L_2) \to \ldots \to (S_m,L_m) \to (FS_k,F_k)$$
by iteratively applying theorem 6.7. Note that each of these constructed resolutions satisfy part (2) of the theorem.

Each such resolution can be encoded using finite amount of data. In particular, there are countably many such resolutions,
and we can order them using the finite encoding. We argue that the collection of all the pair homomorphisms
of $(S,L)$ factor through
a finite collection of these resolutions.

Suppose that finitely many do not suffice. Then there exists a sequence of restricted pair homomorphisms: $\{h_n\}$, such
that for every index $n$, $h_n$ does not factor through the first $n$ resolutions from the ordered countable set
of resolutions that were constructed. By iteratively applying theorem 6.7, from the sequence $\{h_n\}$ we can extract a
subsequence $\{h_{n_r}\}$, that factors through a resolution of form we previously constructed:
$$(S_1,L_1) \to  (S_2,L_2) \to \ldots \to (S_m,L_m) \to (FS_k,F_k)$$
But this resolution appears in the countable set of resolutions that we associated with the pair $(S,L)$. Hence, it appears in the
ordered list. Therefore, for some index $r_0$, and for every $r>r_0$, the subsequence of homomorphisms: $\{h_{n_r}\}$
factors through a resolution that appears in the ordered list of resolutions, a contradiction to the choice of
the homomorphisms $\{h_n\}$.

\hfill$\square$

Theorem 6.8 constructs a diagram that encodes all the homomorphisms from a pair $(S,L)$ into the standard pair
$(FS_k,F_k)$. In addition, the construction of the resolutions in the diagram, that mainly uses the iterative procedure that was
used in proving theorem 6.1 and proposition 6.4, guarantees that there exist sequences of homomorphisms that
factor through them for which the corresponding shortened homomorphisms converge into the pairs that appear along the
resolutions. i.e., the construction of the resolutions  guarantees the existence of (weak) test sequences
for the resolutions in the diagram.

\vglue 1pc
\proclaim{Theorem 6.9} Let $(S,L)$ be a restricted pair, that has no restricted quotient $(\hat S, \hat L)$ in which $\hat L$
admits a free decomposition relative to (the embedding of) $(FS_k,F_k)$. Then there are finitely many resolutions of the form that
is constructed in theorem 6.7:
$$(S_1,L_1) \to  (S_2,L_2) \to \ldots \to (S_m,L_m) \to (FS_k,F_k)$$
where $(S_1,L_1)$ is a restricted quotient pair of $(S,L)$ , such that parts (1) and (2) in theorem 6.8 hold for these resolutions,
and in addition:
\roster
\item"{(3)}" for each of the resolutions in the collection, there exists a sequence of homomorphisms:
$\{h_n:(S,L) \to (FS_k,F_k)\}$, that
can be modified using modular automorphisms and generalized
Dehn twists that are associated with the abelian decompositions: $\Lambda_1,\ldots,\Lambda_m$, to get sequences
of pair homomorphisms $\{h_n=h_n^0\},\{h_n^2\},\ldots,\{h_n^m\}$. Each of these modified sequences of homomorphisms
$\{h_n^i\}$ converges into a faithful action of the pair $(S_i,L_i)$ on a real tree
 with an associated
abelian decomposition (after an appropriate refinement) $\Lambda_i$.
\endroster
\endproclaim

\nfp The construction of the resolutions in the diagram in theorem 6.8. starts with a sequence of pair homomorphisms, $\{h_n:(S,L) \to (FS_k,F_k) \}$,
that converges into the limit pair $(S,L)$. Then it iteratively shorten these homomorphisms and pass to subsequences  that converge into the next
limit pairs in the sequence. Hence, the terminal subsequence satisfies part (3) in theorem 6.9.

Note that to get a (shortened) subsequence that converges into the $i$-th limit group $(S_i,L_i)$, the shortening procedure that is used in the
proof of theorem 6.8 uses the modular groups that are associated with all the prefix of the sequences of pairs up to $(S_i,L_i)$, including $\Lambda_i$
that is associated with $(S_i,L_i)$.

\hfill$\square$

We call a finite collection of resolutions that satisfies properties (1) and (2) in theorem 6.8 and part (3)
in theorem 6.9 a $Makanin$-$Razborov$ $diagram$ of
the restricted pair $(S,L)$. Note that such a diagram is not canonical in general. We view sequences of pair homomorphisms that
satisfy part (3) as weak test sequences for the resolutions in the diagram.
We later use such sequences of homomorphisms as a replacement to $test$ $sequences$ that were used
in [Se2] to construct formal solutions and obtain generalized Merzlyakov theorems for AE sentences and formulas that
are defined over a given variety over a free semigroup.

\vglue 1.5pc
\chapter{A Makanin--Razborov diagram (general case)}
\fullwidthtext
\medskip

In section 8 we  analyzed the collection of homomorphisms from a restricted pair $(S,L)$ into the standard pair $(FS_k,F_k)$,
in case the pair $(S,L)$ has no restricted quotients $(\hat S, \hat L)$, in which the restricted limit group
$\hat L$ is freely decomposable (with respect to the coefficient subgroup $F_k$), i.e., in which $L$ is
freely decomposable as a restricted limit group.

Each of the resolutions in such a diagram terminates with the standard pair $(FS_k,F_k$).

In this section we generalize the construction of the Makanin-Razborov diagram to include all the possible pairs $(S,L)$.
To do that we need to generalize theorems 8.1 and 8.7 to include pairs with freely decomposable limit groups.
We first present such generalizations in case there are no Levitt components in the actions on the real trees
that we need to consider,
and then omit this assumption, and allow Levitt components. In both cases we use the machinery
that was presented in section 8, to construct  JSJ-like decompositions, that in the general case
considers and encodes Levitt components.

\vglue 1pc
\proclaim{Definition 9.1}
Let $(S,L)$ be a pair. We say that $(S,L)$ is $Levitt$-$free$ if every sequence of pair homomorphisms:
$\{h_n:(S,L) \to (FS_k,F_k)\}$ that converges into a faithful action of $L$ on some real tree,
contains no Levitt (thin) components (for a detailed description of Levitt components see [Be-Fe1] who call them
thin).
\endproclaim

\vglue 1pc
\proclaim{Theorem 9.2}
Let $(S,L)$ be a Levitt-free pair, and suppose
that the limit group $L$ contains no non-cyclic abelian subgroup.
Let $\{h_n:(S,L) \to (FS_k,F_k)\}$ be a sequence of pair homomorphisms that converges into
a faithful action of $L$ on a real tree $Y$.

Then there exists a $resolution$:
$$(S_1,L_1) \to  \ldots \to (S_f,L_f)$$
that satisfies the following properties:
\roster
\item"{(1)}" $(S_1,L_1)=(S,L)$, and  $\eta_i:(S_i,L_i) \to (S_{i+1},L_{i+1})$ is an isomorphism for $i=1,\ldots,f-2$ and
$\eta_{f-1}:(S_{f-1},L_{f-1}) \to (S_f,L_f)$ is a
quotient map.

\item"{(2)}" with each of the pairs $(S_i,L_i)$, $1 \leq i \leq f$, there is an associated abelian decomposition that we denote
$\Lambda_i$. The abelian decompositions $\Lambda_1,\ldots,\Lambda_{f-1}$ contain  edges with trivial and cyclic edge stabilizers,
QH vertex groups that are associated with IET components, and rigid vertex groups.


\item"{(3)}" either $\eta_{f-1}$ is a proper quotient map, or the abelian decomposition $\Lambda_f$ contains
a special oriented edge with trivial stabilizer that can have arbitrary values from $FS_k$. i.e., the sequence of homomorphisms that converges to
$(S_f,L_f)$ (that is obtained by iteratively shortening the original sequence of homomorphisms that converge into $(S_1,L_1)$), can be
modified by replacing the value that the edge is mapped into by arbitrary element from $FS_k$, and still remain a pair homomorphism of $(S_f,L_f)$. Such a
special oriented edge in $\Lambda_f$ can be either
separating or not.

\item"{(4)}" there exists a subsequence of the homomorphisms $\{h_n\}$ that factors through the resolution. i.e., each homomorphism
$h_{n_r}$ from the subsequence, can be written in the form:
$$h_{n_r} \ = \ \hat h_r \, \circ \, \varphi^{f-1}_r \, \circ \, \ldots \, \circ \, \varphi^1_r$$
where $\hat h_{r}:(S_f,L_f) \to (FS_k,F_k)$, each of the automorphisms $\varphi^i_r \in Mod(\Lambda_i)$, where $Mod(\Lambda_i)$
is generated by the modular groups that are associated with the QH vertex groups, and by Dehn twists along edge groups with
cyclic stabilizers in $\Lambda_i$.

Each of the homomorphisms:
$$h^i_{n_r} \ = \ \hat h_r \, \circ \, \varphi^{f-1}_r \, \circ \, \ldots \, \circ \, \varphi^i_r$$
is a pair homomorphism $h^i_{n_r}: (S_i,L_i) \to (FS_k,F_k)$.

\item"{(5)}" if $(S_f,L_f)$ is not a proper quotient of $(S,L)$, then the pair homomorphisms $\hat h_r$ are compatible with
$\Lambda_f$. Let $R_1,\ldots,R_v$ be the connected components of $\Lambda_f$ after deleting its oriented special edges.
The homomorphisms $\hat h_r$ are composed from homomorphisms of the fundamental groups of the connected
components $R_1,\ldots,R_v$,
together with assignments of values from $FS_k$ to the oriented separating edges.
The homomorphisms of the fundamental groups of the connected components $R_1,\ldots,R_v$ converge into a faithful action of
these groups on real trees with associated abelian decompositions: $R_1,\ldots,R_v$.
\endroster
\endproclaim

\nfp
We modify the  procedure that was used in proving theorem 8.1. Let $\Lambda$ be the abelian decomposition
that is associated with the action of $L$ on the limit tree $Y$ (that is obtained from the convergent
sequence $\{h_n\}$). If $\Lambda$
contains a segment in its simplicial part, and that segment has a trivial stabilizer, the conclusions
of the theorem follow.

By our assumptions, $(S,L)$ is Levitt free, so
as long as $(S,L)$ is not replaced by a proper quotient, none of the faithful actions of $L$ on the limit
trees that are constructed along the procedure contain Levitt components. Suppose that the action of $L$ on
the corresponding real tree is not geometric (in the sense of [Be-Fe]). In that case we associate with the non-geometric
action (with no Levitt components), an approximating resolution (in the sense of [Be-Fe]) according to the one that we
constructed in theorem 4.9. The graph of groups that is associated with that approximating resolution contains edges
with trivial stabilizers in their simplicial parts, and a subsequence of the given sequence of homomorphisms satisfies parts
(4) and (5) of the theorem according to theorem 4.9.

Therefore, we may assume that the action of $L$ on the limit tree $Y$ is geometric. Hence,
$Y$ contains only
 IET and discrete components,
and every segment in the discrete part of $Y$ can be divided into finitely many non-degenerate
subsegments with  (non-trivial) cyclic stabilizers. Therefore, if $L$ is a free product of non-cyclic,
freely indecomposable limit groups, and $\Lambda$ is the JSJ decomposition of $L$, $L$ is replaced by a proper quotient after
shortening along the modular group of $L$, and the conclusion of the theorem follows.

Hence, $\Lambda$ is not the JSJ of the limit group $L$.
We start by refining the abelian decomposition $\Lambda$ in a similar way to what we did in proving theorems 5.2 and 8.1.
If a non-QH vertex group is connected only to periodic edge groups,
we  restrict the homomorphisms $\{h_n\}$ to such a vertex group and obtain a non-trivial splitting of
it in which all the previous (periodic) edge groups are elliptic. Hence, the obtained abelian decomposition
of the non-QH vertex group can be used to refine the abelian decomposition $\Lambda$. Repeating this refinement
procedure iteratively, we
get an abelian decomposition that we denote $\Lambda_1$.

We fix finite generating sets of all the edge groups and all the non-QH vertex groups in $\Lambda_1$. We divide the edge and
non-QH vertex
groups in $\Lambda_1$ into
finitely many equivalence classes of their growth rates as we did the proof of theorem 8.1.
Two non-QH vertex groups or edge groups (or an edge and a non-QH vertex group)
 are said to be in the same equivalence
class if the maximal length of the images of their finite set of generators have comparable lengths, i.e., the length of a maximal
length image of one is bounded by a (global) positive constant times the maximal length of an image of other, and vice versa.
After possibly passing to a subsequence, there exists a class that dominates
all the other classes, that
we call the $dominant$ class, that includes (possibly) dominant edge groups and (possibly) dominant non-QH vertex groups.

As in the freely indecomposable case, we denote by $Mod(\Lambda_1)$ the modular group that is associated
with $\Lambda_1$. We set $MXMod(\Lambda_1)$ to be the $dominant$ modular group that is generated by Dehn twists along
dominant edge groups
and modular groups of dominant QH vertex groups. i.e., those QH vertex groups that have at least one dominant boundary component.

We start by using the full modular group $Mod(\Lambda_1)$.
For each index $n$, we set the
pair homomorphism: $h_n^1:(S,L) \to (FS_k,F_k)$,
$h_n^1=h_n \circ \varphi_n$, where $\varphi_n \in Mod(\Lambda_1)$, to be a shortest pair homomorphism
that is obtained from $h_n$ by precomposing it with a modular automorphism from $Mod(\Lambda_1)$.
If
there exists a subsequence of the homomorphisms $\{h_n^1\}$ that converges into a proper quotient of the pair
$(S,L)$, or that converges into a non-geometric action   of $L$ on a real
tree, or into an action that contains a segment in its simplicial part, and this segment can not be divided into finitely many
subsegments with non-trivial stabilizers,
we set the limit of this subsequence to be $(S_f,L_f)$, and the conclusion of the theorem follows
with a resolution of length 1.

Therefore, we may assume that every convergent subsequence of the homomorphisms $\{h_n^1\}$
converges into a faithful action of the limit group $L$ on some real tree. In that case we use only elements from the dominant
modular group, $MXMod(\Lambda_1)$.

Before we continue to the next abelian decomposition, we need to check if the original sequence of homomorphisms, $\{h_n\}$,
does not contain a subsequence of $separable$ homomorphisms. We are going to look for subsequences of separable homomorphisms in every
step of the iterative procedure for the construction of the sequence of abelian
decompositions $\Lambda_1,\ldots$. The existence of such subsequence will lead to a termination of the procedure, with an
abelian decomposition that satisfies the conclusions of theorem 9.2.

\vglue 1pc
\proclaim{Definition 9.3}
Let $(S,L)$ be a pair, and let:
$\{u_n:(S,L) \to (FS_k,F_k)\}$ be a sequence of pair homomorphisms that converges into a faithful action of $L$ on some
real tree  $Y$. We say that the sequence $\{u_n\}$ is $separable$ if there exists a (reduced)
graph of groups $\Delta$ with fundamental
group $L$ with the following properties and additional data:
\roster
\item"{(1)}" the graph of groups $\Delta$ is non-trivial and all its edges have trivial stabilizers. We further assign
orientation with each edge in $\Delta$.

\item"{(2)}" with each non-trivial vertex group in $\Delta$ we associate a base point. In addition there is a
 base point for the fundamental group $L$,
that is placed in the interior of one of the edges or it is one of the basepoints that are associated with the vertex groups.

\item"{(3)}" with each edge in $\Delta$ we associate a label. If the basepoint of $L$ is in the interior of an edge, then
with that edge there are two labels that are associated with the two parts of the edge.

\item"{(4)}" the homomorphisms $\{u_n\}$ are composed from homomorphisms of the vertex groups into $F_k$, and assignments
of values in $F_k$ to the labels that are associated with the edges in $\Delta$. Each element in $L$ can be considered as a
path in $\Delta$ that starts and ends in its basepoint. Hence, from the values that are assigned to the labels, and
the homomorphisms of the vertex groups (with their basepoints), it is possible to read (uniquely) a homomorphism of $L$
into $F_k$.

\item"{(5)}" we extend the set of generators of the standard pair $(FS_k,F_k)$, by adding a new free generator to
the standard semigroup $FS_k$ for each label in $\Delta$. We denote the extended standard pair $(FS_m, F_m)$.

Given each of the homomorphisms $u_n$, it is possible
 to replace
the values that are assigned with each of the labels that are associated with the edges in $\Delta$, to values that
contain a single (positively oriented) appearance of the generator that is associated with each label,  and no
appearances of generators that are associated with the other labels,
without changing the homomorphisms
of the vertex groups in $\Delta$ to obtain a homomorphism $\hat u_n$. The homomorphism $\hat u_n$ should coincide
with the original homomorphism $u_n$ if we map the
(new) generators that are associated with the labels to the identity. Furthermore, and for each index $n$, the homomorphism
$\hat u_n$ is a pair homomorphism: $\hat u_n: (S,L) \to (FS_m, F_m)$.
\endroster
\endproclaim

If the sequence of homomorphisms, $\{h_n\}$,
contains a separable subsequence, the conclusion of theorem 9.2 follows.

\vglue 1pc
\proclaim{Proposition 9.4}
Let $(S,L)$ be a Levitt-free pair, and suppose
that the limit group $L$ contains no non-cyclic abelian subgroup.
Let $\{h_n:(S,L) \to (FS_k,F_k)\}$ be a sequence of pair homomorphisms that converges into
a faithful action of $L$ on a real tree $Y$.

Let $\{u_t:(S,L) \to (FS_k,F_k)\}$ be a separable subsequence of the sequence $\{h_n\}$.
Then part (5) of
theorem 9.2 holds for the subsequence $\{u_t\}$.
\endproclaim

\nfp With this separable subsequence there is an
associated abelian decomposition $\Delta$ that satisfies all the properties that are listed in definition 9.3. Given
$\Delta$, and the separable sequence, $\{u_t\}$, for every index $t$, we can associate with the homomorphism $u_t$,
an extended pair homomorphism: $\hat u_t:(S,L) \to (FS_m,F_m)$ that satisfies the properties that are listed in part
(5) of definition 9.3.

If we give the generators that are associated with the labels a very high length, The sequence $\{\hat u_t\}$   converges
into a faithful action of $L$ on a real tree, $\hat Y$, with an associated abelian decomposition $\Delta$. $\Hat Y$ contains segments with trivial stabilizers in its
simplicial part.
Hence, the conclusion of theorem 9.2 holds for
the  abelian decomposition $\Delta$.

\hfill$\square$

For presentation purposes we first assume that $\Lambda_1$, and all the next abelian decompositions that are obtained
along the iterative procedure do not contain any QH vertex groups. In that case we continue iteratively, precisely as we did
in the freely indecomposable case.

First we shorten the homomorphisms $\{h_n\}$ using the dominant modular group, $MXMod(\Lambda_1)$. We denote the shortened
sequence of homomorphisms, $\{h^1_n\}$, and after passing to a subsequence, assume that the obtained sequence converges
into a faithful action on a real tree $Y_2$.
If the action of the limit group $L$ on the real tree $Y_2$ is not faithful, or if there is a segment in the simplicial part of
$Y_2$ that can not be divided into subsegments with non-trivial stabilizers, or if the action of $L$ on $Y_2$ is not geometric, the
conclusions of theorem 9.2 follow.

Hence, we may assume that the action of $L$ on $Y_2$ is geometric and faithful.
With the action of $L$ on $Y_2$ there is an associated abelian decomposition $\Delta_2$.
By our standard procedure $\Delta_2$ can be (possibly) refined to an abelian decomposition $\Lambda_2$. By our assumptions, $\Lambda_2$
contains no QH vertex groups.

We first check if the sequence of abelian decomposition $\{h^1_n\}$ contains a separable subsequence. If it does, theorem 9.2 follows.
Otherwise, we associate with $\Lambda_2$ its modular group, $Mod(\Lambda_2)$, and its dominant modular group, $MXMod(\Lambda_2)$.
We first shorten using the ambient modular group, and if the action is faithful, we shorten using the dominant modular group.

We continue iteratively, and either terminate after finite time, or obtain the sequence of abelian decompositions, $\Lambda_1,\ldots$,
and the corresponding sequences of convergent shortened homomorphisms, that have no separable subsequences. As in the freely indecomposable case,
with the sequence of abelian decompositions, $\Lambda_1,\ldots$ we associate its stable dominant (abelian) decomposition.

\vglue 1pc
\proclaim{Definition 9.5}
Let $(S,L)$ be a  pair, and let $\Lambda_1,\ldots$ be a sequence of abelian decompositions of $L$. With the sequence
we associate its stable dominant abelian decomposition, that generalizes the one that    is  presented in definition 8.2, in the
case of freely decomposable groups.

Note that given two abelian decompositions $\Delta_1$ and $\Delta_2$
of a limit group $L$, their common refinement is the multi-graded abelian JSJ decomposition of $L$
with respect to the collection of
(finitely many conjugacy classes of) subgroups that are elliptic in both splittings, $\Delta_1$ and $\Delta_2$.
The multi-graded abelian JSJ decomposition is an abelian decomposition of a freely decomposable group,
and hence includes free products, and starts with
the multi-graded Grushko decomposition. For the definition of the multi-graded JSJ decomposition see section 12 in [Se1].

Associating an abelian decomposition (a common refinement) with every pair of abelian decompositions of $L$, enables one
to associate such a common refinement with every finite sequence of abelian decompositions of $L$. Given a sequence
of abelian decompositions, $\Lambda_1,\ldots$, the (natural) complexities of the
common refinements of the prefixes, $\Lambda_1,\ldots,\Lambda_{\ell}$, does not decrease with $\ell$, and all
these decompositions are bounded by the abelian JSJ decomposition of the limit group $L$ (the abelian JSJ of a freely
decomposable group that starts with the Grushko decomposition, and then associates the abelian JSJ decomposition with
each of the non-cyclic freely indecomposable factors). Hence, the common refinements of the prefixes of the
sequence $\Lambda_1,\ldots$ stabilize, and it is possible to associate a common refinement with the entire
sequence.

Therefore, it is possible to associate an abelian decomposition (common refinement) with every suffix of the sequence,
$\Lambda_{\ell},\ldots$. The complexities of these  abelian decompositions do not increase with $\ell$,
and every sequence of strictly decreasing
abelian decompositions of $L$ terminates after a finite (in fact, bounded) time. We define the $stable$ $dominant$
abelian decomposition of the sequence $\Lambda_1,\ldots$, to be the minimal abelian decomposition that is associated
with a suffix $\Lambda_{\ell},\ldots$, where the minimum is with respect to all the suffixes of the sequence.
\endproclaim

Let $\Theta_{i_0}$ be the stable dominant abelian decomposition of the sequence $\Lambda_1,\ldots$. Suppose that
 $\Theta_{i_0}$ does not collapse to any  free
decomposition. i.e., $\Theta_{i_0}$ contains no edges with trivial stabilizers (only $QH$ and rigid vertex groups, and edges with
non-trivial cyclic stabilizers).

We assumed that the sequence of abelian decompositions $\Lambda_1,\ldots$ contains no QH vertex groups. Hence, in case
$\Theta_{i_0}$ contains no free products, the procedure that was used to prove theorem 8.1 in case the abelian
decompositions $\Lambda_1,\ldots$ contain no QH vertex groups (cf. proposition 8.3), proves the existence of a sequence
of homomorphisms that converges into an action of $L$ on a real tree, with an associated abelian decomposition, that after
a possible refinement, is identical to $\Theta_{i_0}$. Therefore, in this case the infinite sequence $\Lambda_1,\ldots$ can be replaced by a finite
sequence that ends with $\Theta_{i_0}$, precisely as we did in the freely indecomposable case.

Suppose that $\Theta_{i_0}$ collapses to free products. i.e., that there are edges with trivial stabilizers in $\Theta_{i_0}$.
 In that case, by the structure of the iterative procedure for
the construction of the abelian decompositions $\Lambda_1,\ldots$ (proposition 8.3), there exists
a sequence of shortened homomorphisms that converges into a faithful action of $L$ on a real tree. Since every element in $L$ that is hyperbolic in $\Theta_{i_0}$
has to be hyperbolic in the limit action, the abelian decomposition that
$L$ inherits from its action on the limit tree is $\Theta_{i_0}$.

If the limit action of $L$ is not geometric, the sequence of abelian decompositions:
$\Lambda_1,\ldots$, can be replaced by a finite sequence: $\Lambda_1,\ldots,\Lambda_{i_0-1},\Theta_{i_0}$, together with a subsequence of shortened homomorphisms
that converges to  the limit action that is not geometric, so this sequence of shortened homomorphisms is separable, and theorem 9.2 follows in this case.

Therefore, we may assume that the action of $L$ on the limit tree is geometric.
Since we assume that the pair  $(S,L)$ is Levitt free, this geometric action contains no Levitt components.
In this case the modular groups that are associated with all the abelian decompositions, $\Lambda_{i_0},\ldots$, are contained in
the modular group that is associated with $\Theta_{i_0}$. Therefore, the sequence of abelian decompositions,
$\Lambda_1,\ldots$ can be replaced by the finite sequence:  $\Lambda_1,\ldots,\Lambda_{i_0-1},\Theta_{i_0}$.

\medskip
At this point we still assume that the pair $(S,L)$ is Levitt-free, but allow the abelian decompositions in the sequel
to be arbitrary, i.e., to contain QH vertex groups. We have already assumed that  the sequence $\{h_n\}$, with which we start the construction,
 contains no separable
subsequence, and that
after shortening it using the ambient modular group $Mod(\Lambda_1)$, we get a sequence for which every subsequence is not separable (definition 9.3), and every
convergent subsequence converges into a faithful action of $L$.

We use only the dominant modular group, $MXMod(\Lambda_1)$.
We do what we did in the freely indecomposable case.
We shorten the action of each of the dominant QH vertex groups
using the
procedure that is used in the proof of theorem 4.2.  We shorten edge groups that are stabilized by dominant edge groups using Dehn twists.
We keep the image of the given set of the generators $s_1,\ldots,s_r$ to be in the free semigroup $FS_k$.

We further fix finite generating sets of QH vertex groups in $\Lambda_1$, and as in the freely indecomposable case,
guarantee that a fixed number of translates of the paths that are associated with each generator cover the paths of the other generators. We also
make sure that elements in a ball of radius 1 in the Cayley graph of $L$ w.r.t.\ the generating set $s_1,\ldots,s_r$ that are hyperbolic or contained in
QH vertex groups in $\Lambda_1$, have the property that finitely many translates of the paths that are associated with them cover the paths that are
associated with either the generators of a QH vertex group or the generators of a dominant edge group, if they are hyperbolic w.r.t.\  the QH or contained in the QH vertex
group or hyperbolic w.r.t.\ a dominant edge group.

We possibly pass to a subsequence, and
(still) denote the obtained (shortened) homomorphisms, using modular automorphisms from $MXMod(\Lambda_1)$, $\{h_n^1\}$.
By the shortening procedure, the lengths of the fixed generators of QH vertex group the generators:
$s_1,\ldots,s_r$,
are bounded by some  constant $c_1$ (that is independent of $n$) times the maximal length of the images of
the fixed generators of the dominant
vertex and edge groups (in $\Lambda_1$).

As in the freely indecomposable case, we pass to a subsequence of the homomorphisms $\{h_n^1\}$ that converges into an action of $L$ on some
real tree with an associated abelian decomposition $\Delta_2$. If the action of $L$ is not faithful, the conclusions of
theorem 9.2 follow. If the sequence $\{h_n^1\}$ has a separable subsequence and in particular if the action of $L$ on the limit tree is not geometric,
the conclusion of theorem 9.2 follow.

Let $\Delta_2$ be the abelian decomposition that is associated with the geometric action of $L$ on the limit tree that is obtained from a subsequence of the
shortened homomorphisms, $\{h_n^1\}$.
We further refine $\Delta_2$ as we did in the freely indecomposable case. We restrict (a convergent subsequence of)
the homomorphisms $\{h_n^1\}$ to non-QH vertex groups that are connected only
to periodic edge groups and not to QH vertex groups,  as we did with $\Lambda$, and
construct an abelian decomposition that we denote $\hat \Delta_2$.

Let $Q_1,\ldots,Q_{\ell}$ be the dominant QH vertex groups in $\Lambda_1$, so that all their boundary elements are elliptic in $\hat \Delta_2$.
i.e., every boundary element in these QH vertex groups can be conjugated into either an edge group or into a non-QH vertex group in $\hat \Delta_2$.
In that case, as in the freely
indecomposable case, we refine $\hat \Delta_2$. We start with $\Lambda_1$, and collapse all the edge groups that do not connect the QH vertex groups
$Q_1,\ldots,Q_{\ell}$ to the
 to other vertex groups in $\Lambda_1$. Every vertex group that is not one of the QH vertex groups $Q_1,\ldots,Q_{\ell}$  in the obtained
(collapsed) graph of groups inherits an abelian decomposition from $\hat \Delta_2$. We refine the abelian decomposition that was obtained from collapsing
$\Lambda_1$ using the abelian decompositions that the non-QH vertex groups inherit from $\hat \Delta_2$.  We denote the obtained abelian decomposition,
$\Lambda_2$.

We continue as in the freely indecomposable case. We shorten the homomorphisms, $\{h_n^1\}$, using the modular automorphisms
in $Mod(\Lambda_2)$. If there is a subsequence of shortened automorphisms that converges into a non-faithful action of $L$ on a real tree,
or there
is a subsequence of shortened automorphisms that is a separable subsequence, the conclusions of theorem 9.2 follow.

Otherwise we use the dominant modular group, $MXMod(\Lambda_2)$, and proceed precisely in the same way that we proceeded in the freely
indecomposable case, and in each step precisely as we did with the shortened homomorphisms $\{h_n^1\}$, in the first step of the procedure.
We denote the  abelian decomposition at step $i$, that is constructed from shortened homomorphisms
$\{h^{i-1}_n\}$, $\Lambda_i$.

If at some step $i$, $\{h_n^i\}$ has  a separable subsequence or a subsequence that converges into a
non-faithful action of $L$ on a real tree the conclusions of theorem 9.2 follow.
If for all indices $i$ none of these occur
we get an infinite sequence of
abelian decompositions, $\Lambda_1,\Lambda_2,\ldots$.
Given the infinite sequence of abelian decompositions, we associate with it its stable dominant abelian decomposition (see definition 9.5),
$\Theta_{i_0}$.

For $i \geq i_0$, the modular groups of the abelian decompositions, $\Lambda_i$, are contained in the modular group of $\Theta_{i_0}$. By the arguments
that were used in the freely indecomposable case, and the Levitt-free freely decomposable case with no QH vertex groups, from the sequences of shortened
homomorphisms, $\{h_n^i\}$, it is possible to extract a subsequence that converges into a faithful action of $L$ on a real tree $R$,
 with an abelian decomposition
that can be further refined to the abelian decomposition $\Theta_{i_0}$.

Since the modular groups of all the abelian decompositions $\Lambda_{i_0},\ldots$ is dominated by the modular groups of the abelian
decomposition $\Theta_{i_0}$, we can
replace the infinite sequence of abelian decompositions, $\Lambda_1,\ldots$, with the finite sequence, $\Lambda_1,\ldots,\Lambda_{i_0-1},\Theta_{i_0}$.
If the sequence of homomorphisms that converges to a faithful action of $L$ with an abelian decomposition that can be further refined to $\Theta_{i_0}$
has a separable subsequence, the conclusions of theorem 9.2 follow.

We proceed as we did in the freely indecomposable case. Since the abelian decomposition $\Theta_{i_0}$ was constructed from a sequence of hyperbolic-hyperbolic
or elliptic-hyperbolic splittings, $\Theta_{i_0}$ must contain a QH vertex group or it collapses to a free product of the limit group $L$.

We shorten the sequence of homomorphisms that converges to the action of $L$ on the limit tree $R$ using the modular group $Mod(\Theta_{i_0})$, and if needed,
using the dominant modular group, $MXMod(\Theta_{i_0})$, and continue precisely as we did in the construction of the sequence of abelian decompositions,
$\Lambda_1,\ldots$. If the procedure terminates after finitely many steps, theorem 9.2 follow. Otherwise we replace an infinite tail of the sequence of
decompositions by an abelian decomposition $\Theta_{i_1}$, which is the stable dominant abelian decomposition of the infinite sequence.

We continue iteratively. If the sequence terminates after a finite sequence of abelian decompositions, theorem 9.2 follows. Otherwise we get an infinite sequence
that contains an infinite subsequence of the form: $\Theta_{i_0},\Theta_{i_1},\ldots$, and each such decomposition contains either a QH vertex group or
collapses to a non-trivial free product.

By the same method that enables one to replace the infinite sequence of decompositions: $\Lambda_1,\ldots$, with the finite sequence:
$\Lambda_1,\Lambda_{i_0-1},\Theta_{i_0}$, where $\Theta_{i_0}$ is the stable dominant abelian decomposition of the sequence, it is possible to replace
the sequence: $\Lambda_1,\ldots,\Theta_{i_0},\ldots,\Theta_{i_1},\ldots$ with a finite sequence that terminates with the stable dominant abelian
decomposition of the sequence, that we denote $\hat \Theta_{\hat i_0}$.

Since the infinite sequence of abelian decompositions:
$\Lambda_1,\ldots,\Theta_{i_0},\ldots,\Theta_{i_1},\ldots$ contains hyperbolic-hyperbolic or hyperbolic-elliptic splittings, as well
as the subsequence: $\Theta_{i_0},\Theta_{i_1},\ldots$, and each of the abelian decompositions, $\Theta_{i_j}$ either collapses to a
non-trivial free product, or it contains a QH vertex group, the stable dominant abelian decomposition of the infinite sequence:
$\Lambda_1,\ldots,\Theta_{i_0},\ldots,\Theta_{i_1},\ldots$ either:

\roster
\item"{(1)}"  collapses to a free product with at least 3 factors.

\item"{(2)}" collapses to a  free
product with a free factor which is free and acts freely on the corresponding limit tree (after possible applying the refinement process that
we used to obtain each of the abelian decompositions $\Lambda_i$).

\item"{(3)}" contains a QH
vertex group with Euler characteristic that is bounded above by $-2$

\item"{(4)}" contains at least 2 QH vertex groups.

\item"{(5)}" contains a free decomposition and a QH vertex group in at least one of the factors.
\endroster

As in the proof of theorem 8.1, we repeat the whole construction starting with the (higher complexity)
abelian decomposition that we obtained and the subsequence of homomorphisms
that is associated with it.  Either the construction
terminates in finitely many steps, or  a suffix of an  infinite sequence of abelian decompositions can be replaced with
an abelian decomposition of higher complexity. i.e., properties (1)-(5) are replaced by collapsing to free products with larger number of factors,
or a free factor of higher rank that acts freely, or QH vertex groups with larger topological complexity, or a free products with more factors or more
QH vertex groups.
Repeating this construction iteratively, and applying the accessibility of f.p.\ groups for small splittings [Be-Fe1],
or by acylindrical accessibility, we are left with a finite resolution that satisfies the conclusion of theorem 9.2 (see the
proof of theorem 8.1).

\hfill$\square$

\medskip
So far we assumed that the pairs that we consider $(S,L)$ are Levitt free, and $L$ contains no non-cyclic abelian subgroups.
In section 4 we showed how to shorten sequences of homomorphisms that converge to either an axial actions of a free abelian group on a real tree
(theorem 4.1),
or to an IET action of a (closed) surface group on a real tree (theorem 4.2). In case the limiting action is an indecomposable Levitt action of
a free group on a real tree we showed how to shorten the homomorphisms in the sequence only in case there is a uniform bound on the periodicity of
the homomorphisms in the sequence (theorem 4.10). In proposition 4.14 we used the proof of the Bulitko lemma (lemma 1.3 in [Ma1])
to shorten a general pair homomorphism of a pair $(S,F)$ and obtain a homomorphism with bounded periodicity, where the bound depends only on the
pair $(S,F)$.

Given a sequence of homomorphisms that converges into an indecomposable Levitt
action of a free group, by combining proposition 4.14 with theorem 4.2 we are able to shorten a subsequence of it.
This is what we need in order to generalize theorem 9.2 to general sequences of pair
homomorphisms that converge to general limit pairs, omitting the Levitt free assumption.

To obtain this shortening,
we need the abelian decompositions that we construct and associate with an action of a limit group $L$ on a real tree to have special vertex groups
that are free factors of the ambient limit group $L$, and we call $Levitt$ vertex groups.

If a limit group $L$ admits a
superstable action on a real tree, the action is geometric, and every non-degenerate segment in the simplicial part of the tree can be
divided into finitely many segments with non-trivial abelian stabilizers.
Then using the  Rips'
machine, it is possible to associate with the action a graph of groups (see [Gu]). The vertex groups in this graph
of groups are either QH or axial (abelian) or Levitt (thin) or point stabilizers. The stabilizers of an edge in this graph
of groups are either trivial or abelian. In particular, the edges that are connected to Levitt vertex groups must have
trivial stabilizers.

\proclaim{Theorem 9.6 (cf. theorem 4.10)} Let $F$ be a free group, and let  $(S,F)$ be a pair where $S=<s_1,\ldots,s_r>$ is a f.g.\
subsemigroup that generates $F$
as a group.
Let $\{h_n\}$ be a sequence of pair homomorphisms of $(S,F)$ into $(FS_k,F_k)$ that converges into a faithful indecomposable Levitt action of the free
group $F$
on a real tree $Y$.

By proposition 4.14 there exists a constant $c_p$ that depends only on the pair $(S,F)$, such that each pair homomorphism $h_n$ can be shortened to
a pair homomorphism $\hat h_n:(S,F) \to (FS_k,F_k)$ with periodicity bounded by $c_p$, by iteratively replacing periodic subwords of the images of
the fixed generators of $S$ under the homomorphisms $h_n$ by periodic subwords
with bounded periodicity.

Then the sequence of shortened homomorphisms $\{\hat h_n\}$ contains a subsequence (that we still denote $\{\hat h_n\}$) that satisfies either:
\roster
\item"{(1)}" the subsequence is a separable subsequence.

\item"{(2)}" the subsequence converges into a quotient pair $(\hat S,L)$, and $L$ is a proper quotient of the free group $F$.

\item"{(3)}"
there exists a sequence of automorphisms $\{\varphi_{\ell} \in Aut(F)\}_{\ell=1}^{\infty}$ with the following properties:
\endroster

\roster
%
\item"{(i)}"
for each index $\ell$, there exists $n_{\ell}>0$, so that for every $n>n_{\ell}$, and all the generators of
the subsemigroup $S$, $s_1,\ldots,s_r$, $\hat h_n \circ \varphi_{\ell}(s_j) \in FS_k$.
(note that the automorphisms $\varphi_{\ell}$ do
not preserve  $\nu(S)$ in  general, but the image of $S$ under the twisted homomorphisms $\hat h_n \circ \varphi_{\ell} $ remains
in the standard positive cone $FS_k$ in $F_k$).

\item"{(ii)}" for every $n>n_{\ell}$:
$$\frac {d_T(\hat h_n \circ \varphi_{\ell}(s_j), id)}
{\max_{1 \leq j \leq r}
d_T(\hat h_n(s_j),id.)} < \frac {1}{\ell}$$
where $T$ is a  Cayley graph of the coefficient free group $F_k$ with respect to a fixed set of generators.
\endroster
\endproclaim

\nfp if the sequence $\{\hat h_n\}$ contains a separable subsequence, or a subsequence that converges into a proper limit quotient $L$ of the free
group $F$ then parts (1) or (2) hold. If it contains a uniformly bounded subsequence, then a subsequence of the uniformly bounded subsequence either converges
into a proper limit quotient of $F$ or it is separable.

Hence, we continue with  a convergent subsequence of homomorphisms (still denoted) $\{\hat h_n\}$. Since parts (1) and (2) do not hold for this subsequence, the subsequence
converges to a faithful action of $F$ on a real tree $Y_1$.

If the action of $F$ on $Y_1$ is not geometric, then by the proof of theorem
4.9 the subsequence $\{\hat h_n\}$ contains a separable subsequence. If the simplicial part of the action of $F$ on $Y_1$ contains
a non-degenerate segment with a trivial stabilizer, $\{\hat h_n\}$ contains a separable subsequence.

Therefore, the action of $F$ on $Y_1$ is geometric, and any segment in the simplicial part of $Y_1$ can be divided into finitely
many subsegments with non-trivial cyclic stabilizers. However, the homomorphisms $\{\hat h_n\}$ have bounded periodicity, so the limit tree $Y_1$ contains no subsegments with
non-trivial stabilizers. So the action of $F$ on $Y_1$ contains no simplicial part and no axial part.

By [Gu] we can associate a graph of groups $\Lambda_1$ with the action of $F$ on $Y_1$.  Since the action of $F$ on $Y_1$ is geometric with no
simplicial part and no axial part, the indecomposable
subtrees of $Y_1$ are either indecomposable IET or Levitt components, and the free group $F$ decomposes into a free product with factors that contain
the stabilizers of the IET components and factors that are the stabilizers of the Levitt components.

Now, we construct the sequence of automorphism $\{\varphi_{\ell}\}$ that shorten the homomorphisms from the sequence $\{\hat h_n\}$, and satisfy part (3) in
the theorem, by combining automorphisms for the factors in the free product that is associated with the graph of groups $\Lambda_1$. We use automorphisms of
the stabilizers of the indecomposable IET components applying theorem 4.2, and
automorphisms of the free stabilizers of  the indecomposable Levitt components applying theorem 4.10 (we can apply theorem 4.10
since the homomorphisms $\{\hat h_n\}$ have
bounded periodicity).

\hfill$\square$

The ability to shorten a sequence of pair homomorphisms that converges into an indecomposable Levitt action, enables one to  shorten  a subsequence of the
given sequence of pair homomorphism to obtain a shortened subsequence
that is either a separable sequence or to a shortened sequence that converges into a proper quotient of the free group that is associated with
the original indecomposable Levitt component.

\proclaim{Proposition 9.7}
Let $\{h_n\}$ be a sequence of pair homomorphisms of $(S,F)$ into $(FS_k,F_k)$ that converges into a faithful indecomposable Levitt action of the free
group $F$
on a real tree $Y$.

By proposition 4.14 there exists a constant $c_p$ that depends only on the pair $(S,F)$, such that each pair homomorphism $h_n$ can be shortened to
a pair homomorphism $\hat h_n:(S,F) \to (FS_k,F_k)$ with periodicity bounded by $c_p$, by iteratively replacing periodic subwords of the images of
the fixed generators of $S$ under the homomorphisms $h_n$ by periodic subwords
with bounded periodicity.

Given the sequence of pair homomorphisms $\{\hat h_n\}$,  that has uniformly bounded periodicity, we can further shorten it
free group $F$ into a residually free group and shorten
by automorphisms of the free group $F$ to obtain a subsequence of
shortened homomorphisms $\{\tilde h_{n_t}:(S,F) \to (FS_k,F_k)\}$ that is either a separable subsequence or it converges into a proper (limit) quotient of the free
group $F$.
\endproclaim

\nfp If the sequence $\{\hat h_n\}$ contains a separable subsequence, or a subsequence that converges into a proper quotient of $F$ the proposition follows.
Otherwise, for each homomorphism $\hat h_n$ we choose a pair homomorphism $\tilde h_n$, which is the shortest that is obtained from $\hat h_n$
by pre-composing it with an automorphism of the free group $F$.

If the sequence $\{\tilde h_n\}$ contains a separable subsequence or a subsequence that converges into a proper quotient of $F$ the proposition follows. If not,
then a subsequence of the homomorphisms $\{\tilde h_n\}$ satisfies part (3) in theorem 9.6, which means that for large index $n$ the homomorphisms $\tilde h_n$ can be further shortened,
and that contradicts the way they were constructed.

\hfill$\square$

To generalize theorem 9.2 to pairs that are not necessarily Levitt free, we need to allow the abelian decompositions that
are associated with pairs along a resolution to include $Levitt$ vertex groups. This means that to the generators of the modular
groups that are associated with such abelian decompositions we need to add the operation of replacing (shortening) homomorphisms by iteratively replacing
periodic subwords with periodic subwords of smaller exponent for the parts of the generators that are contained in a Levitt component, and then
pre-composing with  automorphisms  of the free groups that are associated with the Levitt components
(theorems 4.10 and proposition 9.7).

Hence, given a sequence of pair homomorphisms that converges into an action of a limit group on  a real tree,
the abelian decomposition that we associate with such an action will include stabilizers of indecomposable Levitt components, that we call
$Levitt$ vertex groups. These vertex groups are connected to the other parts of the graph of groups decomposition by edges with trivial
stabilizers.

\vglue 1pc
\proclaim{Theorem 9.8}
Let $(S,L)$ be a pair, and suppose
that the limit group $L$ contains no non-cyclic abelian subgroup.
Let $\{h_n:(S,L) \to (FS_k,F_k)\}$ be a sequence of pair homomorphisms that converges into
a faithful action of $L$ on a real tree $Y$.

Then there exists a $resolution$:
$$(S_1,L_1) \to  \ldots \to (S_f,L_f)$$
that satisfies the following properties:
\roster
\item"{(1)}" $(S_1,L_1)=(S,L)$, and  $\eta_i:(S_i,L_i) \to (S_{i+1},L_{i+1})$ is an isomorphism for $i=1,\ldots,f-2$ and
$\eta_{f-1}:(S_{f-1},L_{f-1}) \to (S_f,L_f)$ is a
quotient map.

\item"{(2)}" with each of the pairs $(S_i,L_i)$, $1 \leq i \leq f$, there is an associated abelian decomposition that we denote
$\Lambda_i$. The abelian decompositions $\Lambda_1,\ldots,\Lambda_{f-1}$ contain edges with trivial or cyclic edge stabilizers,
and QH, Levitt and rigid vertex groups.


\item"{(3)}" either $\eta_{f-1}$ is a proper quotient map, or the abelian decomposition $\Lambda_f$ contains
a special oriented edge with trivial stabilizer that can have arbitrary values from $FS_k$, as in part (3) of theorem 9.2.
The sequence of homomorphisms that converges to
$(S_f,L_f)$ (that is obtained by iteratively shortening a subsequence of the original sequence of homomorphisms that converge into $(S_1,L_1)$), can be
modified by replacing the value that the special edge is mapped into by arbitrary element from $FS_k$, and still remain a pair homomorphism of $(S_f,L_f)$. Such a
special oriented edge in $\Lambda_f$ can be either
separating or not.

\item"{(4)}" there exists a subsequence of the homomorphisms $\{h_n\}$ that factors through the resolution. i.e., for a pair homomorphism
$h_{n_r}$ from the subsequence, there exist pair homomorphisms, $h^i_{n_r}:(S_i,L_i) \to (FS_k,F_k)$, $i=1,\ldots,f$. $h_{n_r}^1=h_{n_r}$, and
$h_{n_r}^{i+1}$ is obtained from $h_{n_r}^i$ using the modular groups that are associated with
$\Lambda_i$,
i.e., modular groups that are associated with with the  QH vertex groups in $\Lambda_i$ and with the (non-trivial)
cyclic edge groups in $\Lambda_i$. For the values of the Levitt components,  we first iteratively replace periodic subwords by periodic subwords
with bounded exponent, and then shorten by pre-composing with   automorphism of the free groups that are associated with the Levitt indecomposable
components (as in theorems 4.10 and 9.7).


\item"{(5)}" Part (5) in theorem 9.2 is valid for the constructed resolution.
If $(S_f,L_f)$ is not a proper quotient of $(S,L)$, then the pair homomorphisms $h_{n_r}^f$ are compatible with
$\Lambda_f$. Let $R_1,\ldots,R_v$ be the connected components of $\Lambda_f$ after deleting its oriented special edges.
The homomorphisms $\hat hi^f_{n_r}$ are composed from homomorphisms of the fundamental groups of the connected
components $R_1,\ldots,R_v$,
together with assignments of values from $FS_k$ to the oriented separating edges.
The homomorphisms of the fundamental groups of the connected components $R_1,\ldots,R_v$ converge into a faithful action of
these groups on real trees with associated abelian decompositions: $R_1,\ldots,R_v$.
\endroster
\endproclaim

\nfp We modify the same procedure that was used in proving theorem 9.2  to include Levitt components. Recall that we
started the proof of theorem 9.2
by iteratively constructing abelian decompositions, $\Lambda_1,\ldots$. As we pointed in the construction of these
abelian decompositions in the Levitt-free case, we may assume that all the faithful actions of $L$ on real trees that we
consider are geometric, and in the simplicial parts of these actions every non-degenerate segment can be divided into
finitely many subsegments with non-trivial stabilizers.

Let $\Lambda$ be the abelian decomposition that is associated with the faithful action of $L$ on the real tree $Y$, that is
obtained from a convergent subsequence of the given sequence of homomorphisms, $\{h_n\}$.
$\Lambda$ may contain rigid, QH and Levitt vertex groups, and the stabilizers of edge groups may be trivial or cyclic. A
Levitt vertex group is adjacent only to edges with trivial stabilizers.

We start by refining the abelian decomposition $\Lambda$ in a similar way to what we did in the proof
of theorem 9.2. First we replace each non-QH, non-Levitt vertex group that is connected by non-trivial cyclic edge groups only to
non-QH non-Levitt vertex groups, by restricting the homomorphisms $\{h_n\}$ to such a vertex group, and replace the
vertex group with the obtained (non-trivial) abelian decomposition.
Repeating this refinement
procedure iteratively finitely many times, we
get an abelian decomposition that we denote $\Lambda_1$.

We fix finite generating sets of all the edge groups and all the non-QH, non-Levitt vertex groups in $\Lambda_1$. We divide the edge and
non-QH, non Levitt  vertex
groups in $\Lambda_1$ into
finitely many equivalence classes of their growth rates as we did the proof of theorems 8.1 and 9.2.
After possibly passing to a subsequence, there exists a class that dominates
all the other classes, that
we call the $dominant$ class, that includes (possibly) dominant edge groups and (possibly) dominant non-QH, non-Levitt vertex groups.

If $\Lambda_1$ contains a Levitt vertex group, then by proposition 9.7 we modify the restrictions of the homomorphisms $h_n$ to the generators of the
Levitt vertex groups by iteratively replacing highly periodic subwords with subwords of bounded periodicity. Then we shorten by pre-composing
with an automorphism of the free group that stabilizes the Levitt component.

By proposition 9.7, the obtained shortened homomorphisms contains a subsequence that is either separable or the subsequence converges into a proper
quotient of the limit group $L$. In both cases the subsequences converges into a final pair  $(S_f,L_f)$, the resolution terminates and it satisfies the conclusions of the
theorem.

Hence, if an abelian decomposition along the sequence of pairs and their associated decompositions contains a Levitt vertex group, we can proceed to a final
step of the resolution. If there is no Levitt vertex group we proceed as we did in the Levitt-free case.
We associate with the abelian decomposition a modular group,
$Mod(\Lambda_1)$, the modular group that is associated
with $\Lambda_1$. Since we can assume that $\Lambda_1$ contains no Levitt vertex groups, $Mod(\Lambda_1)$, is generated by the modular groups of the QH vertex groups,
and Dehn twists along the non-trivial cyclic edge groups in $\Lambda_1$.

With $\Lambda_1$ we also associate its $dominant$ modular group, $MXMod(\Lambda_1)$, precisely as what we did in the Levitt-free case. We shorten the given sequence of homomorphisms
 by automorphisms from $Mod(\Lambda_1)$, and
if no subsequence converges into a proper quotient of $L$, we shorten only by automorphisms from
$MXMod(\Lambda_1)$. We pass to a convergent subsequence and associate with it an abelian decomposition $\Lambda_2$.

If $\Lambda_2$ contains a Levitt vertex group, we continue to the final step by shortening the restrictions of the homomorphisms to the values of
the Levitt components, precisely as we did with $\Lambda_1$.
Otherwise we continue iteratively precisely as in the Levitt-free case (theorem 9.2). If a Levitt vertex group appears in one of the abelian decompositions that we construct
along the iterative procedure we continue to the final step of the resolution and the theorem follows. If not, the procedure continues and terminates and the theorem
follows according to the Levitt-free case (theorem 9.2).

\hfill$\square$

For presentation purposes,  in theorems 9.2 and 9.8 we assumed that the limit group $L$ contains no non-cyclic abelian groups. To
include abelian vertex groups we use the same modifications as we used in the freely indecomposable case (see definition 8.6 and
theorem 8.7).

\vglue 1pc
\proclaim{Theorem 9.9 (cf. theorems 9.2 and 9.8)}
Let $(S,L)$ be a pair,  and let $s_1,\ldots,s_r$ be a fixed generating
set of the semigroup $S$.
Let $\{h_n:(S,L) \to (FS_k,F_k)\}$ be a sequence of pair homomorphisms that converges into
a faithful action of $L$ on a real tree $Y$.

Then there exists a $resolution$:
$$(S_1,L_1) \to  (S_2,L_2) \to \ldots \to  (S_f,L_f)$$
that satisfies the following properties:
\roster
\item"{(1)}" $(S_1,L_1)=(S,L)$, and  $\eta_i:(S_i,L_i) \to (S_{i+1},L_{i+1})$ is an isomorphism for $i=1,\ldots,f-2$ and
$\eta_{f-1}:(S_{f-1},L_{f-1}) \to (S_f,L_f)$ is a
quotient map.

\item"{(2)}" with each of the pairs $(S_i,L_i)$, $1 \leq i \leq f$, there is an associated abelian decomposition that we denote
$\Lambda_i$. The abelian decompositions $\Lambda_1,\ldots,\Lambda_{f-1}$ contain edges with trivial or abelian edge stabilizers,
and QH, Levitt, abelian (axial),  and rigid vertex groups.

\item"{(3)}" either $\eta_{f-1}$ is a proper quotient map, or the abelian decomposition $\Lambda_f$ contains
a special oriented edge with trivial stabilizer that can have arbitrary values from $FS_k$, as in part (3) of theorems 9.2 and 9.8.
The sequence of homomorphisms that converges to
$(S_f,L_f)$ (that is obtained by iteratively shortening a subsequence of the original sequence of homomorphisms that converge into $(S_1,L_1)$), can be
modified by replacing the value of the special edge  by arbitrary element from $FS_k$, and still remain a pair homomorphism of $(S_f,L_f)$.

\item"{(4)}" there exists a subsequence of the homomorphisms $\{h_n\}$ that factors through the resolution. i.e., for a pair homomorphism
$h_{n_r}$ from the subsequence, there exist pair homomorphisms, $h^i_{n_r}:(S_i,L_i) \to (FS_k,F_k)$, $i=1,\ldots,f$. $h_{n_r}^1=h_{n_r}$, and
$h_{n_r}^{i+1}$ is obtained from $h_{n_r}^i$ using the modular groups that are associated with
$\Lambda_i$,
i.e., modular groups that are associated with with the  QH and abelian vertex groups in $\Lambda_i$ and with the (non-trivial)
cyclic edge groups in $\Lambda_i$. For the values of the Levitt components,  we first iteratively replace periodic subwords by periodic subwords
with bounded exponent, and then pre-compose with automorphisms of the stabilizers of the indecomposable Levitt components
(as we did in theorems 4.10 and 9.8).

\item"{(5)}" part (5) in theorems 9.2 and 9.8 holds for the subsequence of pair homomorphisms
$\{h^f_n\}$ of $(S_f,L_f)$.
\endroster
\endproclaim

Theorem 9.9 generalizes theorem 8.7 to general pairs. Given a sequence of pair homomorphisms of a pair
$(S,L)$, and a
sequence of pair homomorphisms $\{h_n:(S,L) \to (FS_k,F_k)\}$, that converges into a faithful action of $L$ on some real tree,
it proves that it is  possible to extract a subsequence that factors through
a finite resolution of the pair $(S,L)$ that terminates
with a
proper quotient of the pair $(S,L)$, or with a sequence of separable homomorphisms of $(S,L)$. This allows one to continue
the construction of a resolution iteratively and get a resolution that terminates with a graph of groups in which all the edges are
special edges (that can obtain arbitrary values) and all the vertex groups are trivial.

The existence of these resolutions, enables us to apply a compactness argument, similar to the one that we applied in the
freely indecomposable case, and associate a Makanin-Razborov diagram with
a pair $(S,L)$.

\vglue 1pc
\proclaim{Theorem 9.10} Let $(S,G)$ be a pair of a group $G$ and a f.g.\ subsemigroup $S$ that generates $G$ as a group.
Then there exist finitely many resolutions of the form:
$$(S_1,L_1) \to  (S_2,L_2) \to \ldots \to (S_f,L_f)$$
where $(S_1,L_1)$ is a quotient pair of $(S,G)$, and each $L_i$ is a limit group, such that:
\roster
\item"{(1)}" $(S_i,L_i)$ is a quotient of $(S_{i-1},L_{i-1})$.

\item"{(2)}" With each pair $(S_i,L_i)$ we associate an abelian decomposition $\Lambda_i$ with Levitt, QH, abelian and rigid vertex groups and trivial
or abelian edge groups. Some of the edges with trivial stabilizers are declared to be special edges, and these edges are carried to the
abelian decompositions in the next steps of the resolution.

\item"{(3)}" $L_f$ is a free group and its associated abelian decomposition is a graph of groups with trivial vertex and edge stabilizers.

\item"{(4)}" every  pair homomorphism, $h:(S,L) \to (FS_k,F_k)$, factors through at least one of these
finitely many resolutions.

\item"{(5)}" for each of the resolutions in the collection, there exists a sequence of homomorphisms:
$\{h_n:(S,L) \to (FS_k,F_k)\}$, that converges into a faithful action of the initial pair $(S_1,L_1)$ on
a real tree with an associated abelian decomposition (after the refinement that is used in the proof of theorem
9.2) $\Lambda_1$.

Furthermore, the sequence of homomorphisms $\{h_n\}$ can be modified using modular automorphisms and generalized
Dehn twists that are associated with QH and abelian vertex groups and with non-trivial abelian edge groups in the abelian decompositions: $\Lambda_1,\ldots,\Lambda_m$,
and by shortening the restrictions of the homomorphisms to the Levitt vertex groups in these abelian decompositions as we detail in part (4) in theorem 9.8,
and by adding new generators to the coefficient group
that are associated with the images of the labeled edge groups,
to get sequences
of pair homomorphisms $\{h_n^1\},\ldots,\{h_n^f\}$. Each of these modified sequences of homomorphisms
$\{h_n^i\}$ converges into a faithful action of the pair $(S_i,L_i)$ on a real tree with an associated
abelian decomposition (after an appropriate refinement) $\Lambda_i$.
\endroster
\endproclaim

The finite collection of resolutions that exists by theorem 9.10 is called a $Makanin$-$Razborov$ diagram of a pair $(S,G)$.
It encodes all the
pair homomorphisms from $(S,G)$ to the standard pair $(FS_k,F_k)$. Note that the diagram was constructed by a compactness argument,
so it is not canonical.

However, at least in the single ended case, in which all the subgroups along the diagram are freely indecomposable,
or at least the abelian decompositions that are associated with the pairs along the resolutions in the diagram do not contain
edges with trivial stabilizers,
it is possible to modify the construction to obtain a canonical Makanin-Razborov diagram of a pair. This is done
in [Se8], by starting with  the construction that appears in this book, associating a complexity with a diagram,
and proving that in the single ended case, a minimal complexity diagram is canonical.

We stated theorem 9.10 for pair homomorphisms. It is straightforward to generalize theorem 9.10 to restricted pair homomorphisms, or in parallel
to a diagram that encodes all the solutions to systems of equations with coefficients. The Makanin-Razborov diagram for restricted homomorphisms
(and pairs) is similar to what we constructed in the non-restricted case, just the terminal limit group $(S_f,L_f)$ is a free group with an
associated graph of groups with trivial
edge groups, and vertex groups that are all trivial besides a distinguished vertex that is stabilized by the standard pair $(FS_k,F_k)$

\end{document}